\documentclass{amsart}

\usepackage{mathrsfs}
\usepackage{amsfonts}
\usepackage{amsmath}
\usepackage{amssymb}
\usepackage{mathtools}
\usepackage{float}
\usepackage{xcolor}
\usepackage{enumitem}
\usepackage{csquotes}
\usepackage{subcaption}
\usepackage{tikz}
\usepackage{tikz-cd}
\usetikzlibrary{arrows,automata}
\usepackage{hyperref}

\usepackage{tgtermes}

\newtheorem{theorem}{Theorem}[section]
\newtheorem{lemma}[theorem]{Lemma}

\newtheorem{corollary}[theorem]{Corollary}

\newtheorem*{question}{Question}

\theoremstyle{definition}

\newtheorem*{claim}{Claim}

\newtheorem{examplex}[theorem]{Example}
\newenvironment{example}
  {\pushQED{\qed}\examplex}
  {\popQED\endexamplex}

\theoremstyle{remark}
\newtheorem{remark}[theorem]{Remark}

\numberwithin{equation}{section}

\newcommand{\pres}[3]{\textnormal{#1} \langle #2 \mid #3 \rangle}

\newcommand{\lra}[1]{\xleftrightarrow{}^\ast_{#1}}
\newcommand{\xra}[1]{\xrightarrow{}^\ast_{#1}}
\newcommand{\xr}[1]{\xrightarrow{}_{#1}}

\newcommand{\tc}[1]{\textsc{#1}}

\newcommand{\trev}{\text{rev}}

\DeclareMathOperator{\FP}{FP}
\DeclareMathOperator{\FDT}{FDT}

\begin{document}

\title{The Word Problem for One-relation Monoids: A Survey}

\author{Carl-Fredrik Nyberg-Brodda}
\address{Department of Mathematics, University of East Anglia, Norwich, England, UK}
\email{c.nyberg-brodda@uea.ac.uk}
\thanks{The author is currently a Ph.D. student at the University of East Anglia, United Kingdom}


\date{\today}

\dedicatory{Dedicated to the memory of S. I. Adian (\tc{1931}--\tc{2020}).}


\begin{abstract}
This survey is intended to provide an overview of one of the oldest and most celebrated open problems in combinatorial algebra: the word problem for one-relation monoids. We provide a history of the problem starting in \tc{1914}, and give a detailed overview of the proofs of central results, especially those due to Adian and his student Oganesian. After showing how to reduce the problem to the left cancellative case, the second half of the survey focuses on various methods for solving partial cases in this family. We finish with some modern and very recent results pertaining to this problem, including a link to the Collatz conjecture. Along the way, we emphasise and address a number of incorrect and inaccurate statements that have appeared in the literature over the years. We also fill a gap in the proof of a theorem linking special inverse monoids to one-relation monoids, and slightly strengthen the statement of this theorem. 
\end{abstract}

\vspace*{-1cm}
\maketitle

\vspace{-1cm}

\tableofcontents

\section*{Introduction}

The word problem for one-relation monoids is one of the most fundamental open problem in combinatorial algebra. The problem itself is deceptively simple to state.

\begin{question}
Is the word problem decidable for every one-relation monoid $\pres{Mon}{A}{u=v}$?
\end{question}

The fact that this problem has remained open since its conception more than a century ago is in stark contrast to the same situation in the theory of one-relator groups; among the first known results in this latter case was Magnus' \tc{1932} theorem proving that the word problem is decidable for all one-relator groups. P. S. Novikov is quoted as saying that the word problem for one-relation monoids ``contains something transcendental'', and part of the aim of this survey is to illustrate this. Although A. I. Maltsev \cite{Maltsev1965} wrote in his \tc{1965} monograph on the theory of algorithms that the problem ``has nearly been solved by Adian'', we shall see that the mysterious and complex world in which the problem lives had only just begun to unfurl at that point. This survey is intended to provide a history of the above question, and the numerous attempts to attack, simplify, and solve it. It is intended to be readable by researchers with a graduate level of experience in combinatorial algebra.

An overview of the structure of the survey is as follows. In \S1, we first present a brief rundown on some elementary concepts necessary to appreciate the statement of the question and some of terminology of the methods by which it will be attacked. Then, in \S2, an exposition of the early history of the problem and early results of decidability is presented, in which the special and cancellative cases are treated. In \S3, we then present the two types of compression, which together with a further reduction theorem can be used to prove a reduction result of the problem to some particular difficult cases. In \S4, we present Adian's algorithm $\mathfrak{A}$, which, if its behaviour could be properly understood, would solve the word problem for all one-relation monoids. We also discuss the ramifications of a result of Sarkisian's which was thought to be proved, but where a gap was later discovered. In \S5, we present some sporadic results that have appeared in various publications and contexts. Finally, in \S6, we present some modern and future directions for the problem, including links with inverse monoids, in which a gap in a proof of a theorem from \tc{2001} by Ivanov, Margolis \& Meakin is fixed, as well as links with undecidability and the Collatz conjecture.

There have been some other surveys on the word problem for one-relation monoids, which we begin by mentioning. The \tc{1984} survey by Adian \& Makanin \cite{Adian1984} deals with general algorithmic questions in algebra, and mentions some results on the one-relation case. There is also a \tc{1988} survey by Lallement \cite{Lallement1988}, which later appeared with only minor modifications in a conference proceedings \cite{Lallement1993} and as part of lecture notes \cite{Lallement1995}. However, this survey is rather brief, and does not detail much of the history or ideas behind many of the proofs; furthermore, it includes some results which are now only considered conjecture (as we shall see in \S\ref{Subsec: An incorrect proof}). Adian's \tc{1993} brief survey \cite{Adian1993} suffers from this too, giving many statements which are conditional. That survey additionally focuses primarily on the algorithm $\mathfrak{A}$. The only survey the author is aware of which addresses the now-conjectural results is the \tc{2000} survey by Adian \& Durnev \cite{Adian2000}. The scope of this survey is general decision problems, and the word problem for one-relation monoids only occupies a comparatively small part. Finally, Cain \& Maltcev \cite{Cain2009} have produced an extensive and excellent collection on the status of various miscellaneous decision problems for one-relation monoids, but it gives no details whatsoever regarding the word problem. 

The contributions of authors writing in Russian to the area are numerous. For this reason, we make a linguistic remark. For the reader unfamiliar with Russian-to-English transliteration conventions, certain names which are written in the Cyrillic alphabet can and have been transliterated to the Roman alphabet in several different ways. This is generally done phonetically, with many different standard methods of transliteration existing. As Schein notes, ``it is like transliterating the name \textit{Poincar\'{e}}, written in Russian, as \textit{Puankare} in the Roman alphabet'' \cite{Schein2015}. This means that one person can at times be split into several in the literature; this notoriously happened to the famous semigroup theorist V. V. Wagner, who himself preferred the German spelling of his name, although standard transliteration dictates that it ought to be \textit{Vagner}. Personal preference of the author is also important; for example, S. I. Adian published under \textit{Adian, Adjan}, and \textit{Adyan}, but seems to have favoured \textit{Adian} in later years. We give below, for reference, a fixed set of spellings used in the survey of author names who are affected by the above issues. The alternative transliterations of the same names can be used to inform the reader of the correct pronunciation of the names. 

\begin{figure}[h]
\begin{center}
\begin{tabular}{l|l}
\textbf{Transliteration used} & \textbf{Alternative transliteration(s)} \\
\hline
Adian & Adjan, Adyan \\
Anisimov & An\={\i}s\={\i}mov, Anisimoff \\
Greendlinger\footnotemark & Grindlinger \\
Maltsev & Malcev, Mal'cev, Mal'tsev \\
Markov & Markoff\\
Matiyasevi\v{c} & Matijasevich, Matijasevic \\
Novikov & Novikoff \\
Oganesian & Hovhannisyan, Oganessjan, Oganesyan \\
Sarkisian & Sarkisjan, Sarkisyan\\
Sushkevi\v{c} & Sushkevich, Suschkewitsch \\
Tseitin & Ceitin, Tse\u{\i}tin, Tsejtin \\
Wagner & Vagner
\end{tabular}
\end{center}
\end{figure}
\footnotetext{Martin Greendlinger was born in the US, but defected to the USSR in the \tc{1960}s, and his name was transcribed in a non-involutive way between English and Russian.}
A remark on the source material that forms the backbone of this survey is necessary. In general, most articles on the word problem for one-relation monoids are rather self-contained, and not difficult to read on their own. On the contrary, the English translations of certain Russian articles are rather poor, and at times completely change theorems as written. Remarks have been added in this survey to alert the reader of this. The author wishes to emphasise the contributions of S. I. Adian in this area of research. His results and general interest in this problem, and related areas, have been influential beyond measure. Even a cursory glance through the survey or its bibliography will make this clear. This survey does not aspire to replace its sources; the original proofs are all perfectly readable, and particular care has been taken to make precise attributions of theorems and results. However, an overview of the main ideas behind a given proof has been given at times, to aid in exposition. This has been done in part for reasons of brevity, and in part because there is little to add to the original proofs. The author hopes that the interested reader will pursue these articles and experience these excellent proofs for themselves. After all, as N. H. Abel said, one should study the masters and not the pupils.

The author wishes to express his gratitude to R. D. Gray, V. S. Guba, D. Jackson, J.~Meakin, M. V. Volkov, and the reviewer for many helpful comments. Finally, a special thanks is extended to G. Watier for his detailed and careful reading of an early version of the survey, and comments which significantly improved the exposition in numerous places. 

\clearpage

\section{Preliminaries}

In this section, we shall give some background information, including defining all terms used later and some remarks on the notation used.

\subsection{\textit{Arghmgog}}

A finite set of symbols $A$ is called an \textit{alphabet}. Then $A^\ast$ denotes the \textit{free monoid} on $A$, which consists of all words of finite length over the alphabet $A$, together with the operation of word concatenation. For example, using the charming example given by A. Turing \cite{Turing1950}, if $u, v \in A^\ast$ are two such words (e.g. \textit{arghm} and \textit{gog}), then $uv$ represents the result of writing one after the other (i.e. \textit{arghmgog}). The empty word in $A^\ast$, being the identity element, is denoted interchangeably either by $\varepsilon$ or $1$, depending on the context. For two words $u, v \in A^\ast$, the expression $u \equiv v$ indicates \textit{graphical equality}, i.e. that the words are spelled the same. The \textit{length} $|u|$ of a word $u \in A^\ast$ is defined inductively by setting $|\varepsilon| = 0$, and $|u \cdot a| = |u| + 1$ for $a \in A$ and $u \in A^\ast$. For $u \in A^\ast$ the \textit{reversal} $u^\trev$ is just $u$ written backwards. 

A \textit{rewriting system} $T$ (also called a \textit{semi-Thue system}, named after the Norwegian mathematician A. Thue \cite{Thue1914}) on an alphabet $A$ is a subset of $A^\ast~\times~A^\ast$. All rewriting systems considered, as well as the presentations associated to them (see below), will be assumed to be finite unless explicitly stated otherwise. A rewriting system $T$ induces a relation $\xrightarrow{}_T$ on $A^\ast$ as follows: if $u, v \in A^\ast$, then $u \xr{T} v$ if and only if there exist $x, y \in A^\ast$ and some rule $(\ell, r) \in T$ such that $u \equiv x\ell y$ and $v \equiv x r y$. The reflexive and transitive closure of $\xr{T}$ is denoted $\xra{T}$. The symmetric and transitive closure of $\xra{T}$ is denoted $\lra{T}$. If $(\ell, r) \in T$, then replacing an occurrence of $\ell$ by $r$ (or vice versa) in some word $u \in A^\ast$ is called an \textit{elementary transformation} in $T$.

A rewriting system $T$ on $A$ is called \textit{terminating} (also sometimes called \textit{Noetherian}) if there exists no infinite chain $u_1 \xr{T} u_2 \xr{T} \cdots$. The system is called \textit{locally confluent} if for all $u, v, w \in A^\ast$, we have $u \xr{T} v$ and $u \xr{T} w$ together imply that there exists some $z \in A^\ast$ such that $v \xra{T} z$ and $w \xra{T} z$. The system is called \textit{confluent} if for all $u, v, w \in A^\ast$, we have $u \xra{T} v$ and $u \xra{T} w$ together imply that there exists some $z \in A^\ast$ such that $v \xra{T} z$ and $w \xra{T} z$. If a rewriting system $T$ is terminating and confluent, then we say that $T$ is \textit{complete} (also sometimes called \textit{convergent}). A word is $w \in A^\ast$ is \textit{irreducible} (modulo $T$) if it does not contain any subword that is a left-hand side of some rule of $T$. The set of irreducible elements of $T$ is denoted $\operatorname{Irr}(T)$. If $T$ is terminating, we can for every word $w \in A^\ast$ find an element $w' \in \operatorname{Irr}(T)$ such that $w \xra{T} w'$ by ``rewriting'' $w$, i.e. continuously removing any left-hand sides of rules we find as subwords of $w$ until this cannot be done any further. In a complete rewriting system, there exists a unique such $w'$. Hence any complete rewriting system has unique normal forms for all elements. The \textit{word problem} for a rewriting system $T$ over $A$ is the decision problem of, given any two words $u, v \in A^\ast$, deciding in finite time whether $u \lra{T} v$. Naturally, for a finite complete rewriting system, this problem can be solved by computing the normal forms of the input words, followed by graphical comparison.\footnote{However, the time complexity of the word problem even for finite complete rewriting systems can be arbitrarily difficult, though decidable; i.e. for every Grzegorczyk complexity class $\mathscr{C}$ there is a finite complete rewriting system for which the word problem lies in $\mathscr{C}$, cf. \cite{Bauer1984}.}

A \textit{monoid presentation} $\pres{Mon}{A}{u_i = v_i \: (1 \leq i \leq p)} = \pres{Mon}{A}{T}$ is the ordered pair $(A, T)$, where $T$ is a rewriting system on $A$. We will abuse notation and take such a monoid presentation to denote the monoid $A^\ast / \lra{T}$. This quotient, called the \textit{monoid defined by the presentation $\pres{Mon}{A}{T}$} is well-defined, as $\lra{T}$ is clearly the smallest congruence containing $T$. We will abuse notation and substitute ``let $M = \pres{Mon}{A}{u_i = v_i \: (1 \leq i \leq p)}$'' for ``let $M$ be the monoid defined by the presentation $\pres{Mon}{A}{u_i = v_i \: (1 \leq i \leq p)}$''. We say that $M$ is \textit{finitely presented}. If two words $u, v \in A^\ast$ are equal in $A^\ast / \lra{T}$, then we say that $u=v$ in $M$. If $u, v \in A^\ast$ are such that $u$ can be obtained from $v$ by an elementary transformation in $T$, then we shall say that $u$ can be obtained from $v$ by an elementary transformation in $M = \pres{Mon}{A}{T}$.

We shall also speak of a \textit{group presentation} $\pres{Gp}{A}{T}$, which is a shorthand for the monoid presentation $\pres{Mon}{A \cup A^{-1}}{T \cup \{ a_ia_i^{-1} = 1, a_i^{-1}a_i = 1 \mid a_i \in A \}}$, where $A^{-1}$ is a set in bijective correspondence with $A$ such that $A \cap A^{-1} = \varnothing$, and $T$ is a subset of $(A \cup A^{-1})^\ast \times (A \cup A^{-1})^\ast$.

Unless explicitly specified, all (monoid or group) presentations in this survey will be assumed to be finite.

Let $M = \pres{Mon}{A}{T}$. We say that a word $u \in A^\ast$ is \textit{right invertible} in $M$ if there exists some $v \in A^\ast$ such that $uv = 1$ in $M$. We define \textit{left invertibility} analogously. We say that $u \in A^\ast$ is \textit{invertible} in $M$ if it is left and right invertible. We say that $u$ is \textit{right divisible} by $v$ if there exists $w \in A^\ast$ such that $u = wv$ in $M$, and \textit{left divisibility} is defined analogously. We say that $M$ is \textit{right cancellative} if for all $u, v, x \in A^\ast$ we have that $ux = vx$ in $M$ implies $u = v$ in $M$. We define \textit{left cancellative} analogously. We say that $M$ is \textit{cancellative} if it is left and right cancellative. Note that every group is cancellative. 

\subsection{Notational remarks}

We make some remarks on the notation used in this survey in contrast to other notation for the same concepts found elsewhere in the literature. As mentioned, we use $\equiv$, to denote equality of words, i.e. equality in the free monoid. This is sometimes denoted $\eqcirc$ in older articles, particularly Soviet ones. We sometimes use $:=$ to denote a ``definitional'' equality, i.e. that the equality in question is also a definition. This is sometimes denoted $\rightleftharpoons$ in older articles, particularly Soviet ones. We use $|u|$ to denote the length of a word in the free monoid. This is sometimes denoted $[u^\partial$ or $\partial(u)$ in older articles, particularly Soviet ones, where $\partial$ is used to represent the first letter in the Russian word \textit{dlina}, meaning \textit{length}. We denote the \textit{empty word} as $1$ or $\varepsilon$, depending on notational convenience. This is sometimes denoted $\Lambda$ in older articles, as well as articles in theoretical computer science and set theory, with $\Lambda$ indicating the first letter of the German \textit{leer}, meaning \textit{empty}. 

\subsection{Decision problems}

We give some examples of decision problems which are of central importance to this survey. Let $M = \pres{Mon}{A}{R}$. Then the \textit{word problem} for $M$ has as input two words $u,v \in A^\ast$, and outputs \tc{yes} if $u = v$ in $M$, and otherwise outputs \tc{no}. The \textit{left divisibility problem} for $M$ has as input two words $u, v \in A^\ast$ and outputs \tc{yes} if $u$ is left divisible by $v$ in $M$, and otherwise outputs \tc{no}. The \textit{right divisibility problem} is defined entirely analogously. In general, these three problems are pairwise independent from one another; indeed, the divisibility problems are trivially solvable whenever $M$ is a group. However, if $M$ is given by a presentation in which all defining relations are non-empty, and such that $M$ is left cancellative, then it is not hard to show that decidability of the left divisibility problem implies decidability of the word problem, by induction on word length and noting that no non-empty word is equal to the empty word in such a monoid. The analogous statement is true substituting \textit{right} for \textit{left}.

We note that as, in the context of this survey, a monoid $M$ is always assumed to be given by a finite presentation, say $\pres{Mon}{A}{R}$, we have that if two words $u, v \in A^\ast$ are equal in $M$, then we can find a sequence of single applications of relations from $R$ which transforms $u$ into $v$. Thus, if the left (right) divisibility problem is decidable for $M$, and one finds that the word $u$ is left divisible by the word $v$, then one can always effectively construct a ``witness'' word $w$ such that $u = vw$ in $M$ (resp. $u = wv$ in $M$). 

We also introduce the following useful piece of notation. Let $M = \pres{Mon}{A}{R}$. Let $R^\trev~\subseteq~A^\ast \times A^\ast$ be the rewriting system with rules $(u^\trev, v^\trev)$ whenever $(u, v) \in R$. Let 
\[
M^\trev = \pres{Mon}{A}{R^\trev}.
\]
Then it is easy to see that the word problem for $M$ reduces to the word problem for $M^\trev$, and vice versa; for $u = v$ in $M$ if and only if $u^\trev = v^\trev$ in $M^\trev$. More importantly, the \textit{left} divisibility problem for $M$ reduces to the \textit{right} divisibility problem for $M^\trev$, for given $u, v \in A^\ast$, there exists a word $w \in A^\ast$ such that $u = wv$ in $M$ if and only if $u^\trev = v^\trev w^\trev$ in $M^\trev$. This trick will often be used. 

In \tc{1947}, and almost simultaneously, Markov \cite{Markov1947, Markov1947a} and Post \cite{Post1947} proved the existence of a finitely presented monoid with undecidable word problem. This was quite a remarkable theorem, and can be, as noted by Crvenkovi\v{c} \cite{Crvenkovic1995}, considered the first undecidability result outside the foundations of mathematics. Providing monoids with an undecidable word problem is also no mere idle pursuit if one is interested in providing \textit{groups} with an undecidable word problem, which was at times seen as a primary motivation. Indeed, A. Turing's famous proof \cite{Turing1950} of the existence of a finitely presented cancellative monoid with undecidable word problem plays a key r\^ole in Novikov's \tc{1955} detailed proof of the existence of a finitely presented group with undecidable word problem.\footnote{Novikov's original construction did not use Turing's construction, but upon writing down the detailed proof he realised the proof could be simplified in this way \cite{Adian1984}.} Turing's proof is at times rather inaccurate, and should best be read with the accompanying \tc{1958} analysis of these issues by Boone \cite{Boone1958}. We note that because of the large number of alterations required to make Turing's proof correct, Adian \& Novikov \cite{Adian1958} gave an argument in \tc{1958} which modifies Novikov's original argument to circumvent any reference to cancellative monoids with an undecidable word problem. 

A contrasting theorem had been known for decades. This was a \textit{decidability} result by W. Magnus \cite{Magnus1932}, a student of M. Dehn's, who translated Dehn's geometric intuition about the structure of one-relator groups into a purely combinatorial result \cite{Chandler1982}.

\begin{theorem}[Magnus, \tc{1932}]
The word problem is decidable for every one-relator group.
\end{theorem}

Here a \textit{one-relator} group is one that can be defined by a group presentation with a single defining relation $\pres{Gp}{A}{w=1}$. By contrast, the best known \textit{undecidability} result for groups is, to this day, a $12$-relator group with undecidable word problem due to Borisov \cite{Borisov1969}, and the word problem for $k$-relator groups when $2 \leq k \leq 11$ remains open in general. On the other hand, a much smaller gap is known for monoids. In the sequel to his \tc{1947} paper, Markov provided an example of a monoid with $33$ defining relations and undecidable word problem \cite{Markov1947a}. This was subsequently improved, in \tc{1956} and \tc{1958}, respectively, by D. Scott and G. S. Tseitin, who both provided examples of monoids with seven very short defining relations and undecidable word problem \cite{Scott1956, Tseitin1958}. Tseitin's example remains the shortest, with respect to total length of the defining relations, known example. 

The number of defining relations sufficient for presenting a monoid with undecidable word problem continued to creep down. In \tc{1966}, G. S. Makanin provided an example showing that five (short!) defining relations suffice \cite{Makanin1966}; Ju. V. Matiyasevi\v{c} \cite{Matijasevic1967b} provided an example with the same number of relations (one of which is rather long) in \tc{1967}. This record would not last for long; that same year, Matiyasevi\v{c} improved this to give an example of a monoid with only three defining relations and with undecidable word problem \cite{Matijasevic1967}.\footnote{Adian \cite{Adian2018} recalls that at the end of a \tc{1966} seminar in Moscow given by Makanin regarding his five-relation example, A. A. Markov conjectured that the number of relations could be reduced to three, and suggested Makanin write to Matiyasevi\v{c}. Apparently, Matiyasevi\v{c} had already found such an example, as it was published the next year.} The first two relations of this monoid are very short; the third is very long (several hundred letters in either word). Three relations remains the smallest number of defining relations known to suffice to present a monoid with undecidable word problem. In this way, we have at this point arrived at the question at the heart of this survey.

\begin{question}
Is the word problem decidable for every one-relation monoid $\pres{Mon}{A}{u=v}$?
\end{question}

The deceptively simple nature of the question is a large part of what makes the word problem for one-relation monoids such a fascinating problem; the fact that it remains open even today makes it all the more intriguing.

\clearpage

\section{Early results (\tc{1914}--\tc{1960})}

When one is presented with a one-relation monoid $\pres{Mon}{A}{u=v}$ and faced with the task of solving its word problem, one of the first questions one ought to ask is: what are the trivial cases? We shall begin with these, and then present the theory of cancellative and special one-relation monoids, two classes which were quickly dealt with. The ``father of semigroup theory'', A. K. Sushkevi\v{c}, considered a problem about semigroups -- or indeed monoids -- solved if it could be reduced to a problem about groups \cite[\S38]{Sushkevich1922} (see also \cite{Gluskin1972}). We shall see that this theme is very much present in these early results, in which reductions to Magnus' result on the word problem for one-relator groups are made. 

\subsection{Equal length and self-overlap free words}\label{Subsec: Easy cases}

The first observation one might make when presented with a one-relation monoid $\pres{Mon}{A}{u=v}$ is that if $|u| = |v|$, then any elementary transformation of any word will keep its length fixed. In particular, two words are equal only if their lengths are equal; hence one can effectively enumerate all finitely many words equal to a given word, giving an immediate solution to the word problem.\footnote{While the obvious algorithm produces an exponential time solution, a detailed analysis due to M\'{e}tivier \cite{Metivier1985} shows that this word problem can in fact be solved in polynomial time.} This observation was already made by Thue in \tc{1914}, in the very paragraph following his introduction of the word problem for monoids \cite[Problem~I]{Thue1914}.

In fact, that same paragraph by Thue provides another trivial case. Suppose $|u| > |v|$, and that $u$ is \textit{self-overlap free}, i.e. no non-trivial prefix of $u$ is also a suffix of $u$ (such a word is often called a \textit{hypersimple} word in the Soviet literature). Then the rewriting system with the single rule $u \to v$ is locally confluent, as there are no critical pairs; furthermore, it is terminating as $|u| > |v|$. By Newman's lemma \cite{Newman1942}, it is a finite complete rewriting system for the monoid, thus solving the word problem. This is essentially the idea behind Thue's proof, although this obviously has no reference to the \tc{1942} Newman's lemma.

We pause at this moment to address a potentially misleading comment which has appeared in the literature. In \tc{1984}, Book \& Squier \cite{Book1984} proved that ``almost all'' one-relation monoids have decidable word problem. This claim is made specific in the following sense: for a positive integer $k$, fix an alphabet $A$ of size $k$. For a positive integer $n>1$ let $u_k(n)$ be the number of self-overlap free words in $A^\ast$. Then Book \& Squier proved that the ratio $u_k(n) / k^n$ tends to $1$ as $k, n \to \infty$ (though note that this result was already observed by Nielsen \cite{Nielsen1973}, outside the context of one-relation monoids). A quick thought together with this result yields that ``almost all'' one-relation monoids have decidable word problem.

This result is less exciting than it first appears, as the following analysis will indicate. Indeed, for fixed $k$, the ratio $u_k(n) / k^n$ does not tend to $1$ as $n \to \infty$. Indeed, for $k=2$, the ratio is approximately $0.2677868$ (see OEIS sequence $A094536$ \cite{OeisBifix}, and also Nielsen \cite{Nielsen1973} for other values of $k$). Hence this argument can only be used to yield that around $27\%$ of two-generated one-relation monoids $\pres{Mon}{a,b}{u=v}$ have decidable word problem, which is not quite as exciting -- by comparison, the fraction of two-generated one-relation monoids $\pres{Mon}{a,b}{u=v}$ satisfying $|u|=|v|$ is quickly seen, by summing a geometric series, to be $\frac{1}{3}$. The asymptotic argument by Book \& Squier (as beautiful as the statement might be) also ignores the wild and complex behaviour of one-relation monoids that we shall detail presently. It is hence not particularly useful as a tool for gaining insight into the word problem for \textit{all} one-relation monoids. In spite of this, Book \cite{Book1987} frames the above density result as strong evidence that the word problem is decidable for all one-relation monoids; we reiterate that, in view of the above analysis, this framing is not accurate.\footnote{Adian \cite[p.294]{Adian1993} has made a brief remark to the same effect.}

Outside these trivial cases\footnote{Thue, with remarkable foresight, also provides some examples (see \cite[\S VIII]{Thue1914}) of other solvable word problems using what is clearly recognisable as a prototypical form of the Knuth-Bendix completion algorithm, half a century before this would be defined.}, the word problem for one-relation monoids appears to have laid untouched for some decades, with seemingly little interest in the problem. However, in the early \tc{1960}s, the Soviet mathematician S. I. Adian would begin working on this problem, and would thus begin the development that would transform the area into its modern form.

\subsection{Cancellativity and embeddability}\label{Subsec: Cancellative}

The study of cancellativity of monoids and their embeddability into groups goes back to the very beginning of semigroup theory. Any submonoid of a group is cancellative, and in the commutative case it is not hard to see (analogous to constructing a field of fractions) that a cancellative monoid can be embedded in a group. Furthermore, it is clear by universal considerations that if $M = \pres{Mon}{A}{R}$ is group-embeddable, then $M$ can be embedded in $\pres{Gp}{A}{R}$, i.e. the group with the ``same presentation'' as $M$, by the identity map $a \mapsto a$. For example, the cancellative commutative monoid $\pres{Mon}{a,b}{ab=ba}$, isomorphic with $\mathbb{N} \times \mathbb{N}$, can be embedded by the identity map in $\pres{Gp}{a,b}{ab=ba}$, isomorphic with $\mathbb{Z} \times \mathbb{Z}$.

Sushkevi\v{c} studied the problem of embedding cancellative monoids in groups, and in \tc{1935}, he published a ``proof'' that being cancellative is also sufficient for embeddability into a group \cite{Sushkevich1935}!\footnote{The author thanks Christopher Hollings for providing him with a copy of this paper.} This ``proof'', however, would not be long-lived; in \tc{1937}, Maltsev found a counterexample to Sushkevi\v{c}'s ``theorem'', i.e. an example of a cancellative monoid which is not group-embeddable \cite{Maltsev1937}. Sushkevi\v{c} later that same year wrote a monograph\footnote{Very few physical copies of this monograph remain, most having been destroyed during the many battles in the city of Kharkiv, Ukraine, during World War II (see \cite{Hollings2009}). The author of the present survey is in possession of one of these physical copies, and is currently producing an English translation of the monograph \cite{NybergBrodda2021d}.} on the theory of generalised groups, in which he (unsuccessfully) attempts to fix his erroneous proof, while simultaneously, slightly perplexingly, acknowledging Maltsev's counterexample \cite{Sushkevich1937}. Maltsev was correct, and would later produce a countable list of necessary and sufficient conditions for a monoid to embed in a group, such that no finite sublist is also necessary and sufficient \cite{Maltsev1939, Maltsev1940}. Later, Adian would provide a monoid which is finitely presented as a cancellative monoid, but which is not finitely presented as a monoid \cite[Theorem~1]{Adian1957}. Hollings has written an excellent and thorough overview of the history of embedding monoids into groups, to which we refer the interested reader \cite{Hollings2014}.

As concluded above, being cancellative is not in general sufficient for a monoid to embed in a group. One might instead ask what conditions \textit{are} sufficient. S. I. Adian seems to be have become interested in this problem -- and cancellativity in general -- at an early stage; indeed, in \tc{1955}, in one of his first published papers, he proved the existence of a finitely presented cancellative semigroup with undecidable divisibility problems \cite[Theorem~1]{Adian1955}. Five years later, Adian \cite{Adian1960b} introduced a very simple criterion which is sufficient for group-embeddability, and which furthermore can easily be read off a presentation for the monoid. We now present this criterion.

Let $M = \pres{Mon}{A}{R}$. The \textit{left graph} $\mathcal{L}(M)$ of $M$ is defined as the undirected (not necessarily simple) graph with vertex set $A$, and an edge $(a_i, a_j)$ for $a_i, a_j \in A$ in $\mathcal{L}(M)$ for \textit{every} occurrence of a relation $r=s$ in $R$ in which $a_i$ is the initial letter of $r$, and $a_j$ is the initial letter of $s$. We define the \textit{right graph} $\mathcal{R}(M)$ in an analogous way, substituting terminal letters for initial letters. We emphasise that we permit both loops and multiple edges; see the examples below. We say that (the presentation for) $M$ is \textit{left cycle-free} if $\mathcal{L}(M)$ is a forest, i.e. a disjoint union of trees, and otherwise we say that $M$ has \textit{left cycles}.\footnote{The property \textit{left cycle-free} has at times been translated from Russian to English as either \textit{left non-cancellable} or \textit{irreducible from the left}; however, in poor translations, \textit{left cancellative} has sometimes become \textit{reducible from the left}, which would be the opposite meaning, all combining to make for rather confusing reading. Context, however, always makes such statements discernible and fixable with little difficulty.} We define \textit{right cycle-free} analogously, and we say that $M$ is \textit{cycle-free} if it has no left or right cycles. We extend this definition to group presentations in which the defining relations are all written over a positive alphabet. Thus we may speak of e.g. the cycle-free group $\pres{Gp}{a,b}{ab=ba}$.

\begin{example}
Let $M_1 = \pres{Mon}{a,b,c}{ab = ba^2, ac = c^2b}$. Then the left and right graphs of $M_1$ are given below.
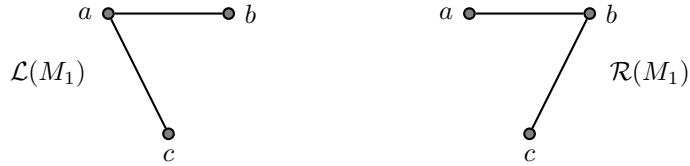
\[
\begin{tikzpicture}[>=stealth',thick,scale=0.8,el/.style = {inner sep=2pt, align=left, sloped}]%

\node (l0)[label=below:$c$] [circle, draw, fill=black!50,
                        inner sep=0pt, minimum width=4pt] at (-3,0) {};
\node (l1)[label=right:$b$] [circle, draw, fill=black!50,
                        inner sep=0pt, minimum width=4pt] at (-2,2) {};
\node (l2)[label=left:$a$] [circle, draw, fill=black!50,
                        inner sep=0pt, minimum width=4pt] at (-4,2) {};
\path[-] 
    (l2)  edge (l1)
    (l0)  edge (l2);
    
    \node (r3)[label=below:$c$][circle, draw, fill=black!50,
    inner sep=0pt, minimum width=4pt] at (3,0) {};
    \node (r4)[label=left:$a$] [circle, draw, fill=black!50,
    inner sep=0pt, minimum width=4pt] at (2,2) {};
    \node (r5)[label=right:$b$] [circle, draw, fill=black!50,
    inner sep=0pt, minimum width=4pt] at (4,2) {};

\path[-] 
    (r4)  edge  (r5)
    (r5)  edge  (r3);

\node at (-5,1) {$\mathcal{L}(M_1)$};
\node at (5,1) {$\mathcal{R}(M_1)$};

\end{tikzpicture}
\]
Hence $M_1$ is both left and right cycle-free. In particular (see Theorem~\ref{Thm: M cycle-free embeds in group} below), $M_1$ is cancellative and embeds in the cycle-free group $G_1 = \pres{Gp}{a, b, c}{ab = ba^2, ac=c^2b}$. By a simple Tietze transformation, this latter group is a one-relator group isomorphic with an HNN-extension of a free group of rank three, and the word problem is hence straightforward to solve using the Britton-Novikov lemma and the Nielsen procedure for decidability of the membership problem for subgroups of free groups; alternatively, one can use Magnus' breakdown procedure (cf. e.g. \cite{McCool1973}). In either case, having solved the word problem in $G_1$ we hence conclude that $M_1$ has decidable word problem, as $M_1 \leq G_1$. 
\end{example}

\begin{example}
Let $M_2 = \pres{Mon}{a,b,c}{b^2ca  =   bca b, a^2b = ba^2 c, ac^2 = cb^2a}$. Then the left and right graphs of $M_2$ are given below.  

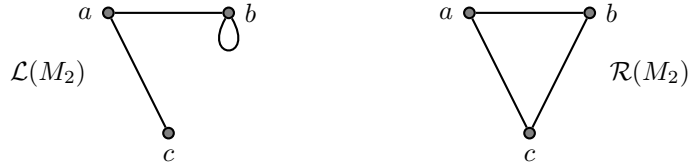
\[
\begin{tikzpicture}[>=stealth',thick,scale=0.8,el/.style = {inner sep=2pt, align=left, sloped}]%

\node (l0)[label=below:$c$] [circle, draw, fill=black!50,
                        inner sep=0pt, minimum width=4pt] at (-3,0) {};
\node (l1)[label=right:$b$] [circle, draw, fill=black!50,
                        inner sep=0pt, minimum width=4pt] at (-2,2) {};
\node (l2)[label=left:$a$] [circle, draw, fill=black!50,
                        inner sep=0pt, minimum width=4pt] at (-4,2) {};

\path[-] (l1) edge [loop below,out=300,in=240, looseness=20] (l1);
\path[-] 
    (l2)  edge (l1)
    (l0)  edge (l2);
    
    \node (r3)[label=below:$c$][circle, draw, fill=black!50,
    inner sep=0pt, minimum width=4pt] at (3,0) {};
    \node (r4)[label=left:$a$] [circle, draw, fill=black!50,
    inner sep=0pt, minimum width=4pt] at (2,2) {};
    \node (r5)[label=right:$b$] [circle, draw, fill=black!50,
    inner sep=0pt, minimum width=4pt] at (4,2) {};

\path[-] 
    (r3)  edge (r4)
    (r4)  edge  (r5)
    (r5)  edge  (r3);

\node at (-5,1) {$\mathcal{L}(M_2)$};
\node at (5,1) {$\mathcal{R}(M_2)$};

\end{tikzpicture}
\]
Hence $M_2$ has both left and right cycles. We can conclude nothing about the cancellativity of $M_2$ based on these graphs, nor anything about its group-embeddability. Solving the word problem in this monoid is left as a potentially somewhat interesting challenge. 
\end{example}

\begin{example}\label{Ex: double edges cycles}
Let $M_3 = \pres{Mon}{a,b,c,d}{ab = cd, aeb = ced}$. Then the left and right graphs of $M_3$ are
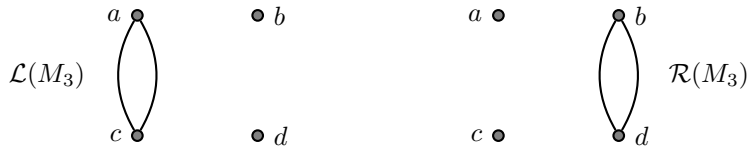
\[
\begin{tikzpicture}[>=stealth',thick,scale=0.8,el/.style = {inner sep=2pt, align=left, sloped}]%

\node (l0)[label=left:$c$] [circle, draw, fill=black!50,
                        inner sep=0pt, minimum width=4pt] at (-4,-0) {};
\node (l4)[label=right:$d$] [circle, draw, fill=black!50,
                        inner sep=0pt, minimum width=4pt] at (-2,0) {};
\node (l1)[label=right:$b$] [circle, draw, fill=black!50,
                        inner sep=0pt, minimum width=4pt] at (-2,2) {};
\node (l2)[label=left:$a$] [circle, draw, fill=black!50,
                        inner sep=0pt, minimum width=4pt] at (-4,2) {};

\path[-] 
    (l0)  edge[bend left] (l2)
    (l0)  edge[bend right] (l2);

    \node (r3)[label=left:$c$][circle, draw, fill=black!50,
    inner sep=0pt, minimum width=4pt] at (2,0) {};
    \node (r2)[label=right:$d$][circle, draw, fill=black!50,
    inner sep=0pt, minimum width=4pt] at (4,0) {};
    \node (r4)[label=left:$a$] [circle, draw, fill=black!50,
    inner sep=0pt, minimum width=4pt] at (2,2) {};
    \node (r5)[label=right:$b$] [circle, draw, fill=black!50,
    inner sep=0pt, minimum width=4pt] at (4,2) {};

\path[-] 
    (r2)  edge[bend left] (r5)
    (r2)  edge[bend right] (r5);

\node at (-5.5,1) {$\mathcal{L}(M_3)$};
\node at (5.5,1) {$\mathcal{R}(M_3)$};

\end{tikzpicture}
\]
Hence $M_3$ has both left and right cycles. We can conclude nothing about its cancellativity based on these graphs; however, it can be shown that it is simultaneously both cancellative and not group-embeddable (see below).
\end{example}

\begin{example}
Let $\Gamma$ be an (undirected) finite graph with vertex set $a_1, \dots, a_k$. Let 
\[
A(\Gamma) = \pres{Gp}{a_1, \dots, a_k}{a_i a_j = a_j a_i \: \textnormal{whenever} \: (a_i, a_j) \in E(\Gamma) }.
\]
Then $A(\Gamma)$ is called a \textit{right-angled Artin group} (RAAG). Evidently, the left and right graphs of $A(\Gamma)$ are both isomorphic to $\Gamma$, and hence the above presentation for $A(\Gamma)$ is a cycle-free presentation if and only if $\Gamma$ is a finite forest. Hence, given any $\Gamma$, as the left and right divisibility problems are trivially decidable in the monoid presentation with the same generators and defining relations (often called a \textit{trace monoid}), it follows by a result due to Sarkisian~\cite[Theorem~3]{Sarkisian1979} that the word problem is decidable for the RAAG $A(\Gamma)$ whenever $\Gamma$ is a finite forest. We remark that this is quite a contorted method of solving the word problem; in fact the word problem is relatively straightforward to solve in $A(\Gamma)$ for arbitrary finite graphs $\Gamma$ (see especially Crisp et al. \cite{Crisp2009}). Furthermore, it is a consequence of a more general due to Paris \cite{Paris2002} that any trace monoid embeds in its corresponding right-angled Artin group; this is also easy to prove using the theory of rewriting systems, cf. Chouraqui \cite{Chouraqui2009}. Right-angled Artin groups play a key r\^ole in modern geometric group theory due to their rich subgroup structure, as demonstrated in the work of Wise \& Haglund on special cube complexes (an area far beyond the scope of this survey; the reader is directed to Wise's monograph \cite{Wise2012}) and the recent interest in solving equations over right-angled Artin groups, cf. e.g. \cite{Shestakov2005, Shestakov2006, CasalsRuiz2011}.
\end{example}

We emphasise that, given a presentation, it is easily decidable whether it has left or right cycles, or indeed whether it is cycle-free. The key property of cycle-free monoids is that one can show that any such monoids are in fact group-embeddable. 

\begin{theorem}[{\cite[Theorem~5]{Adian1960b}}]\label{Thm: M cycle-free embeds in group}
Let $M = \pres{Mon}{A}{R}$ be a finitely presented monoid. If $M$ is left (right) cycle-free, then $M$ is left (right) cancellative. If $M$ is cycle-free, then $M$ is cancellative and can furthermore be embedded in $\pres{Gp}{A}{R}$ via the identity map.
\end{theorem}

We make some remarks on this theorem before proceeding, and some developments following the work by Adian. The theorem has been generalised to all cycle-free monoids (not just finitely presented), by Remmers \cite{Remmers1980}, who used used the \textit{diagram} method of geometric semigroup theory. See Higgins \cite[1.\S 7 and 5.\S 3]{Higgins1992} for an excellent overview of these methods. Recently, these diagram methods have been used to generalise Adian's theorem to certain situations in which some relations of the form $w=1$ (see \S\ref{Subsec: Special monoids}) are also permitted \cite{Kashkarev2013}. Diagram methods have also been used for studying left or right cycle-free monoids from the point of view of \textit{asphericity} e.g. in the work by Kilibarda \cite{Kilibarda1997} (see also the monograph by Guba \& Sapir \cite{Guba1997b}). A. I. Valitskas also strengthened the second half of Theorem~\ref{Thm: M cycle-free embeds in group} to prove that if a monoid $M = \pres{Mon}{A}{R}$ is (1) left (right) cycle-free; and (2) right (left) cancellative; then $M$ is group-embeddable. This is a non-trivial strengthening, as being right cycle-free implies being right cancellative, but the converse does not hold. Valitskas' proof of the strengthening was never published; a proof was given later by Guba \cite[Theorem~4]{Guba1994}. Gerasimov \cite{Gerasimov1982} has also given some necessary and sufficient conditions for a monoid to embed in a group. For a broad overview of various embeddability criteria for monoids and general algebra, we refer the reader to the survey by Bokut' \cite{Bokut1987}.

Returning to cycle-free presentations and their relation to the word problem, we note that as it is decidable whether or not a monoid (presentation) is cycle-free, being cycle-free is \textit{not} a necessary condition for cancellativity; indeed, it the problem of deciding whether a given finitely presented monoid is cancellative is undecidable in general. However, this can also be seen by way of concrete example, as provided by Adian \cite{Adian1960b}\footnote{In that article, the second relation is misprinted as $aed = ced$. This is corrected in \cite[Theorem~II.7]{Adian1966}.} as follows
\[
\pres{Mon}{a,b,c,d,e}{ab=cd, aeb=ced}.
\]
As witnessed above in Example~\ref{Ex: double edges cycles}, it has both left and right cycles, but it is not hard to show that it is cancellative. One can also check that $ae^2b \neq ce^2d$ in this monoid, and so it cannot possibly embed in the group with the same presentation. 

Thus we have three properties for monoids, which in general are not equivalent:
\begin{enumerate}
\item being cancellative;
\item being group-embeddable;
\item being cycle-free.
\end{enumerate}
Hence, as the examples above show, we only have $(3) \implies (2) \implies (1)$, and the reverse implications can fail already in the case of two defining relations. However, in the case of one-relation monoids, it is not very hard to show that $(1) \implies (3)$, by induction on the number of elementary transitions required. Hence these three properties are all equivalent in the class of one-relation monoids. Summarising, we have the following theorem. 

\begin{theorem}[{\cite[Corollary~4 \& 5]{Adian1960b}}]
Let $M = \pres{Mon}{A}{u=v}$ be a one-relation monoid. Then the following are equivalent: 
\begin{enumerate}
\item $M$ is cancellative;
\item $M$ embeds in $\pres{Gp}{A}{u=v}$;
\item $M$ is cycle-free.
\end{enumerate}
Hence, if $M$ is cycle-free then $M$ has decidable word problem.
\end{theorem}

As a concrete example, any monoid of the form $\pres{Mon}{a,b}{aub = bva}$ is cancellative, and the word problem is decidable in any such monoid.\footnote{One-relation monoids of this form have been called \textit{Adian monoids} by various authors (see e.g. \cite{Ivanov2001, Inam2017}). As we shall see, this name could equally be applied to a number of other families of one-relation monoids.} This monoid embeds in $\pres{Gp}{a,b}{auba^{-1}v^{-1}b^{-1} = 1}$, to which we can readily apply Magnus' procedure for solving the word problem. We hence see that our first example of a non-trivial solvable word problem comes about as an example of Sushkevi\v{c}'s principle of reducing a semigroup problem to a group problem. With one non-trivial class of one-relation taken care of, there is one class that stands out among the rest as being particularly special. 

\subsection{Special monoids}\label{Subsec: Special monoids}

A monoid is called \textit{special} if it admits a presentation of the form $\pres{Mon}{A}{r_1 = 1, r_2 = 1, \dots, r_k =1}$. Of course, all groups are special monoids; in fact, given a $k$-relator group, it always admits a $(k+1)$-relation special monoid presentation by a simple trick introduced already by von Dyck \cite{Dyck1882}: if $G$ is given by the presentation $\pres{Gp}{a_1, \dots, a_n}{r_1 = 1, \dots, r_k = 1}$, where the $r_i$ are words in the $a_j$ and their inverses, we can add a single generator $x$ and a defining relation $a_1 a_2 \cdots a_{n} x a_{n} \cdots a_2 a_1 = 1$. We can now clearly rewrite the $w_i$ as words over the positive alphabet; for example, $a_1^{-1}$ is equal to $a_2 \cdots a_{n} x a_{n} \cdots a_2 a_1$. The resulting special monoid is isomorphic with $G$. There are, however, $k$-relator groups which are not $k$-relation special monoids (see \S\ref{Subsec: String rewriting sporadic}).

Now, not all special monoids are groups; the simplest example is the bicyclic monoid $\pres{Mon}{b,c}{bc=1}$. There $b$ is right (but not left) invertible, and $c$ is left (but not right) invertible. It is, in fact, not too hard to show that $M$ has no non-trivial invertible elements. Special monoids were first properly\footnote{However, Thue very explicitly solves the word problem for the special monoids  $\pres{Mon}{a,b,c}{abbcab=1}$ and $\pres{Mon}{a,b}{ababa=1}$ (both of which are easily seen to be isomorphic with free groups). For more details, see \cite[Examples~2 \& 3]{Thue1914}.} introduced to the literature and given their name by G. S. Tseitin in \tc{1958}, who named them \textit{special associative systems} (\textit{spetsialnaia assotsiativnaia sistema}). However, the first systematic study of special monoids came two years later, by Adian \cite{Adian1960}. 

We say that a special monoid $\pres{Mon}{A}{r_1 = 1, r_2 = 1, \dots, r_k=1}$ is \textit{$\ell$-homogeneous} if there exists $\ell \in \mathbb{N}$ such that $|r_i| = \ell$ for all $1 \leq i \leq k$. One defines $\ell$-homogeneous groups in precisely the same manner. In his paper, Adian gives a very thorough overview of the proofs of the following central result regarding homogeneous special monoids; the full proofs were provided six years later in his celebrated monograph, see \cite[III.Theorem~1]{Adian1966}. 

\begin{theorem}[Adian, \tc{1960}]
Let $M = \pres{Mon}{A}{r_1 = 1, r_2 = 1, \dots, r_k=1}$ be an $\ell$-homogeneous special monoid. If the word problem is decidable for all $k$-relator $\ell$-homogeneous groups, then the word problem and the divisibility problems are decidable for $M$.
\end{theorem}

Note that for fixed $k$ and $\ell$, there are only finitely many $k$-relator $\ell$-homogeneous groups (up to a free factor of a free group). Furthermore, every one-relation special monoid is $\ell$-homogeneous, yielding the following immediate corollary by applying Magnus' result.

\begin{corollary}[Adian, \tc{1960}]
Let $M = \pres{Mon}{A}{w=1}$. Then the word problem and divisibility problems for $M$ is decidable. 
\end{corollary}

Adian also proves certain undecidability results for $\ell$-homogeneous special monoids. In particular, he proves that for every $\ell > 3$, there exists an $\ell$-homogeneous special monoid with undecidable word problem. Note that in the case $\ell=2$, the word problem is trivially decidable, as every $2$-homogeneous group is a free product of finitely many cyclic groups (cf. \cite[III.\S5,~Theorem~11]{Adian1966}). The subgroup $U(M)$ of a monoid $M$ consisting of all invertible elements of $M$ is called the \textit{group of units} of $M$. When giving the full proofs in his monograph, Adian also proves (see \cite[III.\S4,~Theorem~8]{Adian1966}) that the group of units of an $\ell$-homogeneous $k$-relation special monoid $M$ is isomorphic with an $\ell$-homogeneous $k$-relator group, and that the word and divisibility problems for $M$ reduce to the word problem for $U(M)$ (though this latter reduction is not in general constructive). In particular, this gives the following beautiful result.

\begin{theorem}[Adian, \tc{1966}]
Let $M = \pres{Mon}{A}{w=1}$. Then the group of units $U(M)$ of $M$ is a one-relator group. 
\end{theorem}

In fact, Adian gives an algorithm for computing a presentation for the group of units of a one-relation special monoid \cite[III.\S4,~Theorem~7]{Adian1966} which Zhang \cite{Zhang1992b} noticed had unnecessary steps (this is discussed in greater detail below). We give a brief overview of the simplified algorithm here, following Kobayashi \cite{Kobayashi2000}. Let $M = \pres{Mon}{A}{w=1}$. Let $C_0 = \{ w \}$, and suppose, for induction, that $C_i$ has been defined for some $i \geq 0$. Let $x, y \in C_i$. If $x \equiv vu$ and $y \equiv uw$ for some words $u, v, w \in A^+$ (i.e. if $x$ and $y$ overlap in the word $u$) then we set 
\[
C_{i+1} = (C_i \setminus \{x, y \} ) \cup \{ u, v, w \}.
\]
As the length of $u, v$ and $w$ are all less than the lengths of $x$ and $y$, we must have that this process will stabilise eventually, giving us a finite (but not uniquely determined) sequence $C_0, C_1, \dots, C_k$. In this case, $C_k$ is a biprefix code (i.e. none of the words of $C_k$ overlap), and in fact it is not difficult to see that although the sequence is not, the set $C_k$ is uniquely determined by $w$. We set $C(w)$ to be this set, and call this the \textit{self-overlap free code generated by $w$}. In an entirely analogous fashion we can define the self-overlap free code $C(W)$ generated instead by a set of words $W$. 

\begin{example}
Let $w \equiv abbaab$. Let $C_0 = \{ abbaab\}$. Then choosing $x \equiv abbaab$ and $y \equiv abbaab$, we find that we can take $x \equiv (ab)(baab)$ and $y \equiv (abba)(ab)$, so we set 
\[
C_1 = (C_0 \setminus \{ abbaab \}) \cup \{ ab, baab, abba\} = \{ ab, baab, abba\}.
\]
Now we can pick $x \equiv ab$ and $y \equiv baab$, and find $x \equiv  (a)(b)$ and $y \equiv (b)(aab)$, so we set
\[
C_2 = (C_1 \setminus \{ ab, baab \}) \cup \{ a, b, aab\} = \{ a, b, aab, abba\}.
\]
Now we can repeatedly pick $x \equiv a$ followed by repeatedly picking $x \equiv b$ to remove both $aab$ and $abba$, and we eventually end up with a set 
\[
C(w) = \{ a, b \}.
\]
Thus $C(w) = \{ a, b \}$ is the self-overlap free code generated by $w \equiv abbaab$. 
\end{example}

\begin{example}
Let $w \equiv abcabdab$. Then the self-overlap free code generated by $w$ is 
\[
C(w) = \{ ab, cabd\}.
\]
In particular, $C(w)$ is not an \textit{infix} code in general, i.e. one might find words in $C(w)$ appearing as subwords of other words in $C(w)$. 
\end{example}

Returning to our one-relation monoid $\pres{Mon}{A}{w=1}$, if $C(w) = \{ w_1, w_2, \dots, w_n\}$, then let $X = \{ x_1, \dots, x_n \}$ be a set in bijective correspondence with $C(w)$ via the map $\varphi \colon w_i \mapsto x_i$. As $C(w)$ is a biprefix code, this can be uniquely extended to a homomorphism 
\[
\varphi \colon C(w)^\ast \to X^\ast.
\]
As $w$ is clearly a word in $C(w)^\ast$, we can consider the word $\varphi(w)$, and the monoid 
\[
\pres{Mon}{X}{\varphi(w) = 1}.
\]
It is not hard to see that this monoid is, in fact, a (one-relator) group. Clearly, every word in $C(w)^\ast$ is invertible, but the key insight by Adian is that every invertible word is equal in $M$ to some word from $C(w)^\ast$. In particular, $U(M)$ is isomorphic with the group with the above presentation.

\begin{example}
Let $M = \pres{Mon}{a,b,c,d}{abcabdab = 1}$. Then 
\[
C(w) = \{ ab, cabd\} = \{ w_1, w_2 \}.
\]
We factor $abcabdab$ uniquely as $(ab)(cabd)(ab)$, and hence find that 
\[
U(M) \cong \pres{Mon}{x_1, x_2}{x_1x_2x_1 = 1} \cong  \pres{Gp}{x_1, x_2}{x_1x_2x_1 = 1} \cong \mathbb{Z}.
\]
Thus the group of units of $M$ is infinite cyclic, and the word and divisibility problems are rather straightforward to solve in $M$ (see below). 
\end{example}

This algorithm is very simple to use in practice. We shall present the main idea of why the word problem can be reduced to the word problem for the group of units below, in the general setting of $k$-relation special monoids. 

Makanin \cite{Makanin1966b} in his Ph.D. thesis extended Adian's results from $\ell$-homogeneous special monoids to all special monoids (the results were announced in a bulletin article \cite{Makanin1966}). Specifically, he proved that the group of units of a $k$-relation special monoid $M$ is a $k$-relator group, and that the word problem and divisibility problems for a $k$-relation special monoid in which all defining relations have length $\leq \ell$ reduce to the word problem for all $k$-relator groups in which all defining relators have length $\leq \ell$. This solution is constructive, but Makanin also proves that the word and divisibility problems in $M$ reduce non-constructively to the word problem for $U(M)$. The non-constructibility comes from the difficulty of actually computing a presentation for $U(M)$. The author of the survey clarifies this issue in a forthcoming expository article on the subject \cite{NybergBrodda2021c}. We also remark that the author has recently translated Makanin's Ph.D. thesis into English \cite{NybergBrodda2021b}. We summarise the key theorem below. 

\begin{theorem}[Makanin, \tc{1966}]
Let $M = \pres{Mon}{A}{r_1 = 1, r_2 = 1, \dots, r_k=1}$ be a special monoid. If the word problem is decidable for $U(M)$, then the word problem and the divisibility problems are decidable for $M$. 
\end{theorem}

The fundamental idea behind the study of special monoids (not just the one-relation case) can be heuristically explained as follows: suppose that we have a word $w$ containing $r_1$ and $r_2$ as subwords, where $r_1=1$ and $r_2=1$ are some two defining relations of a special monoid. Suppose that these two occurrences have a non-trivial overlap. Then we can write $w \equiv w' r_1' s r_2'' w''$, where $r_1 \equiv r_1' s$ and $r_2 \equiv s r_2''$. As $s$ is a suffix of $r_1$, it is left invertible, and as it is a prefix of $r_2$, it is right invertible. Hence any overlap of defining relations must be invertible: in particular, if we factor $r_1$ and $r_2$ (necessarily uniquely) into minimal invertible factors as $r_1 \equiv \delta_1 \cdots \delta_k$ and $r_2 \equiv \delta_1' \cdots \delta_\ell'$, then we must have 
\[
s \equiv \delta_i \delta_{i+1} \cdots \delta_k \delta_1' \delta_2' \cdots \delta_j'
\]
for some $i, j \geq 1$. Hence, all resolutions of overlaps in a special monoid are actually resolutions of equalities of invertible words. Because of this crucial point, the reader familiar with the importance of resolving overlaps when solving the word problem in rewriting systems will at this point, perhaps, feel more confident in accepting that the word problem of a special monoid could somehow be reduced to the same problem for its group of units (while noting that the above heuristic is far from a proof!). 

A little more rigorously, there is certainly always a way one can factorise any relator word $r_i$ into certain \textit{minimal invertible pieces}, i.e. words which do not have any left or right factors which are invertible. To find this set, Adian and Makanin use certain \textit{extension operations} to compute this set, starting from $C(\cup_i \{ r_i \})$ (i.e. the self-overlap free code generated by the set of left-hand sides of the defining relations). These operations are based on using equalities in the group with the $k$-relator presentation obtained from factorising the pieces into words from $C(\cup_i \{ r_i \})$, which can produce new, shorter invertible words, giving a presentation for another group, and the process repeats. Eventually, one finds that this process terminates (though of course not constructively, unless the word problem is decidable in all groups one encounters along the way). This results in a finite set $\Delta$ consisting of all invertible words $w$ with no non-trivial left or right factors in $\Delta$ with the property that $w$ is equal to some invertible factor of some relation word $r_i$, and such that additionally $|w| \leq |r_i|$.\footnote{The reader familiar with Zhang \cite{Zhang1992a} will perhaps wish to substitute $|w| \leq \max_i |r_i|$ for this final inequality, but this ``global'' bound is superfluous.} One can then factor the relation words into words from $\Delta$ uniquely, and obtain a presentation for a group. This group turns out to be isomorphic with the group of units of $M$. 

Let us say we now wish to solve the word problem for $M$. We do this by finding normal forms of words in the following manner. Let $w \in A^\ast$ be a word. The key point in Adian's and Makanin's proofs is to find subwords of $w$ in $\Delta^\ast$ and replace them with equal and shorter words in $\Delta^\ast$, until this cannot be done anymore. They prove that the resulting form is unique, and hence reduce the word problem for $M$ to the comparison of words in $\Delta^\ast$. By construction, this is the word problem for comparison of words in $\Lambda^\ast$, which is just the word problem for $U(M)$. 

At this point, it is important to discuss the contributions to this area by L. Zhang in the \tc{1990}s, coming from theoretical computer science and the theory of rewriting systems. He noticed that the extension operations of Adian and Makanin, used to compute the set $\Delta$, is not necessary in the one-relation case. This observation reduces to the fact that given a one-relation monoid $\pres{Mon}{A}{w=1}$, it is not possible to have non-trivial equalities of distinct pieces with one another, cf. \cite[Lemma~1]{Zhang1992b}. This follows from Magnus' \textit{Freiheitssatz} for the one-relator group $\pres{Gp}{X}{\varphi(w) = 1}$. Beyond this, the proof -- which is now phrased in terms of rewriting systems -- is identical.

Zhang would thereafter write several other papers on special monoids, rephrasing many of the old results using rewriting systems. This brought the attention of Adian \cite{Adian1993}, who writes: 
\blockquote{
\textit{I was surprised recently seeing several papers of L. Zhang [...] published in well known mathematical journals. The large part of these papers looks like a result of rewriting from [Adian's monograph] in a direct meaning of the word. Of course, to refresh the ideas and the technique of [the monograph] may be useful, but it should be done in a more decent way.}
} 

The papers referenced by Adian are \cite{Zhang1991, Zhang1992a, Zhang1992b, Zhang1992c}; also relevant to this discussion are \cite{Otto1991, Zhang1992d, Zhang1996}. There is merit to Adian's criticism. For example, none of the results appearing in \cite{Zhang1992b} are new. They are all rather straightforward corollaries of Makanin's work, whether they appear implicitly or explicitly therein. For example, Zhang proves that the submonoid of right invertible elements of a special monoid is isomorphic with a free product of a free monoid by the group of units; this observation, though not phrased algebraically, is precisely what is used by Makanin to reduce the right divisibility problem to the word problem for the group of units \cite[I.\S3]{Makanin1966b}.

One feature common to Zhang's papers on special monoids (especially \cite{Zhang1992b}) is that they may appear succinct and elegant when compared to the difficult and lengthy inductive arguments of Adian and Makanin. This characterisation is also not entirely accurate; the difficult and lengthy inductive arguments are still included in Zhang's approach, but are hidden in the references to various confluence lemmas and results for general rewriting systems. On the contrary, the work by Adian and Makanin on special monoids is extensively self-contained. This can, of course, make a cursory reading of the material rather difficult: the chapter of Adian's monograph \cite{Adian1960} which deals with the word problem for special monoids ends with Lemma~111 (!), and the corresponding chapter of Makanin's Ph.D. thesis, which is rather streamlined by comparison, ends with Lemma~31. Thus the main benefit of reading \cite{Zhang1992b} comes from its expository nature. 

The other papers by Zhang on the topic do, however, contain some new results for special monoids, including the one-relation case. To this end we only expand on the \textit{conjugacy problem}. There are a number of generally inequivalent ways to define this problem for special monoids; Zhang \cite{Zhang1992a} proved that several natural such definitions coincide in the special case. The one we shall adopt here is the following: decide for two given words $u, v \in A^\ast$ whether $ux = xv$ holds for some $x \in A^\ast$. Using rewriting techniques, Zhang reduces the conjugacy problem to the same problem for the group of units. In particular, this shows that the conjugacy problem is decidable in $\pres{Mon}{A}{w^n=1}$, when $n >1$, as the conjugacy problem is decidable for one-relator groups with torsion by applying the B. B. Newman Spelling Theorem \cite{Newman1968} (cf. also \cite{NybergBrodda2021a}). 

We mention briefly that special monoids have recently been investigated by the author from the point of view of their geometry and formal language theory, cf. \cite{NybergBrodda2020a, NybergBrodda2020b}. 

\

In summary, we have thus seen two examples (special and cycle-free) of families of one-relation monoids with decidable word problem, in which the solution to the word problem reduces fairly directly to Magnus' result that the word problem is decidable in one-relator groups. At this point, however, there yet remained many cases and difficult reductions to make, which began to illustrate the complex nature of the problem. 

\clearpage

\section{Compression and reductions (\tc{1974}--\tc{1987})}

Following the initial success of Adian's work on cycle-free and special one-relation monoids, Adian and his doctoral student G. U. Oganesian would provide the next major step forward in the form of \textit{(weak) compression} in \tc{1978}. However, this method had in fact been discovered independently several years earlier by G. Lallement, for slightly different (but closely related) purposes; the authors seem to have been unaware of the others' work. The method allows one to reduce each of the word and divisibility problems for one-relation monoids with left and right cycles to the same problem for a one-relation monoid with a shorter defining relation. Lallement focused on the case when this resulting monoid was special -- this case would later be called \textit{subspecial} one-relation monoids by Kobayashi. The methods in this section ultimately lead to the situation of reducing the word and divisibility problems to one-relation monoids without left cycles or without right cycles. 

There are two types of compression. One is ``weak'', and is based on finding a self-overlap free word as a prefix and suffix of both words in the defining relation. This method has been used and analysed in-depth by Lallement \cite{Lallement1974}, a \tc{1978} paper by Adian \& Oganesian \cite{Adian1978}, Kobayashi \cite{Kobayashi2000}, Gray \& Steinberg \cite{Gray2019}, and Nyberg-Brodda \cite{NybergBrodda2020c}. There is also another, more general, form of compression, which bears some similarity to the first, and is ``stronger'' in the sense it only requires that the two words in the defining relation have some shared prefix and some shared suffix; this is featured only in a \tc{1987} paper by Adian \& Oganesian \cite{Adian1987} and the surveys by Lallement \cite{Lallement1988, Lallement1993, Lallement1995}. The titles of the two papers \cite{Adian1978, Adian1987} are very similar, and their years of publications are, for obvious reasons, rather easy to mistake for one another. For this reason, it is not uncommon to see references in the literature to the two papers confused. We remark that the terminology ``weak'' resp. ``strong'' compression does not appear elsewhere in the literature.

\subsection{Weak compression}\label{Subsec: Strong compression}

We shall present the theory of weak compression essentially as it appears in the work of Adian \& Oganesian \cite{Adian1978}, while blending in some notation from Kobayashi \cite{Kobayashi2000} (see \S\ref{Sec: Modern and Future} for more details on Kobayashi's work). We remark that compression can easily be defined for arbitrary monoids, not necessarily one-relation, but for brevity we only deal with the one-relation case below; see \cite{NybergBrodda2020c} for full details.

Let $M = \pres{Mon}{A}{u=v}$ be a one-relation monoid. Let $\alpha \in A^+$ be a self-overlap free word. If $u, v \in \alpha A^\ast \cap A^\ast \alpha$, then we say that $M$ is \textit{(weakly) compressible}. Suppose that $M$ is compressible. Consider the language $\Sigma(\alpha) = \alpha (A^\ast \setminus A^\ast \alpha A^\ast)$. Let
\[
X(\alpha) = \{ x_{w_1}, x_{w_2}, \dots, x_{w_i}, x_{w_{i+1}}\dots, \}
\] be a set of new symbols in bijective correspondence with the (infinitely many) elements of $\Sigma(\alpha)$ via the bijection $\alpha w_i \mapsto x_{w_i}$. It is not hard to show that $\Sigma(\alpha)$ is a suffix code (see Kobayashi \cite[Lemma~3.4]{Kobayashi2000} and Lallement \cite{Lallement1974} for details), i.e. no element of $\Sigma(\alpha)$ is a suffix of any other; thus, for any word in $\alpha A^\ast$, we can decode it uniquely, moving from left to right, as a product of elements of $\Sigma(\alpha)$. In particular, as $u, v \in \alpha A^\ast$, we can \textit{uniquely} factor the words in the defining relation as
\begin{align*}
u &\equiv w_{1,1} w_{1,2} \cdots w_{1,n} \alpha, \\
v &\equiv w_{2,1} w_{2,2} \cdots w_{2,m} \alpha
\end{align*}
where $w_{i, j} \in \Sigma(\alpha)$. We define the \textit{compressor function} $\psi_\alpha$ on the defining relation as:
\begin{align*}
\psi_{\alpha}(u) &= x_{w_{1,1}} x_{w_{1,2}} \cdots x_{w_{1,n}}, \\
\psi_{\alpha}(v) &= x_{w_{2,1}} x_{w_{2,2}} \cdots x_{w_{2,m}}.
\end{align*}
Let $X \subseteq X(\alpha)$ be the (finite) set of all $x_{i, j}$ which appear in one of the above two factorisations. Then we define the \textit{left monoid} $L(M)$ associated to $M$ to be
\[
L(M) = \pres{Mon}{X}{\psi_{\alpha}(u) = \psi_{\alpha}(v)}.
\]
This definition looks rather abstract, but in practice it is very simple to find. 

\begin{example}
Consider the one-relation monoid 
\[
M = \pres{Mon}{a,b}{ab ba ab bb ab bb ab = abbaab}.
\]
Then the self-overlap free word $\alpha = ab$ is both a prefix and a suffix of each word in the defining relation, and obviously it is the unique such word. We factor each word uniquely into elements from $\Sigma(ab) = ab (A^\ast \setminus A^\ast ab A^\ast)$ as:
\begin{align*}
abbaabbbabbbab &\equiv ab\underline{ba} \cdot ab\underline{bb} \cdot ab\underline{bb} \cdot ab, \\
abbaab &\equiv ab\underline{ba} \cdot ab.
\end{align*}
Thus, the compressor $\psi_{ab}$ as applied to each word gives:
\begin{align*}
\psi_{ab}(abbaabbbabbbab) &= x_{ba} x_{bb} x_{bb}, \\
\psi_{ab}(abbaab) &= x_{ba}.
\end{align*}
The left monoid of $M$ thus has the presentation 
\[
L(M) = \pres{Mon}{x_{ba}, x_{bb}}{x_{ba} x_{bb} x_{bb} = x_{ba}} \cong \pres{Mon}{c,d}{cdd = c}.
\]
The word problem in $L(M)$ is now easily decidable by using the complete rewriting system $(cdd \to c)$; as we shall see below (Theorem~\ref{Thm: Strong compression}), this implies that the word problem for $M$ is decidable. 
\end{example}

We note some immediate properties of $L(M)$. The defining relation $\psi_{\alpha}(u) = \psi_{\alpha}(v)$ is shorter than the defining relation $u=v$, so there is basis for the name \textit{compression}. The relation $\psi_{\alpha}(u) = \psi_{\alpha}(v)$ could itself, of course, be compressible. There is also a natural way to extend $\psi_\alpha$ to other words than just the defining relations; this behaves similarly to a homomorphism, there is a close connection between words equal in $M$ and words equal in $L(M)$. We do not give the full details of this here, as it is already very well written in Adian \& Oganesian's \cite{Adian1978} original paper. 

\begin{theorem}[Lallement \tc{1974}, Adian \& Oganesian \tc{1978}]\label{Thm: Strong compression}
Let $M$ be a compressible one-relation monoid. Then each of the word problem and the divisibility problems reduces to the respective problem in $L(M)$. 
\end{theorem}

The proof of this theorem, while occasionally notationally somewhat technical, is not conceptually difficult at all, and bears some similarity to the situation of special monoids (see \S\ref{Subsec: Special monoids}). Indeed, a resolution of an overlap of defining relations in $M$ will clearly be a resolution of an overlap of defining relations in $L(M)$ -- as $\alpha$ is self-overlap free and $u, v \in \alpha A^\ast \cap A^\ast \alpha$ -- and the proof heuristic given for special monoids applies in this case too. We do not give the full details of the proof here, but give an instructive example.

\begin{example} Let $M = \pres{Mon}{a,b}{abaabbab = abbabaab}$. Then $M$ is weakly compressible with respect to $\alpha \equiv ab$, and we find that 
\begin{align*}
\psi_\alpha(abaabbab) &= x_a x_b, \\
\psi_\alpha(abbabaab) &= x_b x_a.
\end{align*}
Hence we expect $M$ to behave as the free commutative monoid 
\[
L(M) = \pres{Mon}{x_a, x_b}{x_a x_b = x_b x_a}.
\]
Indeed, it is not hard to find a structure coarsely resembling the right Cayley graph of the free commutative monoid on two generators inside that of $M$. The geometric aspects of weak compression are interesting in their own right, but are beyond the scope of the present survey. Full details will soon appear in work by the author.
\end{example}

We remark that Theorem~\ref{Thm: Strong compression} and the general method of weak compression as it appears in Adian \& Oganesian's article \cite[Theorem~3]{Adian1978} appears to have been solely due to Adian; a footnote on the first page of the article specifies precisely which of the two authors contributed which theorem.

The astute reader will have noticed the partial attribution of Theorem~\ref{Thm: Strong compression} to Lallement. We give a brief overview of weak compression as used by him. Lallement's primary interest appears to have been in characterising which one-relation monoids have non-trivial elements of finite order, similar to Fischer et al.'s \tc{1972} characterisation in the one-relator group case \cite{Fischer1972}. Lallement considered monoids of the form
\[
\pres{Mon}{A}{u = v} \quad \text{such that } u \in v A^\ast \cap A^\ast v,
\]
and proved that the word problem and divisibility problems in this case reduce to the word problem for a special monoid. The monoids of the above form are called \textit{subspecial}. The special monoid obtained from this reduction is just the left monoid of the monoid in question, although Lallement's compression technique, which is in a ``single step'', differs slightly from the method using $\psi_\alpha$ as detailed above, and is rather technical in nature. Clearly, a subspecial monoid (with $v$ non-empty) can be weakly compressed in the earlier sense, by considering the maximal self-overlap free word which is a prefix and a suffix of $v$. By repeatedly compressing in the usual sense, one eventually obtains the same monoid as described by Lallement; see Gray \& Steinberg's treatment \cite{Gray2019} of compression for this connection. We do not detail this method here, but mention that Zhang \cite{Zhang1992c} rewrote Lallement's proof using the language of rewriting system, and used this to solve the conjugacy problem in many cases of one-relation monoids with non-trivial elements of finite order. 

As mentioned, one of the reasons for Lallement's interest in subspecial monoids came from characterising the one-relation monoids which have non-trivial elements of finite order. First, the following result (see \cite[Corollary~2.5]{Lallement1974}) completely characterises the existence of non-trivial (i.e. $\neq 1$) idempotents in one-relation monoids. 

\begin{theorem}[Lallement, \tc{1974}]
Let $M = \pres{Mon}{A}{u=v}$. Then $M$ has a non-trivial idempotent if and only if $M$ is subspecial or special, but not a group.
\end{theorem}

Using this, he is able to completely characterise the one-relation monoids with elements of finite order as the subspecial monoids of a particular form, see \cite[Theorem~2.7]{Lallement1974}. 

\begin{theorem}[Lallement, \tc{1974}]
A one-relation monoid has non-trivial elements of finite order if and only if its presentation is of the form $\pres{Mon}{A}{(uv)^mu = (uv)^nu}$, where $uv$ is not a proper power in the free monoid, and $m > n \geq 0$. 
\end{theorem}

Lallement's goal is clearly algebraic in nature, as opposed to Adian \& Oganesian, who focused on the decision problems for the monoids involved. He also proves many statements regarding the residual finiteness of one-relation monoids, which we shall not expand on here, beyond mentioning that one can derive many more powerful algebraic statements and properties regarding one-relation monoids using weak compression. This algebraic connection has been exploited by Kobayashi \cite{Kobayashi2000} (see \S\ref{Sec: Modern and Future} and Gray \& Steinberg \cite{Gray2019}. The author of the present survey has also investigated the language-theoretic aspects of compression \cite{NybergBrodda2020c}. For example, it is decidable whether a one-relation monoid with a non-trivial idempotent has context-free word problem (see \cite{NybergBrodda2020c} for the relevant definitions).

At this point, we make the crucial point that weak compression is not enough to reduce the word problem for all one-relation monoids to the left cancellative case. Indeed, the one-relation monoid 
\[
\pres{Mon}{a,b}{aab = ab}
\]
is not (weakly) compressible, but it has both left and right cycles. To deal with these cycles, we therefore need a more powerful tool. 

\subsection{Strong compression}\label{Subsec: Weak compression}

For strong compression, we shall use an encoding function $\tau_k$, which is applicable to any one-relation monoid with left and right cycles, unlike weak compression. This encoding function already (rather confusingly) appeared in the \tc{1978} paper by Adian \& Oganesian. We follow the definition as given by Adian \& Oganesian in their \tc{1987} article introducing this encoding \cite{Adian1987}. Let $E_1, E_2, \dots, E_n$ be all words of length $k>0$ in an alphabet $A$ with $|A|=m$; note that $n = m^k$. Let $e_1, e_2, \dots, e_n$ be letters in some new alphabet. We shall define $\tau_k(w)$ for words $w \in A^\ast$. If $|w| < k$, then we set $\tau_k(w) = 1$. If $w \equiv E_i$ for some $1 \leq i \leq n$, then we set $\tau_k(w) = e_i$. If $w$ begins with $E_i$, and $w \equiv aw'$ for some letter $a \in A$, then we set $\tau_k(w) \equiv e_i \tau_k(w')$. 

\begin{example}
Let $A = \{ a, b\}$ and $k = 2$. Let $E_1 \equiv aa, E_2 \equiv ab, E_3 \equiv ba, E_4 \equiv bb$. Then we consider the encoding function $\tau_2$, which operates as follows:
\begin{align*}
\tau_2(aababab) &\equiv e_1 \tau_2(ababab) \\
&\equiv e_1 e_2 \tau_2(babab) \\
&\equiv e_1 e_2 e_3 \tau_2(abab) \\
&\equiv e_1 e_2 e_3 e_2 \tau_2(bab) \\
&\equiv e_1 e_2 e_3 e_2 e_3 \tau_2(ab) \\
&\equiv e_1 e_2 e_3 e_2 e_3 e_2.
\end{align*}
Obviously, the particular ordering chosen on $\{ aa, ab, ab, bb \}$ is of no significance. 
\end{example}

Now the encoding function $\tau_k$ is not, in general, a homomorphism, but it has many similar properties: let $u, v \in A^\ast$ be arbitrary words such that $|v| \geq k-1$. Write $v \equiv v' v''$ such that $|v'| = k-1$. Then we clearly have
\[
\tau_k(uv) \equiv \tau_k(uv') \tau_k(v''),
\]
and a similar statement is easy to make regarding $\tau_k(uv)$ when $|u| \geq k-1$. Using this property of $\tau_k$, one finds the following quite striking lemma, which is almost (see the remark following the lemma) the statement of \cite[Lemma~1]{Adian1987}.

\begin{lemma}\label{Lem: weak compression.}
Suppose $M = \pres{Mon}{A}{u=v}$ has both left and right cycles. Let $C$ be the maximal common prefix of $u$ and $v$, and let $D$ be the maximal common suffix of $u$ and $v$. Let $k = 1 + \min(|C|, |D|)$. Then the word problem and the divisibility problems for $M$ can be reduced to the corresponding problems for the monoid
\[
M_\tau := \pres{Mon}{e_1, e_2, \dots, e_n}{\tau_k(u) = \tau_k(v)},
\]
The relation $\tau_k(u) = \tau_k(v)$ has no left cycles if $|C| \leq |D|$, and no right cycles if $|C| \geq |D|$. \\ If $|C| = |D|$, then the word problem for $M$ is decidable. 
\end{lemma}
\begin{remark}
The statement of \cite[Lemma~1]{Adian1987} concludes that in the case $|C| = |D|$ the divisibility problems for $M$ are decidable, as in this case $M_\tau$ is cycle-free. It was, at the time, believed that there was a proof that any cycle-free monoid has decidable divisibility problems; as we shall expand on in \S\ref{Subsec: An incorrect proof}, this belief should be regarded as conjecture. 
\end{remark}

The monoid $M_\tau$ in the statement of Lemma~\ref{Lem: weak compression.} is said to be obtained from $M$ by \textit{strong compression}; this terminology is justified by the fact that any weakly compressible monoid is also strongly compressible. The proof of the lemma is not particularly difficult or long, and is self-contained enough that there is little value in reproducing it here. The interested reader may find the original proof repeated word for word in Adian \& Durnev's survey \cite[Lemma~4.5]{Adian2000}. We remark that (just as in the case of weak compression) there is an interesting language-theoretic interpretation of strong compression. This will appear in future work by the author.

\begin{example}\label{Example: tau compression}
Let $M = \pres{Mon}{a,b}{abaababb = abbaabb}$. Then $M$ has left and right cycles. The maximal common prefix of $abaababb$ and $abbaabb$ is $ab$, and the maximal common suffix is $abb$. Thus we will consider the encoding function $\tau_k$, where $k = 1 + \min(|ab|, |abb|) = 3$. Suppose we enumerate the $2^3$ words of length $3$ as 
\[
(E_1, E_2, E_3, E_4, E_5, E_6, E_7, E_8) = (aaa, aab, aba, abb, baa, bab, bba, bbb).
\]
Using this, we apply the encoding function $\psi_k$ to the defining relation to find
\begin{align*}
\psi_3(abaababb) &= e_3 \psi_3(baababb) \\
&= e_3 e_5 \psi_3(aababb) \\
&= e_3 e_5 e_2 \psi_3(ababb) \\
&= e_3 e_5 e_2 e_3 \psi_3(babb) \\
&= e_3 e_5 e_2 e_3 e_6 \psi_3(abb) \\
&= e_3 e_5 e_2 e_3 e_6 e_4,
\end{align*}
and, with the details assigned to the diligent reader, also
\[
\psi_3(abbaabb) = e_4 e_7 e_5 e_2 e_4.
\]
Thus the weakly compressed monoid $M_\tau$ has the presentation 
\[
M_\tau = \pres{Mon}{e_1, \dots, e_7}{e_3 e_5 e_2 e_3 e_6 e_4 = e_4 e_7 e_5 e_2 e_4},
\]
which obviously has no left cycles. The word problem and the divisibility problems for $M$ reduce to the same problems for $M_\tau$. Note that the word problem is trivially solvable in $M_\tau$, as the left-hand side of the relation is self-overlap free. Hence the word problem is decidable in $M$. In fact, the left divisibility problem is also decidable in $M_\tau$, by a result of Oganesian (see \S\ref{Subsec: Applications of A}).
\end{example}

By using Lemma~\ref{Lem: weak compression.} one easily sees that the word problem for any one-relation monoid can be reduced to the word problem for a one-relation monoid with no left cycles or to one with no right cycles. By using the standard symmetry argument, we hence have the following theorem.

\begin{theorem}[Adian \& Oganesian \tc{1978}]\label{Thm: Reduction (but not 2-gen)}
The word problem resp. the left divisibility problem for an arbitrary one-relation monoid reduces to the same problem for a one-relation monoid of the form
\[
\pres{Mon}{A}{bua = ava}
\]
or of the form
\[
\pres{Mon}{A}{bua = a},
\]
where $a, b \in A$ are distinct letters and $u$ and $v$ are arbitrary words. 
\end{theorem}

We remark that the statement of this theorem already appears as \cite[Theorem~5]{Adian1978}, and is proved by a method devised by Oganesian\footnote{The aforementioned footnote in that article specifies that this method is due solely to Oganesian.}. As the above formulation using Lemma~\ref{Lem: weak compression.} becomes significantly cleaner and quicker, we use this instead. Theorem~\ref{Thm: Reduction (but not 2-gen)} is also (essentially) presented by Howie \& Pride \cite[Corollary~4]{Howie1986}, though without a detailed proof and with some non-essential cases included. The paper of Howie \& Pride uses geometric techniques of diagrams, and their results were arrived at independently of the work by Adian \& Oganesian. 

\subsection{Two generators}

We now have two types of compression, with strong compression being the more general, and allowing us to reduce the word problem for all one-relation monoids to the same problem for left-cycle free one-relation monoids. However, there is a further significant reduction that can be made, which is to reduce to the two-generator case. 

Now, \textit{a priori} it is not difficult to reduce the word problem for a finitely presented $k$-generator monoid to the word problem for a two-generator monoid; indeed, if one has a finitely presented monoid on the generators $a_1, \dots, a_k$, then the mapping 
\[
a_i \mapsto aba^{i+1}b^{i+1}
\] will define an injective homomorphism into a finitely presented monoid given by the generators $a, b$. Thus the word problem for a $k$-generator monoid reduces to the word problem for a two-generator monoid, as decidability of the word problem is inherited by submonoids. However, even if one starts with a left cycle-free one-relation monoid, the embedding is not necessarily into a left cycle-free monoid, the defining relations of the resulting monoid might be significantly longer, and decidability of the divisibility problems is not inherited by taking submonoids. Thus we will require a more careful reduction if we are to reduce the word and divisibility problems to two-generated one-relation monoids without left cycles. This is precisely what is provided in the first part of Adian \& Oganesian's  article \cite{Adian1987}. 

\begin{theorem}[Adian \& Oganesian, \tc{1987}]\label{Theorem: Reduce n-rel to n+1-gen}
Suppose that 
\[
M = \pres{Mon}{A}{u_1 = v_1, \dots, u_n = v_n}
\]
is a left cycle-free monoid. Then one can find a monoid 
\[
\overline{M} = \pres{Mon}{a_1, \dots, a_{n+1}}{\overline{u}_1 = \overline{v}_1, \dots, \overline{u}_n = \overline{v}_n}
\]
such that $|u_i| = |\overline{u}_i|$ and $|v_i| = |\overline{v}_i|$ for all $1 \leq i \leq n$, and each of the equality and left divisibility problems for $M$ reduces to the same problem for $\overline{M}$. 
\end{theorem}

We shall not give a proof of this theorem, as it uses the algorithm $\mathfrak{A}$ (see \S\ref{Subsec: A}) and is in any case already very clearly written in Adian \& Oganesian's article \cite{Adian1987}. However, we will give the method by which one finds the words $\overline{u}_i$ and $\overline{v}_i$.

Let $M$ be as given in the statement of Theorem~\ref{Theorem: Reduce n-rel to n+1-gen}. Consider the left graph $\mathcal{L}(M)$ of $M$ (see \S\ref{Subsec: Cancellative}). By assumption $\mathcal{L}(M)$ has no cycles. If the letter $a$ appears at the beginning of some $u_i$ or $v_i$, then we say that the letter $a$ is \textit{leftist}. Of course, any non-leftist letter is isolated. Denote the connected components of $\mathcal{L}(M)$ as $L_1, \dots, L_k$. For every $L_j$ we will choose one letter $c_j$ appearing in this component. Then, for every word $w \in A^\ast$, we let $\overline{w}$ be the result of replacing in $w$ every occurrence of the letter $c_j$ by $c_1$. Then obviously $|w| = |\overline{w}|$, and as $\mathcal{L}(M)$ has no cycles, the monoid $\overline{M}$ with presentation as in the statement of Theorem~\ref{Theorem: Reduce n-rel to n+1-gen} will clearly have that $\mathcal{L}(\overline{M})$ is cycle-free. In fact $\mathcal{L}(\overline{M})$ consists of a single component, and so clearly has $n+1$ letters. 

\begin{example}
We continue our example given in Example~\ref{Example: tau compression}. Recall that 
\[
M_\tau = \pres{Mon}{e_1, \dots, e_7}{e_3 e_5 e_2 e_3 e_6 e_4 = e_4 e_7 e_5 e_2 e_4}.
\]
Now $\mathcal{L}(M)$ has a single edge between the two leftist vertices $e_3$ and $e_4$. We choose $c_3 \equiv e_3$ (though of course the choice $c_3 \equiv e_4$ works the same), and all other choices of $c_j$ are forced, as all other vertices are isolated. Thus replacing the letters $c_j$ by $c_1$, we find 
\[
\overline{M_\tau} = \pres{Mon}{c_1, e_4}{c_1 c_1 c_1 c_1 c_1 e_4 = e_4 c_1 c_1 c_1 e_4} \cong \pres{Mon}{a,b}{bbbbba = abbba}.
\]
Each of the word problem and left divisibility problems for $M_\tau$ reduces to the same problem for $\overline{M_\tau}$ (and both problems are easily decidable).
\end{example}

In summary, by combining the two techniques of strong compression and Theorem~\ref{Theorem: Reduce n-rel to n+1-gen} to pass to the $2$-generated case, we find the following theorem. 

\begin{theorem}[Adian \& Oganesian, \tc{1987}]\label{Thm: Main reduction theorem}
The word problem resp. the left divisibility problem for an arbitrary one-relation monoid reduces to the same problem for a one-relation monoid of the form
\[
\pres{Mon}{a,b}{bua = ava}
\]
or of the form
\[
\pres{Mon}{a,b}{bua = a},
\]
where $u$ and $v$ are arbitrary words. 
\end{theorem}

We remark that Adian \& Durnev \cite[p. 242]{Adian2000} incorrectly claim that Theorem~\ref{Thm: Main reduction theorem} was first proved in \tc{1978}, referencing Adian \& Oganesian's \tc{1978} article \cite{Adian1978}. However, the only theorem of the sort proved in this article is Theorem~\ref{Thm: Reduction (but not 2-gen)}, as discussed above, which does not reduce to the $2$-generated case. Instead, Theorem~\ref{Thm: Main reduction theorem} first occurs in the literature as precisely the statement of \cite[Theorem~2]{Adian1987}, if one accounts for the fact that the monadic case $\pres{Mon}{a,b}{bua = a}$ was at that point believed to have been solved (see \S\ref{Subsec: An incorrect proof} for further discussion).

\begin{example}
We give a full example of combining several reduction techniques for a rather difficult-looking one-relation monoid. Consider the one-relation monoid 
\[
M_4 = \pres{Mon}{a,b,c,d}{abdadadacbaca = abdadabdaca}.
\]
Then $M_4$ is both weakly and strongly compressible; as weak compression reduces relation length significantly quicker, and is a bit more fun to use, we first use this. The (unique) self-overlap free word $\alpha$ such that both words in the relation begin and end with $\alpha$ is $\alpha \equiv a$. We find 
\begin{align*}
abdadadacbaca &\equiv a\underline{bd} \cdot a\underline{d} \cdot a\underline{d} \cdot a\underline{cb} \cdot a\underline{c} \cdot a, \\
abdadabdaca &\equiv a\underline{bd}\cdot a\underline{d} \cdot a\underline{bd} \cdot a\underline{c} \cdot a,
\end{align*}
and hence, when applying the compressor function $\psi_a$ we have
\begin{align*}
\psi_a(abdadadacbaca) &= x_{bd} x_d x_d x_{cb} x_c \\
\psi_a(abdadabdaca) &= x_{bd} x_{d} x_{bd} x_c.
\end{align*}
Hence the divisibility problems and word problem for $M_4$ reduce to the same problems for 
\begin{align*}
M_4' &= \pres{Mon}{x_{bd}, x_d, x_{cb}, x_c}{x_{bd} x_d x_d x_{cb} x_c =x_{bd} x_{d} x_{bd} x_c}\\
&\cong \pres{Mon}{a,b,c,d}{abbcd = abad}.
\end{align*}
This one-relation monoid still has left and right cycles. It is strongly (but not weakly) compressible with respect to the maximal common prefix $ab$ and suffix $d$. Hence we find $k = 1 + \min(2,1) = 2$, and we are using the function $\tau_2$. There are $4^2 = 16$ words of length $2$ in the alphabet $\{ a, b, c, d\}$, which we order as
\[
(aa, ab, ac, ad, ba, bb, bc, bd, ca, cb, cc, cd, da, db, dc, dd) 
\]
where the index $1 \leq i \leq 16$ of each word in this sequence indicates that it is word $E_i$. Thus
\begin{align*}
\psi_2(abbcd) &= e_2 \psi_2(bbcd) = e_2 e_6 \psi_2(bcd) = e_2 e_6 e_7 \psi_2(cd) = e_2 e_6 e_7 e_{12}, \\
\psi_2(abad) &= e_2 \psi_2(bad) = e_2 e_5 \psi_2(ad) = e_2 e_5 e_4.
\end{align*}
Hence the word problem and divisibility problems for $M_4'$ (and hence also for $M_4$) reduce to the same problems for 
\begin{align*}
M_4'' &= \pres{Mon}{e_2, e_4, e_5, e_6, e_7, e_{12}}{e_2 e_6 e_7 e_{12} = e_2 e_5 e_4} \\
&\cong \pres{Mon}{a,b,c,d,e,f}{adef = acb}, 
\end{align*}
where the redundant generators have not been included, as these split off in a free factor. This is now a right cycle-free one-relation monoid, but it has left cycles. Of course, the word problem and the left/right divisibility problem for $M_4''$ reduces to the word problem and the right/left divisibility problem for 
\[
M_4''' = \pres{Mon}{a,b,c,d,e,f}{feda = bca},
\]
by reading the words backwards. This is now a one-relation monoid without left cycles. We now apply Theorem~\ref{Theorem: Reduce n-rel to n+1-gen} to reduce the word and divisibility problems to a two-generated one-relation monoid without left cycles. The only leftist letters are $f$ and $b$. We choose $b$ to represent this component in $\mathcal{L}(M_4''')$. This finally reduces the word and divisibility problems for $M_4$ to the same problems for
\[
M_4'''' = \pres{Mon}{c_1, f}{fc_1c_1c_1 = c_1 c_1 c_1} \cong \pres{Mon}{a,b}{baaa = aaa}.
\]
In this monoid, solving these aforementioned problems is not hard, and so we have a solution to the same problems for our original (rather complicated-looking) monoid $M_4$. 
\end{example}

We have thus found our desired reduction theorem. The above Theorem~\ref{Thm: Main reduction theorem} is still, at the time of the writing of this survey, the strongest general reduction theorem for one-relation monoids. There are a number of important cases of left cycle-free one-relation monoids for which the word and left divisibility problems are known to be decidable. One of the major sources of such results is Adian's algorithm $\mathfrak{A}$, which we shall now present.

\clearpage

\section{Adian's algorithm $\mathfrak{A}$ (\tc{1976}--\tc{2001})}

From the reductions we have seen, the word problem for all one-relation monoids has been reduced to the left divisibility problem for one-relation monoids without left cycles. The algorithm $\mathfrak{A}$ was devised by Adian\footnote{The choice by Adian of the letter $A$ to denote this algorithm is much more likely to have been chosen for pragmatic reasons (being the first letter of the alphabet) rather than as an act of perceived self-importance. Nevertheless, it provides a clear precedent for pronouncing $\mathfrak{A}$ as \textit{Adian's algorithm.}} in \tc{1976} as an attempt to solve this latter problem. 

\subsection{Two generators, one relation}\label{Subsec: A} The algorithm $\mathfrak{A}$ is very general, and will be defined for all left cycle-free monoids; to get a feel for how the algorithm operates, we shall begin by considering how the algorithm operates on left cycle-free two-generated one-relation monoids. In view of Theorem~\ref{Thm: Main reduction theorem}, this is no (direct) restriction if one wishes to solve the word problem in one-relation monoids using the algorithm $\mathfrak{A}$. We borrow some aspects of the excellent exposition of the algorithm $\mathfrak{A}$ by Lallement \cite{Lallement1993}. Let 
\[
M = \pres{Mon}{a,b}{u=v}
\] be a left cycle-free one-relation monoid. Let $w \in \{a, b\}^+$. As $M$ has no left cycles, we can uniquely factor $w$ from the left into a product of maximal prefixes of $u$ and $v$. Using this, we describe the \textit{prefix decomposition} of $w$. We factor $w$ from the left into maximal prefixes of $u$ and $v$. If at some point during this factorisation a prefix happens to be $u$ or $v$, then the decomposition stops, and we call this prefix the \textit{head} of the prefix decomposition. A prefix decomposition without a head is called a \textit{headless} prefix decomposition.

We illustrate this by an example. Let $M = \pres{Mon}{a,b}{babba = aba}$, and consider the word $w_1 \equiv abbababab$. Then $w_1$ has $ab$ as a prefix, which is the maximal prefix that is also a prefix of either $babba$ or $aba$. Thus the decomposition begins with this prefix, which we denote as $ab \mid bababab$. The decomposition now continues on the word $bababab$. There, the next prefix is $bab$, so the decomposition continues as $ab \mid bab \mid abab$. Finally, the next prefix is $aba$, which is one of the words in the defining relation. This ends the decomposition. Thus $aba$ is the head, and we denote the prefix decomposition of $w_1$ as 
\[
ab \mid bab \mid \fbox{$aba$}\: b.
\]
Similarly, we can find the prefix decomposition of $w_2 \equiv abbabbba$ to be
\[
ab \mid babb \mid ba
\]
which is hence a headless prefix decomposition. Using the prefix decomposition, we can describe the algorithm $\mathfrak{A}$. Let $x \in \{ a, b \}$ be a letter. The algorithm $\mathfrak{A}$ attempts to decide whether $w$ is left divisible by $x$ as follows.

\begin{enumerate}
\item If $w$ begins with $x$, then \tc{stop}.
\item Otherwise, compute the prefix decomposition of $w$.
\begin{enumerate}
\item If this prefix decomposition is headless, then \tc{stop}.
\item Otherwise, replace the head by the other side of the defining relation, obtaining a new word $w'$, and \tc{goto} (1) with the new word $w'$ and the letter $x$.
\end{enumerate}
\end{enumerate}

It is not \textit{a priori} clear that this algorithm should do any better than the na\"ive procedure for checking if two words are equal in a monoid by successively attempting all possible elementary transformations. However, there is much more power in $\mathfrak{A}$ than it first appears. We first give a straightforward example of using this algorithm below.

\begin{example}[{\cite[Example~2.1]{Lallement1993}}]
Let $M = \pres{Mon}{a,b}{baababa = aba}$. We may ask the question: is $w \equiv abbaaababab$ left divisible by $b$? The prefix decomposition of $w$ is: 
\[
ab\mid baa \mid \fbox{$aba$} \: bab
\]
As the defining relation is $aba = baababa$, we replace the head $aba$ by $baababa$, and repeat. This gives the sequence:
\begin{align*}
ab\mid baa \mid \fbox{$aba$} \: bab &\longrightarrow abbaabaabababab  \\
ab\mid baaba \mid \fbox{$aba$}\: babab &\longrightarrow abbaababaababababab \\
ab\mid \fbox{$baababa$}\: ababababab &\longrightarrow ababaababababab \\
ab\mid \fbox{$aba$}\: ababababab &\longrightarrow abbaababaababababab \\
ab\mid \fbox{$baababa$}\: ababababab &\longrightarrow ababaababababab \\ 
\fbox{$aba$}\: baababababab &\longrightarrow baabababaababababab.
\end{align*}
Thus the algorithm has terminated with a witness for the fact that $z$ is left divisible by $b$, i.e. in $M$ we have $w = b \cdot aabababaababababab$.
\end{example}

It is not always the case that $\mathfrak{A}$ terminates. For this reason, it is perhaps slightly abusive to use the term ``algorithm'' for $\mathfrak{A}$. We give an easy example of this non-terminating behaviour happening below. 

\begin{example}[{\cite[Example~29]{Lallement1995}}]\label{Ex: Looping A}
Let $M = \pres{Mon}{a,b}{baabbaa = a}$, and consider the word $w \equiv bbaaa$. Applying $\mathfrak{A}$, we find:
\begin{align*}
b \mid baa \mid \fbox{$a$}\: &\longrightarrow bbaabaabbaa \\
b \mid baab \mid \fbox{$a$}\: abbaa &\longrightarrow bbaabbaabbaaabbaa \\
b \mid \fbox{$baabbaa$}\: bbaaabbaa &\longrightarrow babbaaabbaa \\
ba \mid b \mid baa \mid \fbox{$a$}\: bbaa &\longrightarrow \cdots 
\end{align*}
Thus we have entered a loop, for we have found an occurrence of the first prefix decomposition $b \mid baa \mid \fbox{$a$}\:$ inside the decomposition $ba \mid b \mid baa \mid \fbox{$a$}\: bbaa$.
\end{example}

It seems as if this looping behaviour would be difficult to deal with, as one might not know whether falling into an infinite loop would imply divisibility by a letter or not. However, the remarkable result obtained by Adian regarding $\mathfrak{A}$ is the following theorem, which is a combination of the arguments given in \cite[\S6]{Adian1976}. 

\begin{theorem}[Adian, \tc{1976}]\label{Thm: Algo A three cases}
The result of applying the algorithm $\mathfrak{A}$ to a word $w$ and a letter $x$ will be exactly one of the following three cases:
\begin{enumerate}
\item $\mathfrak{A}$ transforms $w$ into a word $xw'$ after finitely many steps.
\item $\mathfrak{A}$ ends in a headless prefix decomposition after finitely many steps.
\item $\mathfrak{A}$ loops indefinitely. 
\end{enumerate}
The word $w$ is left divisible by $x$ if and only if we are in case (1). \\ Furthermore, in this case $\mathfrak{A}$ produces the shortest proof of this.  
\end{theorem}

As an example of how this theorem deals with loops, we see that as $\mathfrak{A}$ loops in Example~\ref{Ex: Looping A} on input $w$ and $a$, we can conclude that $w$ is not left divisible by $a$. 

Now a solution to the left divisibility problem by a letter is easily equivalent to a solution to the left divisibility problem in left cancellative monoids, by induction on the length of the word by which one divides. Thus we \textit{almost} have a solution to the left divisibility problem for all two-generated left cycle-free one-relation monoids, and hence \textit{almost} have a solution to the word problem for all one-relation monoids by Theorem~\ref{Thm: Main reduction theorem}! Of course, the difficulty comes from the fact that there is at present no known general algorithm which detects when case $(3)$ of Theorem~\ref{Thm: Algo A three cases} occurs. Adian, however, conjectured that such an algorithm (which he calls $\mathfrak{B}$) exists, and therefore indirectly conjectures that the word problem is decidable for all one-relation monoids. This conjecture appears throughout his published works, and was repeated by Adian as late as \tc{2018} at a conference at the Euler International Mathematical Institute in Saint Petersburg. 

There have been some attempts to derive general methods for detecting loops in $\mathfrak{A}$ (i.e. attempts at producing $\mathfrak{B}$). For example, Lallement \cite[pp. 38--39]{Lallement1995} provides an algorithm which can always detect whether or not $\mathfrak{A}$ loops on an input word in the monoid given in Example~\ref{Ex: Looping A}, thus solving the word (and left divisibility) problem for this monoid. This detection is via regular languages and finite state automata. These methods of detecting loops have been expanded by Lallement's student J. Bouwsma in her Ph.D. thesis \cite{Bouwsma1993}, solving a number of more cases of the form $\pres{Mon}{a,b}{bua=a}$. In fact, Lallement \cite{Lallement1993} conjectures that their methods (using finite state automata) might be sufficient to detect all loops that can appear when applying $\mathfrak{A}$ to $\pres{Mon}{a,b}{bua = a}$. If this is true, then this would imply decidability of the word problem for such monoids, which remains open (see \S\ref{Subsec: An incorrect proof}).

\subsection{The general case}

We now describe the algorithm $\mathfrak{A}$ in the general case, as it is presented by Adian \cite{Adian1976}. Let $M = \pres{Mon}{A}{T}$ be a finitely presented monoid without left cycles. For all triples $(a, b, aW)$, where $a, b \in A$ are distinct letters and $W \in A^\ast$ is an arbitrary word, we will define the \textit{prefix decomposition} $R_\ell(aW, b)$ of $aW$ with respect to $b$ inductively on the length of $aW$. This prefix decomposition will not always exist. We will first reduce the definition of $R_\ell(aW, b)$ to when $a$ and $b$ are adjacent in $\mathcal{L}(M)$, and then present a definition similar to the two-generator one-relation case. For expositional reasons, we will simultaneously with the general case present the material in the special case of the left cycle-free monoid 
\[
\Pi = \pres{Mon}{\alpha, \beta, \gamma, \delta}{\alpha \gamma = \gamma \delta \alpha,  \gamma \gamma \beta = \beta\alpha\beta\delta\gamma},
\]
which has left graph $\mathcal{L}(\Pi)$ as below.
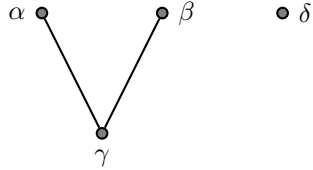
\[
\begin{tikzpicture}[>=stealth',thick,scale=0.8,el/.style = {inner sep=2pt, align=left, sloped}]%

\node (l0)[label=below:$\gamma$] [circle, draw, fill=black!50,
                        inner sep=0pt, minimum width=4pt] at (-3,0) {};
\node (l1)[label=right:$\beta$] [circle, draw, fill=black!50,
                        inner sep=0pt, minimum width=4pt] at (-2,2) {};
\node (l2)[label=left:$\alpha$] [circle, draw, fill=black!50,
                        inner sep=0pt, minimum width=4pt] at (-4,2) {};
\node (l3)[label=right:$\delta$] [circle, draw, fill=black!50,
                        inner sep=0pt, minimum width=4pt] at (0,2) {};
\path[-] 
    (l2)  edge (l0)
    (l0)  edge (l1);

\end{tikzpicture}
\]
First of all, if the left graph $\mathcal{L}(M)$ contains no path from the vertex $a$ to $b$, then we say that the prefix decomposition $R_\ell(aW, b)$ does not exist. Thus, for example, $R_\ell(\alpha\beta\gamma\delta, \delta)$ is not defined (for $\Pi$ as above). On the other hand, if there is such a path, then there is a unique shortest such path, as $M$ is left cycle-free. Let $c$ be the letter adjacent to $a$ on the shortest path joining $a$ to $b$ in $\mathcal{L}(M)$. Let $aE = cD$ be the defining relation corresponding to the edge between $a$ and $c$ on the specified path, accommodating an interchange of the left- and right-hand sides of the relation. We define 
\[
R_\ell(aW, b) := R_\ell(aW, c).
\]
For example, we have $R_\ell(\alpha \gamma^3 \delta^2, \beta) := R_\ell(\alpha \gamma^3 \delta^2, \gamma)$, as the unique shortest path joining $\alpha$ and $\beta$ in $\mathcal{L}(\Pi)$ has $\alpha$ adjacent to $\gamma$. 

We must now define $R_\ell(aW, c)$, which will eventually be a definition quite familiar to the reader who is accustomed to the two-generator one-relation case. Let $aF$ be the maximal common prefix of the words $aE$ and $aW$. Write 
\begin{align*}
aE &\equiv (aF)E_1 \\
aW &\equiv (aF)W_1. 
\end{align*}

We will now define $R_\ell(aW, c)$ based on the properties of the words $E_1, W_1$. 
\begin{enumerate}
\item If $E_1 \equiv \varepsilon$, i.e. if $aE \equiv aF$, we then define 
\[
R_\ell(aW, c) := \fbox{$aE$} \: W_1,
\]
and we say that $aE$ is the \textit{head} of the prefix decomposition, which is indicated by the box. Note that a head is always one side of a defining relation! We say that this head $aE$ is \textit{associated} to the relation $aE = cD$. For example, 
\begin{align*}
R_\ell(\alpha \gamma \gamma \gamma \beta, \gamma) &= \fbox{$\alpha \gamma$} \:  \gamma \gamma \beta \\
R_\ell(\gamma \delta \alpha \beta \beta, \alpha) &= \fbox{$\gamma \delta \alpha$} \: \beta \beta,
\end{align*}
and the heads are associated to the relations $\alpha \gamma = \gamma \delta \alpha$ resp. $\gamma \delta \alpha = \alpha \gamma$.
\item If $E_1 \not\equiv \varepsilon$ but $W_1 \equiv \varepsilon$, we then define 
\[
R_\ell(aW, c) := aE \mid 
\]
where $\mid$ indicates that $aE$ is the first prefix in the prefix decomposition. This decomposition is \textit{headless}, i.e. it has no head. For example, 
\begin{align*}
R_\ell(\beta \alpha \beta, \gamma) &= \beta \alpha \beta \mid \\
R_\ell(\gamma \delta, \alpha) &=  \gamma \delta \mid
\end{align*}
as $\beta \alpha \beta$ and $\gamma \delta$ are proper prefixes of their corresponding relation words (i.e. $\beta \alpha \beta \delta \gamma$ resp. $\gamma \delta \alpha$).
\item If\footnote{In both the Russian original and the English translation of Adian \cite{Adian1976}, this case has a typo; it reads ``if the words $E_1$ and $F_1$ are non-empty'', but should read ``if the words $E_1$ and $Z_1$ are non-empty'' (p. 616 resp. p. 382).}  $E_1 \not\equiv \varepsilon$ and $W_1 \not\equiv \varepsilon$, then we can write 
\begin{align*}
E_1 &\equiv qE_2 \\
W_1 &\equiv d W_2,
\end{align*}
for some letters $q, d \in A$ and words $E_2, W_2 \in A^\ast$. Then $d$ and $q$ are distinct. As $|dW_2| < |aW|$, we can by the inductive hypothesis determine whether or not $R_\ell(dW_2, q)$ exists. 
\begin{enumerate}[label=(\roman*)]
\item If $R_\ell(dW_2, q)$ does not exist, then we say that $R_\ell(aW, c)$, and consequently also $R_\ell(aW, b)$, does not exist.
\item If $R_\ell(dW_2, q)$ exists, and is of the form 
\[
H_1 \mid H_2 \mid \cdots \mid H_k \: \fbox{$R$} \: W'
\]
where the $H_i$ are the prefix components of the decomposition, and the head $R$ is associated to the relation $R=S$, then we define 
\[
R_\ell(aW, c) := aF \mid H _1 \mid H_2 \mid \cdots \mid H_k \: \fbox{$R$} \: W'.
\]
\end{enumerate}
\end{enumerate}
This completes the description of the prefix decomposition $R_\ell(aW, b)$. It is not difficult to check that this is uniquely defined. Note that the prefix decomposition $R_\ell(aW, b)$ can always be algorithmically computed for any pair of letters $a, b$ and any word $W \in A^\ast$, regardless of whether we have e.g. a solution to the word problem or not. 

\begin{example}\label{Ex: Pi-examples of PD}
Continuing our example with $\Pi$ as above, which for ease of access was defined by 
\[
\Pi = \pres{Mon}{\alpha, \beta, \gamma, \delta}{\alpha \gamma = \gamma \delta \alpha,  \gamma \gamma \beta = \beta\alpha\beta\delta\gamma}
\]
we give a complete computation of a prefix decomposition below. The reader is encouraged to follow along (perhaps with their own example).
\begin{align*}
R_\ell(\alpha \beta \alpha \gamma \gamma \beta \alpha \delta \delta, \beta) &= R_\ell(\alpha \beta \alpha \gamma \gamma \beta \alpha \delta \delta, \gamma) \\
&= \alpha \mid R_\ell(\beta \alpha \gamma \gamma \beta \alpha \delta \delta, \gamma) \\
&= \alpha \mid \beta \alpha \mid R_\ell(\gamma \gamma \beta \alpha \delta \delta, \beta) \\
&= \alpha \mid \beta \alpha \mid \: \fbox{$\gamma \gamma \beta$} \: \alpha \delta \delta.
\end{align*}
We can also compute 
\begin{align*}
R_\ell(\gamma \beta \gamma \delta \gamma \delta, \beta) &= \gamma \mid R_\ell(\beta \gamma \delta \gamma \delta, \gamma) \\
&= \gamma \mid \beta \mid R_\ell(\gamma \delta \gamma \delta, \alpha) \\
&= \gamma \mid \beta \mid \gamma \delta \mid R_\ell(\gamma \delta, \alpha) \\
&= \gamma \mid \beta \mid \gamma \delta \mid \gamma \delta \mid 
\end{align*}
which is a headless decomposition. On the other hand, 
\begin{align*}
R_\ell(\beta \alpha \beta \alpha, \gamma) = \beta \alpha \beta \mid R_\ell(\alpha, \delta)
\end{align*}
and as $R_\ell(\alpha, \delta)$ is not defined, it follows that $R_\ell(\beta \alpha \beta \alpha, \gamma)$ is not defined. 
\end{example}

We emphasise that it is not hard to check that the prefix decomposition $R_\ell$ as defined here coincides with the earlier defined prefix decomposition for two generators and one relation. We may now state Adian's algorithm $\mathfrak{A}$ in full generality for a left cycle-free monoid $M = \pres{Mon}{A}{R}$. Let $w \in A^+$ be an arbitrary word, and let $b \in A$ be any letter. The algorithm $\mathfrak{A}$ will attempt to decide if $w$ is left divisible by $b$ in $M$.

\begin{enumerate}
\item If $w$ begins with $b$, then \tc{stop}.
\item Otherwise, determine if the prefix decomposition $R_\ell(w, x)$ exists.
\begin{enumerate}
\item If $R_\ell(w, x)$ does not exist, then \tc{stop}.
\item If $R_\ell(w, x)$ exists, but is headless, then \tc{stop}.
\item If $R_\ell(w, x)$ exists and has a head, let $R=S$ be the relation associated to its head $R$. Replace the head $R$ by the word $S$ in $w$, obtaining a new word $w'$, and \tc{goto} (1) with the new word $w'$ and the letter $x$.
\end{enumerate}
\end{enumerate}

This is very similar to the algorithm $\mathfrak{A}$ presented earlier in the two-generated one-relation case; the only difference is the added step of ensuring that the prefix decomposition exists. Note that -- obviously! -- for every transformation $w \longrightarrow w'$ done by a step of $\mathfrak{A}$, we have $w = w'$ in $M$. We can now state Adian's Theorem~\ref{Thm: Algo A three cases} in full generality.

\begin{theorem}[Adian, \tc{1976}]\label{Thm: Algo A three cases General}
The result of applying the algorithm $\mathfrak{A}$ to a word $w$ and a letter $x$ will be exactly one of the following four cases:
\begin{enumerate}
\item $\mathfrak{A}$ transforms $w$ into a word $xw'$ (after finitely many steps).
\item $\mathfrak{A}$ transforms $w$ into a word $w'$ for which $R_\ell(w', x)$ does not exist.
\item $\mathfrak{A}$ transforms $w$ into a word $w'$ for which $R_\ell(w', x)$ is headless.
\item $\mathfrak{A}$ does not terminate. 
\end{enumerate}
The word $w$ is left divisible by $x$ if and only if we are in case (1). \\ Furthermore, in this case $\mathfrak{A}$ produces the shortest proof of this.  
\end{theorem}

The proof is rather short, and does not depend on any significant results not contained in the paper itself. We direct the reader to the paper for the proof (and pleasant reading).

\begin{example}
We again consider 
\[
\Pi = \pres{Mon}{\alpha, \beta, \gamma, \delta}{\alpha \gamma = \gamma \delta \alpha,  \gamma \gamma \beta = \beta\alpha\beta\delta\gamma}.
\]
We provide a few fully worked examples of applying the algorithm $\mathfrak{A}$. We begin with $w \equiv \alpha \beta \alpha \gamma \gamma \beta \alpha \delta \delta$ and the letter $\beta$. That is, we will try and decide whether or not $w$ is left divisible in $\Pi$ by $\beta$. Now as computed in Example~\ref{Ex: Pi-examples of PD}, we have 
\[
R_\ell(\alpha \beta \alpha \gamma \gamma \beta \alpha \delta \delta, \beta) = \alpha \mid \beta \alpha \mid \: \fbox{$\gamma \gamma \beta$} \: \alpha \delta \delta \longrightarrow \alpha \beta \alpha \beta \alpha \beta \delta \gamma \alpha \delta \delta,
\]
as the head $\gamma \gamma \beta$ corresponds to the relation $\gamma \gamma \beta = \beta \alpha \beta \delta \gamma$. We now apply $\mathfrak{A}$ to the word $\alpha \beta \beta \alpha \beta \delta \gamma \alpha \delta \delta$ and the (same as before!) letter $\beta$.

We compute 
\[
R_\ell(\alpha \beta \alpha \beta \alpha \beta \delta \gamma \alpha \delta \delta, \beta) = \alpha \mid \beta \alpha \beta \mid R_\ell(\alpha \beta \delta \gamma \alpha \delta \delta, \delta) 
\]
and since $R_\ell(\alpha \beta \delta \gamma \alpha \delta \delta, \delta)$ is not defined -- as $\alpha$ and $\delta$ are in distinct connected components of $\mathcal{L}(\Pi)$ -- it follows that $R_\ell(\alpha \beta \alpha \beta \alpha \beta \delta \gamma \alpha \delta \delta, \beta)$ is not defined. By Theorem~\ref{Thm: Algo A three cases General} we conclude that in $\Pi$ the word $w \equiv \alpha \beta \alpha \gamma \gamma \beta \alpha \delta \delta$ is \textbf{not} left divisible by $\beta$.

The other two examples in Example~\ref{Ex: Pi-examples of PD} are quicker to find conclusions about. As $R_\ell(\gamma \beta \gamma \delta \gamma \delta, \beta)$ is headless, it follows that $\gamma \beta \gamma \delta \gamma \delta$ is \textbf{not} left divisible by $\beta$. Similarly, as $R_\ell(\beta \alpha \beta \alpha, \gamma)$ does not exist, it follows that $\beta \alpha \beta \alpha$ is \textbf{not} left divisible by $\gamma$. 

Now, if we wish to decide whether $\gamma \gamma \beta \delta \alpha$ is left divisible by $\beta$, we simply compute
\[
R_\ell(\gamma \gamma \beta \delta \alpha, \beta) = \fbox{$\gamma \gamma \beta$} \: \delta \alpha,
\]
and thus applying a step of $\mathfrak{A}$, we transform $\gamma \gamma \beta \delta \alpha$ into $\beta \alpha \beta \delta \gamma \delta \alpha$. Hence, we conclude that $\gamma \gamma \beta \delta \alpha$ is left divisible by $\beta$ in $\Pi$. 

We leave as an exercise to the reader a proof that the algorithm $\mathfrak{A}$ always terminates for $\Pi$, solving the left divisibility problem for this monoid. 
\end{example}

\begin{example}
We can adapt the one-relation monoid $\pres{Mon}{a,b}{baabbaa = a}$ from Example~\ref{Ex: Looping A} to give the left cycle-free monoid 
\[
\Pi' = \pres{Mon}{\alpha, \beta, \gamma}{\beta \alpha \gamma \alpha \gamma \beta \beta \alpha \gamma \alpha \gamma = \alpha \gamma, \gamma = \alpha \gamma}.
\]
We leave the reader to check that $\mathfrak{A}$ will loop when checking whether $\beta \beta \alpha \gamma \alpha \gamma \alpha \gamma$ is left divisible by $\alpha$, just as $\mathfrak{A}$ loops in the monoid of Example~\ref{Ex: Looping A} when checking whether $bbaaa$ is left divisible by $a$. 
\end{example}

What is remarkable is that this algorithm $\mathfrak{A}$ has the potential to solve the left divisibility problem and the word problem in \textit{any} finitely presented monoid without left cycles -- not just the one-relation case. If this potential turns out to be fulfillable, this would indicate that the class of left cycle-free monoids is significantly more well-behaved than simply being left cancellative (as left cancellative monoids can have undecidable word problem). However, no criterion is yet known for deciding when $\mathfrak{A}$ loops; the general case is thus, in a sense, at present understood to approximately the same degree as the one-relation case. 

\subsection{Applications of $\mathfrak{A}$}\label{Subsec: Applications of A}

There have been a number of appearances of the algorithm $\mathfrak{A}$ in the literature, and it has spurred a good deal of research, having appeared in a large number of publications \cite{Adian1978, Adian1987, Oganesian1978, Oganesian1978b, Oganesian1979, Oganesian1982, Oganesian1984, Sarkisian1976, Sarkisian1979, Sarkisian1981, Bouwsma1993, Adian1994, Kobayashi1998, Watier1996, Watier1997}. It is also present in the brief articles by Adian \cite{Adian1993, Adian2005}, but no new results are presented therein. In this section, we shall mention some major classical results, the proofs of which depend critically on the algorithm $\mathfrak{A}$. The first concerns ``monadic one-relation monoids with torsion'', and was proved (see \cite{Oganesian1978}) only two years after the introduction of $\mathfrak{A}$. 

\begin{theorem}[Oganesian, \tc{1978}]
Let $M = \pres{Mon}{a,b}{(bu)^na = a}$ for some $n>1$ and $u$ arbitrary. Then the left divisibility (and hence the word) problem for $M$ is decidable.  
\end{theorem}

A result in the non-monadic case was proved by the same author in the same year (see \cite[Theorem~1]{Oganesian1978b}), but does not appear to have received much attention in the literature. Let $\sigma_a$ be the function which counts the number of occurrences of the letter $a$ in a word. 

\begin{theorem}[Oganesian, \tc{1978}]
Let $M = \pres{Mon}{a,b}{bua = ava}$. If either
\begin{align*}
\sigma_a(bua) &= \sigma_a(ava), \quad \text{or if} \\
\sigma_b(bua) &= \sigma_b(ava)
\end{align*}
then the left divisibility problem (and hence also the word problem) is decidable for $M$.
\end{theorem}

As an example, this shows that the word problem is decidable for the one-relation monoid $\pres{Mon}{a,b}{ababa = baba}$, as both sides have equally many occurrences of the letter $b$. The above theorem is one of the most general results in the non-monadic case. 

There are two further important contributions in the non-monadic case, both due to G. Watier. Consider $\pres{Mon}{a,b}{bua=ava}$. The case when $bua$ is self-overlap free received some attention by Adian \& Oganesian \cite{Adian1987}. In fact, a proof was claimed that whenever $bua$ is a factor of $ava$, the word problem is decidable; this result, however, should be regarded as conditional (see \S\ref{Subsec: An incorrect proof}). Nevertheless, Watier \cite{Watier1996} studies this case in his first of two articles on the subject, and is able to prove the following remarkable theorem.

\begin{theorem}[Watier, \tc{1996}]
Let $M = \pres{Mon}{a,b}{bua=ava}$. If the leftmost sequence of $b$'s in $bua$ is strictly longer than the others, then the word problem is decidable for $M$. 
\end{theorem}

This is a very general theorem. As noted by Watier, this solves the word problem in $\pres{Mon}{a,b}{b^ma^n = ava}$ for $m, n > 0$, with no conditions on the word $ava$. Note that $bua$ is always self-overlap free in the cases that the above theorem applies. The second article\footnote{At this point, Watier had been made aware of the significant gap in a result by Sarkisian, see \S\ref{Subsec: An incorrect proof}. The article itself is communicated by Adian, and Watier graciously thanks him at the end.} by Watier on the subject again concerns the case when $bua$ is self-overlap free. As mentioned, the case when $bua$ is a factor of $ava$ has been studied to some extent; Watier's article studies the case when $bua$ is \textit{not} a factor of $ava$. He is then able to show the following very general result.

\begin{theorem}[Watier, \tc{1997}]
Let $M = \pres{Mon}{a,b}{bua = ava}$. Suppose that $bua$ is self-overlap free and not a subword of $ava$. If $|ava| \geq |bua|^2$, then the word problem for $M$ is decidable. 
\end{theorem}

The article itself makes for very pleasant reading, and many of the results are phrased in terms of formal language theory and the theory of codes, while still being centered on the algorithm $\mathfrak{A}$. Watier's two results, along with the results by Oganesian mentioned above, represent the bulk of the progress hitherto made on the case $\pres{Mon}{a,b}{bua=ava}$. We mention in passing that Kashintsev \cite{Kashintsev1978b} proved that if $u$ is self-overlap free, and $|u| > |v|$, then the conjugacy problem is decidable in $\pres{Mon}{A}{u=v}$. Unlike for groups, however, decidability of the conjugacy problem is not sufficient for decidability of the word problem.

\subsection{An incorrect proof}\label{Subsec: An incorrect proof}

Up to this point, all results mentioned have been unconditional. However, around the early \tc{1990}s a gap was discovered in the proof of a theorem that had up to that point been generally accepted. The ``theorem'' was a result by O. A. Sarkisian, another doctoral student of Adian's, and is simple to state. 

\begin{claim}[Sarkisian, \tc{1981}]
The divisibility problems are decidable for all cycle-free monoids.
\end{claim}

This is a rather remarkable result; indeed, outside of the above claim, even the word problem is not known to be decidable for cycle-free monoids. There are many remarkable consequences of the result, which we shall list below. However, Oganesian discovered a gap in Sarkisian's proof \cite{Guba2020}. The presence of this gap is first mentioned in \tc{1994} by Adian \cite{Adian1994}, but by then the result had already been used in the proofs of a number of further results. As we shall see, this has led to a number of results in the literature either being entirely conditional, or else needing a repair (see e.g. \S\ref{Subsec: Inverse monoids} where such a repair is done in connection with one-relation special inverse monoids). 

The author has collected all results known to him to be conditional on this result. We define the \textit{Sarkisian hypothesis} \tc{sh} to be the above claim, i.e. that the divisibility problems are decidable for all cycle-free monoids.

\begin{theorem}
Suppose $\textsc{sh}$ holds. Then:
\begin{enumerate}
\item The word problem is decidable for every cycle-free semigroup.
\item The word problem is decidable for every cycle-free group \cite{Sarkisian1981}.\footnote{In the English translation of the article in which this result appears (\cite[Theorem~2]{Sarkisian1981}), the word ``group'' has been mistranslated as ``semigroup'' (!). The Russian original is correct.}
\item The isomorphism problem for one-relation monoids is decidable \cite{Oganesian1984}.
\item\label{Item: Monadic dec} The word problem for $\pres{Mon}{a,b}{bua = a}$ is decidable \cite{Oganesian1982}. 
\item Suppose $M$ is a strongly compressible one-relation monoid. If $|C|=|D|$, then the divisibility problems are decidable for $M$ (see \S\ref{Subsec: Weak compression}) \cite{Adian1987}.
\item\label{Item: SOF dec} If $u, v$ are such that $v$ is self-overlap free and $v$ is a factor of $u$, then the word and at least one of the divisibility problems are decidable for $\pres{Mon}{A}{u=v}$ \cite{Adian1987}.
\end{enumerate}
\end{theorem}

Without \tc{sh}, all of the above results should only be regarded as conditional, as no other proof of them is known. Note also that $(\ref{Item: SOF dec})$ in the above theorem generalises $(\ref{Item: Monadic dec})$ in virtue of the fact that the word and divisibility problems are easily decidable for $\pres{Mon}{a,b}{b^n = a}$. We remark that there are other statements which are highly specialised (see e.g.  \cite[Theorem~7]{Adian1993}) and rather involved to state, which are also conditional on \tc{sh}. We instead refer the reader to Adian \& Durnev \cite{Adian2000}, in which slightly detailed corrections that can be made to such statements.\footnote{In that survey, the crucial condition of left cycle-freeness is sometimes mistakenly omitted from theorems.}

V. S. Guba \cite[p. 1142]{Guba1997} has stated that, in his opinion, decidability of the divisibility problems for cycle-free monoids ``is likely to be still more complicated than the word problem for one-relator semigroups''. We shall see some results by Guba in \S\ref{Subsec: monadic case} which illustrates the potential difficulty in solving the word problem for monadic one-relation monoids (and thereby also indirectly the difficulty in proving that \tc{sh} holds). 

\clearpage

\section{Sporadic results}

In this section, we shall present some approaches to the word problem for one-relation monoid which are ``sporadic'' in nature, in that they are not necessarily founded in an attempt to solve the word problem for all one-relation monoids, but have nevertheless been fruitfully applied to many classes. One quick example is that Magnus' \textit{Freiheitssatz}\footnote{If $G = \pres{Gp}{A}{w=1}$ with $w$ cyclically reduced and $Y \subseteq A \cup A^{-1}$ excludes some letter appearing in $w$, then the subgroup of $G$ generated by $Y$ is freely generated by $Y$. }, which was integral in solving the word problem for one-relator groups, can with little difficulty be generalised to one-relation monoids; an elementary proof is given by Squier \& Wrathall \cite{Squier1983}. However, this does not yet appear to have led to any direct new insights regarding the word problem. We mention one related insight. 

Magnus classified all one-relator groups satisfying some non-trivial identity \cite{Magnus1930}, and Adian classified all special one-relation monoids satisfying some non-trivial identity \cite{Adian1966}.\footnote{A monoid $M$ is said to satisfy the identity $U(x_1, x_2, \dots, x_n) = V(x_1, x_2, \dots, x_n)$ if all equalities of the form $U(A_1, \dots, A_n) = V(A_1, \dots, A_n)$, obtained from replacing the variables $x_1, \dots, x_n$ by arbitrary elements $A_1, \dots, A_n$ from $M$, are true in $M$. An identity is said to be non-trivial if it does not hold in the free semigroup on two generators. We refer the reader to e.g. the survey by Shevrin \& Volkov \cite{Shevrin1985} or \cite[Chapter~II]{Burris1981} for more information on identities.} Using the \textit{Freiheitssatz} and the work by Adian and Magnus, L. Shneerson \cite{Shneerson1972, Shneerson1972a} was able to completely classify the one-relation monoids which satisfy some non-trivial identity.\footnote{The author thanks Mikhail Volkov for bringing these papers to his attention and sending him copies, and Lev Shneerson for his interest in the author's forthcoming translation of the two papers.} Using this classification, Vazhenin \cite{Vazhenin1983} proved that a non-special one-relation monoid has decidable first-order theory if and only if it is monogenic, or else generated by $a, b$ and with defining relation one of 
\[
ab = ba \: \mid \: ab = b^k \: (k \geq 1) \: \mid \: ba = aba \: \mid \: ab = bab^2 \: \mid \: a = bab \: \mid \: a^2 = b^2,
\]
or the right cycle-free equivalent of one of the above defining relations. This coincides with the class of non-special one-relation monoids satisfying some non-trivial identity. Cain et al \cite[Proposition~9.1]{Cain2009b} showed that the above cases are also precisely the cases in which a one-relation monoid admits an \textit{automatic} presentation. This gives an efficient solution to the word problem in all the above cases, and is of independent interest (though note that directly solving the word problem in any of the above specified one-relation monoids requires very little effort).

We shall now present some families of one-relation monoids which have been studied in a rather detailed manner using normal forms in a more or less ``sporadic'' fashion.

\subsection{Normal forms}

One example of using normal forms to solve the word problem comes from Jackson \cite{Jackson1986} who, seemingly inspired by the example $\pres{Mon}{a,b}{baaba = a}$ which escaped the methods of Howie \& Pride \cite{Howie1986}, proved that the one-relation monoids
\[
\pres{Mon}{a,b}{ba^nba = a} \quad (n \geq 0)
\]
admit a particular nice solution to their word problem via normal forms; the article is only two pages long, and consists of a quick proof via van der Waerden's trick that the normal forms as specified are correct. We remark that although Jackson modestly notes that ``the result here should be regarded as a special case of a more general result of Oganesian'', this modesty is unfounded, in view of the gap mentioned in \S\ref{Subsec: An incorrect proof}. Zhang \cite[\S6]{Zhang1992a} later remarked that the one can solve the word problem in any monoid as above by considering its right cancellative analogue $\pres{Mon}{a,b}{aba^nb = a}$ and noting that this admits a finite complete rewriting system with the two rules 
\[
\{ (aba^nb \to a), (a^{n+1}b \to aba^n)\}.
\]
A similar result regarding certain infinite families comes from Lallement \& Rosaz \cite{Lallement1994}, who prove that the monoids
\[
\pres{Mon}{a,b}{ba = a(ba)^na} \quad (n \geq 0)
\]
have decidable word problem, again using normal forms. The special case with $n=1$, i.e. $\pres{Mon}{a,b}{ba = abaa}$ was then studied by Jackson \cite{Jackson2001}, who proved that such monoids even have decidable submonoid membership problem. Jackson \cite{Jackson2002} also studied the submonoid membership problem for the \textit{Baumslag-Solitar semigroups} 
\[
\pres{Mon}{a,b}{ab^k = b^\ell a} \quad (k, \ell > 0)
\]
and proved that this problem is decidable. Note, however, that the word problem is easily solvable in such monoids as they are all cycle-free. 

Finally, Yasuhara \cite{Yasuhara1970} proved that the word problem for one-relation monoids
\[
\pres{Mon}{A \cup \{ t \}}{u = vtw}
\]
where $t$ does not appear in $u, v$, or $w$, and $|u| > \max(|v|, |w|)$. This result was strengthened by Oganesian \cite{Oganesian1979} to solve the word problem in the above situation, with no condition on the lengths $|u|, |v|$ and $|w|$. However, Yasuhara's proof shows that the equivalence class of any word is finite when the length condition holds, which can be used to solve a number of other decision problems in this case (such as the membership problem). We note that in the two-generator case, the above monoids are all of the form
\[
\pres{Mon}{a,b}{a^n = a^p b a^q}, \quad n > p + q.
\]
Normal form results are closely connected to rewriting systems by a result of Squier \cite{Squier1987}, see also Brown \cite{Brown1992}. Thus, we shall give some brief remarks regarding the interface between the theory of rewriting systems and the word problem for one-relation monoids. 

\subsection{Rewriting systems}\label{Subsec: String rewriting sporadic}

While the theory of rewriting systems has often been limited to solving the word problem for very particular examples of one-relation monoids, there are some important points to be made on their general applicability. We begin by noting that occasionally, the algebraic structure of the rewriting systems in question seems to be somewhat unduly ignored. For example, the \textit{Jantzen monoid} $\pres{Mon}{a,b}{abbaab = 1}$ has been the subject of a number of investigations after its \tc{1981} introduction by Jantzen \cite{Jantzen1981}. An explicit linear representation for it is found, and it is proved that no congruence class of words is a context-free language. However, it is not hard to see that this monoid is a group: $ab$ is a prefix and a suffix of $abbaab$, and hence $ab$ and $ba$ are invertible. Furthermore, this group is easily seen to be the Baumslag-Solitar group $\operatorname{BS}(1, -2)$, by using the free group automorphism induced by $a \mapsto ab^{-1}$ and $b \mapsto b$. Seeing that no congruence class of words is a context-free language in $\operatorname{BS}(1, -2)$ is very straightforward. The ``Jantzen monoid'' has, however, been studied to good effect from the point of view of admitting a finite complete rewriting system, where the situation has been seen to be similar to the same for the \textit{Greendlinger group} $\pres{Gp}{a,b,c}{abc = cba}$, see \cite{Greendlinger1960, Greendlinger1960b, Otto1984}.

Similarly, using rewriting techniques, Otto \cite{Otto1988} proves in \tc{1988} that there exists a one-relator group that is not a one-relation monoid. This group is just $\mathbb{Z} \times \mathbb{Z}$, but Magnus \cite{Magnus1930} -- in the very first paper on one-relator group theory -- uses the \textit{Hauptform des Freiheitssatzes} to prove that the only one-relator group presentation for this group is $\pres{Gp}{a,b}{[a,b]=1}$. Hence it is already clear that there can be no one-relation monoid presentation for $\mathbb{Z} \times \mathbb{Z}$, for such a presentation would also be a one-relator group presentation for the group, but the defining relation would be a positive word. 

However, one of the major benefits of rewriting is that their techniques can often provide a strong link between formal language theory and monoids. A full description of this is certainly beyond the scope of this survey, but we point the reader to some articles which the author found to be of particular relevance to the word problem for one-relation monoids, viz. \cite{Book1983, Otto1984a, Potts1984, Metivier1985, McNaughton1987, Jantzen1988, Bucher1988, Kurth1990, Narendran1991, Wrathall1992, Otto1992,Otto1995, Wrathall1995, Choffrut2010}. We also refer to the excellent survey by Book, Jantzen \& Wrathall \cite{Book1982} and to the monograph by Jantzen \cite{Jantzen1988} for a general introduction to the interface between these areas. We remark that the author of the present survey has also studied the language-theoretic properties of special monoids (see \S\ref{Subsec: Special monoids}), and proved that a special monoid has context-free word problem if and only if its group of units is virtually free \cite{NybergBrodda2020b}, generalising the \textit{Muller-Schupp theorem} from groups to all special monoids. 

Regarding obtaining finite complete rewriting systems for one-relation monoids, Pedersen \cite{Pedersen1984, Pedersen1989} introduced \textit{morphocompletion}. This is a method that automatically introduces new generators, and can be regarded as a rather powerful completion procedure. The idea is to attempt standard completion (e.g. Knuth-Bendix completion) and, if a confluent system is not thereby attained, it backtracks, and adds new generators to stop the non-confluent branching from appearing. His morphocompletion solves the word problem in a number of one-relation monoids. For example, he solves the word problem in the cases 
\[
\pres{Mon}{a,b}{buba^m = a^n} \quad (n \geq m+3)
\]
when the word $buba^m$ is self-overlap free. Pedersen's morphocompletion seems to be one of few methods that also applies to the case $\pres{Mon}{a,b}{bua=ava}$ when $bua$ is not self-overlap free. Even for small examples, it generally produces finite complete rewriting systems with around $30$--$50$ rules. It also seems rather unlikely that his methods will generalise to all one-relation monoids. For a thorough survey on other forms of completion, see Dershowitz \cite{Dershowitz1989}.

\subsection{Small overlap conditions}

We shall briefly mention a particularly general method for solving many decision problems, which has some use also in the one-relation case. Consider a monoid presentation with alphabet $A$ and with defining relations $(u_i, v_i)$, where the relation words $u_i, v_i \in A^+$ are assumed non-empty. A non-empty subword $p$ of a relation word is called a \textit{piece} if it appears in at least two distinct ways as a subword of relation words. For $n \geq 1$, we say that the presentation satisfies the \textit{small overlap hypothesis} $C(n)$ if no relation word can be written as a product of fewer than $n$ pieces. Remmers \cite{Remmers1980} proved that the word problem is decidable in any monoid presentation satisfying $C(3)$. In fact, this has been greatly extended; $C(3)$ is sufficient to imply decidability of the conjugacy problem \cite{Cummings2014}, and Kambites \cite{Kambites2009b} has shown that $C(4)$ is sufficient to solve even the \textit{rational subset} membership problem, which greatly generalises the divisibility problems, submonoid membership problem, and the word problem. In fact, in $C(4)$ monoids, the word problem can be solved in linear time by a result of Kambites \cite{Kambites2009a}, cf. also the recent normal form algorithm for $C(4)$-monoids devised by Mitchell \& Tsalakou \cite{Mitchell2021}. Kashintsev \cite{Kashintsev1992} has explored connections between small overlap conditions and embeddability of monoids into groups, as well as using small overlap techniques for solving the word problem in some classes of special monoids \cite{Kashintsev1978, Kashintsev1993}.

It follows from Kambites' work on genericity in small overlap monoids (see \cite{Kambites2011}\footnote{The author thanks Mark Kambites for bringing his paper to his attention.}) work that for any $n \geq 1$, the probability that a two-generated one-relation monoid satisfies the small overlap condition $C(n)$ tends to $1$ as the length of the relation increases. Hence all of the problems mentioned above are decidable for \textit{almost all one-relation monoids}. It follows that the small overlap argument provides a reasonable argument for conjecturing that the word problem (even the rational subset membership problem) is decidable for all one-relation monoids. We remark that a parallel definition of small overlap monoids was given by V. A. Osipova. She proved that in certain monoids (so-called $\geq \frac{1}{2}$-monoids) satisfying a type of overlap condition the word problem is decidable \cite{Osipova1968}. These methods were also applied by her to partially understand the solvability of equations (the \textit{Diophantine problem}) in $\geq \frac{1}{3}$-monoids \cite{Osipova1972, Osipova1973}. This latter problem is beyond the scope of this survey, but has very recently been investigated in the special one-relation case by Garreta \& Gray \cite{Garreta2021}.

\clearpage

\section{Modern and future results (\tc{1997}--present)}\label{Sec: Modern and Future}

This section will give a high-level overview of some modern results related to the word problem for one-relation monoids. For obvious reasons of space, we will not be giving detailed proofs (or sometimes definitions) of the results and concepts mentioned here, but pointers to further reading will be amply provided.

\subsection{Finite complete rewriting systems}

It is an open problem whether every one-relation monoid admits a finite complete rewriting system, i.e. whether for every one-relation monoid $M = \pres{Mon}{A}{u=v}$ there exists an alphabet\footnote{The surface group $\pres{Mon}{a,b}{abba=1}$ does not admit a finite complete rewriting system over $\{a, b \}$, but does admit one over an alphabet with three letters \cite{Jantzen1985}.} $B$ and a finite complete rewriting system $T \subseteq B^\ast \times B^\ast$ such that $M$ is isomorphic with $B^\ast / \lra{T}$. In general, this property is stronger than having decidable word problem. In this section, we will briefly present two finiteness properties which are closely connected to finite complete rewriting systems, and what is known about these in the one-relation case. 

Kobayashi \cite{Kobayashi1998} investigated homotopy finiteness properties for one-relation monoids. One central such property is \textit{finite derivation type} ($\FDT$), which was introduced in \tc{1987} by Squier \cite{Squier1987} (and, independently, Pride \cite{Pride1995}), though Kobayashi notes that the idea appears implicitly in work by Adian \cite{Adian1966, Adian1976}. The central idea behind $\FDT$ is as a homotopy finiteness property for equivalence classes of derivations in monoids, and is independent of the particular finite presentation chosen for a monoid. While the full details of $\FDT$ are beyond the scope of the current article, our main interest comes from the following result due to Squier: if a monoid admits a finite complete rewriting system, then it has $\FDT$. This gives a potential venue for proving that a monoid does \textit{not} admits a finite complete rewriting system (which otherwise is a very difficult task), or indeed for providing supporting evidence that, say, a given one-relation monoid admits a finite complete rewriting system.

Kobayashi \cite{Kobayashi1998} first used the algorithm $\mathfrak{A}$ in \tc{1998} to prove that every one-relation monoid presented by a non-subspecial relation has $\FDT$. In a subsequent paper, he then proves that subspecial one-relation monoids have $\FDT$, by using weak compression to reduce it to the special one-relation case, which can then be reduced to $\FDT$ for one-relator groups \cite{Kobayashi2000}. As one-relator groups have $\FDT$ by a \tc{1994} result of Cremanns \cite{Cremanns1994b}, it then follows that all one-relation monoids have $\FDT$. Further to this, if a monoid admits a finite complete rewriting system, it also satisfies a certain \textit{homological} finiteness property $\FP_\infty$ (by Kobayashi \cite{Kobayashi1990} and Squier \cite{Squier1987}). It is known that any monoid with $\FDT$ also has $\FP_3$ by Cremanns \& Otto \cite{Cremanns1994}, and for some time it was an open problem whether every one-relation monoid has $\FP_\infty$. This was very recently answered affirmatively by Gray \& Steinberg \cite{Gray2019} in \tc{2019}. This can be regarded as supporting the conjecture that every one-relation monoid admits a finite complete rewriting system. The proofs leading up to these results are also notable in their usage of (weak) compression and the algorithm $\mathfrak{A}$ not to solve the word problem, but to derive strong structural results. 

We end this section by noting that some partial progress has also recently been made in constructing explicit finite complete rewriting systems for certain one-relation monoids. To this end, the results by Cain \& Maltcev \cite{Cain2013a, Cain2013b} bear mentioning, as they show that all one-relation monoids of the form $\pres{Mon}{a,b}{b^\alpha a^\beta b^\gamma a^\delta b^\varepsilon a^\varphi = a}$ admit finite complete rewriting systems, where $\alpha, \beta, \gamma, \delta, \varepsilon, \varphi \geq 0$, i.e. where the ``relative length'' of the left-hand side is at most $6$. This result provides quite explicit solutions to the word problem for such monoids, but their methods do not seem to be easily generalisable to all monoids $\pres{Mon}{a,b}{bua = a}$. We note in passing that the smallest monadic one-relation monoid to which no result in the literature appears to be available to solve the word problem for is $\pres{Mon}{a,b}{bababbbabba = a}$. The author has not found a finite complete rewriting system for this monoid, but has solved the word problem for this monoid by other means. 

\subsection{Special inverse monoids}\label{Subsec: Inverse monoids}

In \tc{2001}, a landmark paper by Ivanov, Margolis \& Meakin appeared; this discovered a strong link between the word problem for one-relation monoids and special one-relator \textit{inverse} monoids. We give a brief summary of this link here. 

A monoid $M$ is said to be \textit{inverse} if for every $x \in M$ there exists a unique $y \in M$ such that $xyx = x$ and $yxy = y$. This ``inverse'' (which need not be a group inverse) is usually denoted $x^{-1}$. The existence of free inverse monoids was first proved in \tc{1961} by V. V. Wagner \cite{Wagner1952}, although the word problem for free inverse monoids was not proved to be decidable until the \tc{1974} groundbreaking paper by D. Munn \cite{Munn1974}. It bears remarking that free inverse monoids on a non-empty finite set of generators are not finitely presented\footnote{In fact, as was recently discovered, they are not even of homological finiteness type $\FP_2$, see \cite{Gray2020b}.} as ordinary monoids \cite{Schein1975}, even though they are (of course!) finitely presented as inverse monoids. Every inverse monoid admits an inverse monoid presentation. Such presentations are commonly written as $\pres{Inv}{A}{R}$, where, as for group presentations, one considers words over $A \cup A^{-1}$. Of course, one defines \textit{special} inverse monoids as is done for ordinary monoids (see \S\ref{Subsec: Special monoids}). Arguably the most useful tool for studying inverse monoid presentations comes from \textit{Stephen's procedure}, defined in J. B. Stephen's Ph.D. thesis \cite{Stephen1987}\footnote{Stephen also defined an ``ordinary'' monoid analogue of his procedure in his thesis, which is often overlooked. The author of the present thesis, however, has made use of this procedure to characterise the geometry of special monoids \cite{NybergBrodda2020a}.
}. This is a graphical procedure defined entirely analogously to M. Dehn's \textit{Gruppenbild} (see Dehn \cite{Dehn1910} and Chandler \& Magnus \cite[pp. 24--25]{Chandler1982}). We refer the reader to subsequent publications by Stephen \cite{Stephen1990} for more details. 

Stephen's procedure is often useful for solving the word problem in finitely presented inverse monoids, and has been used to solve the word problem in a number of cases (see e.g. \cite{Margolis1993, Birget1994, Hermiller2010}). In the case of special monoids, the main interest in solving this problem comes from the following fascinating link proved by Ivanov, Margolis \& Meakin \cite{Ivanov2001}. A word is \textit{reduced} if it does not contain a subword of the form $aa^{-1}$ or $a^{-1}a$, where $a$ is some letter.

\begin{theorem}[Ivanov, Margolis \& Meakin, \tc{2001}]\label{Thm: Inv reduced => OR WP}
If the word problem is decidable for all inverse monoids of the form $\pres{Inv}{A}{w=1}$ where $w$ is some reduced word, then the word problem is also decidable for every one-relation monoid.
\end{theorem}

We make an important remark that the proof of this theorem as given by Ivanov, Margolis \& Meakin is incomplete. This fact does not appear to have been observed in the literature before. Their proof begins by reducing the word problem for $M$ to the word problem for the right cancellative case $\pres{Mon}{a,b}{aub = ava}$, before embedding this monoid in the special one-relator inverse monoid $\pres{Inv}{a,b}{aub(ava)^{-1} = 1}$, from which the result follows. However, this does not account for the fact that the word problem remains open for the monadic case $\pres{Mon}{a,b}{aub = a}$ (as seen in \S\ref{Subsec: An incorrect proof}). Fortunately, the only property of $\pres{Mon}{a,b}{aub = ava}$ used to produce such an embedding as above is that it is right cancellative; in particular, by adding the case $\pres{Mon}{a,b}{aub=a}$ and embedding this into $\pres{Mon}{a,b}{auba^{-1}=1}$, one finds that the proof of the above theorem is fixed.

As the word problem is decidable for all one-relator groups $\pres{Gp}{A}{w=1}$ and all special one-relation monoids $\pres{Mon}{A}{w=1}$ (see \S\ref{Subsec: Special monoids}), one might view the above result with optimism by conjecturing that the word problem is decidable for all one-relator special inverse monoids $\pres{Inv}{A}{w=1}$. However, the following recent and astounding result by Gray \cite{Gray2020} demonstrates that this optimism is unfounded.

\begin{theorem}[Gray, \tc{2020}]
There exists a special one-relator inverse monoid 
\[
I_B = \pres{Inv}{A}{w_B = 1}
\]
such that the word problem for $I_B$ is undecidable.
\end{theorem}

However, the word $w_B$ constructed by Gray is not reduced, so the implication in Theorem~\ref{Thm: Inv reduced => OR WP} remains a valid path to solving the word problem for all one-relation monoids. Of course, the word problem for special one-relator inverse monoids defined by a reduced word could potentially be significantly harder than that for one-relation monoids. However, the author of the present survey has reasons to suspect that the word problem for all one-relation monoids is equivalent to the word problem for all types of special one-relator inverse monoids occurring in the proof of Theorem~\ref{Thm: Inv reduced => OR WP}.

The author notes that it is very straightforward to show that any special one-relator inverse monoid occurring in the proof of Theorem~\ref{Thm: Inv reduced => OR WP} will have trivial group of units. If special inverse monoids behaved anything like ordinary special monoids (see \S~\ref{Subsec: Special monoids}), we might expect this to be strong evidence in favour of decidability. However, very recently, Gray \& Ru\v{s}kuc \cite{Gray2021} have demonstrated that the group of units of even one-relator special inverse monoids can exhibit rather exotic behaviour when compared to the monoid itself; it need not, for example, be a one-relator group. There is, at present, not even an algorithm known for decomposing the relator word $w$ in $\pres{Inv}{A}{w=1}$ into minimal invertible pieces. Adian's overlap algorithm for ordinary special monoids (see \S\ref{Subsec: Special monoids}) fails spectacularly here, as is demonstrated by the \textit{O'Hare monoid} with presentation
\[
\pres{Inv}{a,b,c,d}{(abcd)(acd)(ad)(abbcd)(acd) = 1},
\]
where the defining relation word has no self-overlaps, but the factorisation into minimal invertible pieces is indicated by the parentheses \cite{Margolis1987}. The author of the present survey has recently found a smaller counterexample, namely
\[
\pres{Inv}{a,b}{aabbaabab = 1},
\]
which is a group (the trefoil knot group), despite the fact that the defining relation word has no self-overlaps. Gray \& Ru\v{s}kuc \cite{Gray2021} propose an improved algorithm (the ``Benois algorithm'') for computing the minimal invertible pieces of a special one-relator inverse monoid, and which correctly computes the pieces of the above two examples. However, the author of the present survey has recently found an example showing that this algorithm does not always produce the minimal invertible pieces. This will appear in future work by the author. Thus, until the word problem for one-relator special inverse monoids is better understood, this seems a difficult avenue for tackling the word problem for one-relation monoids. 

\subsection{The monadic case}\label{Subsec: monadic case}

As mentioned, cf. \S\ref{Subsec: An incorrect proof}, the word problem for monadic one-relation monoids $\pres{Mon}{a,b}{bua=a}$ remains an open problem, conditional on the decidability of the divisibility problems for cycle-free monoids. However, two major results have since appeared for the monadic case, and both were proved in \tc{1997} by Guba. 

The first concerns the equivalence of decision problems for such monoids. While it is clear that decidability of the left divisibility problem (even by a single letter) implies decidability of the word problem, it is not \textit{a priori} true that the converse holds. Furthermore, the r\^ole of the right divisibility problems seems unclear at first. However -- surprisingly -- Guba showed that these problems are all, in fact, equivalent in the monadic case. The corresponding statement for general left cycle-free one-relation monoids is not known to hold or not.

\begin{theorem}[Guba, \tc{1997}]
Let $M = \pres{Mon}{a,b}{bua=a}$. The following are equivalent:
\begin{enumerate}
\item The word problem for $M$ is decidable;
\item The left divisibility problem for $M$ is decidable;
\item The right divisibility problem for $M$ is decidable.
\end{enumerate}
Furthermore, each of these problems is equivalent to its restricted variant of considering equality with (resp. left/right divisibility by) the single letter $a$.
\end{theorem}

The proof uses diagrammatic methods, and are rather involved; we refer the reader to \cite[Theorem~4.1]{Guba1997} to begin navigating the proof, and do not expound on it any further. 

We now present Guba's second major result. Consider a monadic one-relation monoid $\pres{Mon}{a,b}{bua=a}$. Oganesian considered the submonoid $S_M$ generated by all suffixes of the word $bua$ in $\pres{Mon}{a,b}{bua = a}$. He then reduces, by a very general result, the left divisibility problem for $M$ to the left divisibility problem for $S_M$ (\cite[Theorem~1]{Oganesian1982}). He then proves the quite remarkable (and non-trivial!) fact that $S_M$ can be defined by a cycle-free presentation (\cite[Theorem~2]{Oganesian1982}). This makes the implication regarding \tc{sh} and the word problem for monadic one-relation monoids clear (see \S\ref{Subsec: An incorrect proof}). 

Guba, however, studied $S_M$ in more depth. In general, $S_M$ can be shown to not always be a one-relation monoid. However, as it is cycle-free, it is of course group-embeddable (see \S\ref{Subsec: Cancellative}), and so one might reasonably ask questions about the group $G_M$ with the same defining relations. Guba, remarkably, shows that $G_M$ is always isomorphic with a one-relator group (which he denotes $\overline{G}(\Pi)$). Furthermore, this one-relator group is a \textit{positive} one-relator group; recall that a one-relator group $\pres{Gp}{A}{w=1}$ is \textit{positive} if no inverse symbols appear in the word $w$ (these have been studied by Baumslag \cite{Baumslag1971}). He then shows that the left divisibility problem reduces to deciding membership in $S_M$ inside $G_M$. This provides the following astounding statement, which is \cite[Corollary~2.1]{Guba1997}. 

\begin{theorem}[Guba, \tc{1997}]
The word problem for $\pres{Mon}{a,b}{bua=a}$ reduces to the membership problem for a certain submonoid of a positive one-relator group.
\end{theorem}

Perrin \& Schupp \cite{Perrin1984} have proved that a one-relator group is a special one-relator monoid if and only if it admits a presentation $\pres{Gp}{A}{w=1}$ where $w$ is a positive word. Hence, we have the following remarkable statement: if the submonoid membership problem is decidable for all special one-relator monoids $\pres{Mon}{A}{w=1}$, then the word problem is decidable for all one-relation monoids of the form $\pres{Mon}{a,b}{bua = a}$. We do not know of a special one-relator monoid with undecidable submonoid membership problem. However, there is some indication that it might be undecidable, which we discuss below.

\begin{theorem}[{\cite[Theorem~B]{Gray2020}}]
The cycle-free one-relator group 
\[
B = \pres{Gp}{a,b}{abba=baab}
\]
has undecidable submonoid membership problem. 
\end{theorem}

The proof is almost immediate, and so we reproduce it here: consider the right-angled Artin group $A(P_4)$ with defining graph the path on $4$ vertices. Then $A(P_4)$ has generators $a_1, \dots, a_4$ and defining relations $[a_i, a_j] = 1$ whenever $j = i + 1$. It is well-known, and follows quickly from results of Aalbersberg \& Hoogeboom \cite{Aalbersberg1986}, that the submonoid membership problem for $A(P_4)$ is undecidable. Consider the HNN-extension $B$ of $A(P_4)$ with stable letter $t$ and associated subgroups $\langle a_1, a_2, a_3 \rangle \cong \langle a_2, a_3, a_4 \rangle$. Then by a few quick Tietze transformations, one finds
\[
B \cong \pres{Gp}{a, b}{[a, bab^{-1}] = 1} \cong \pres{Gp}{a, b}{abba = baab}.
\]
Hence $A(P_4)$ embeds in $B$, and the result follows.\footnote{The fact that the one-relator group $B$ is cycle-free is not observed directly in Gray \cite{Gray2020}, where only the former of the two above presentations for $B$ is presented, but the latter presentation was shown to the author via private communication with Gray; in fact one quickly obtains the latter from the former by the free group automorphism induced by $a \mapsto ab$ and $b \mapsto b$.} However, we note that by \cite[Corollary~2.6]{Margolis2005} the membership problem for all \textit{positively} generated submonoids of $B$ (with respect to the latter of the two presentations) is decidable. Here a submonoid is \textit{positively} generated if it admits a generating set with only positive words. Furthermore, the \textit{subgroup} membership problem is decidable, as $B$ can easily be shown to be an HNN-extension of $\mathbb{Z}^2$ conjugating one generator to the other; thus this problem is decidable by using the results in e.g. Kapovich et al \cite{Kapovich2005}. We also note that $\pres{Mon}{a,b}{abba=baab}$ trivially has solvable submonoid membership problem.

\subsection{The Collatz conjecture}

In recent years, an interesting connection due to Guba between the word problem for monadic one-relation monoids $\pres{Mon}{a,b}{bua=a}$ and the \textit{Collatz conjecture} has appeared. The author is not aware of any place in the literature where this connection has been written down, and so it is fully expanded on here.

The Collatz conjecture (or the $3x+1$ problem) is a famous problem concerning the function $f \colon \mathbb{N} \to \mathbb{N}$ defined as 
\[
f(x)=
\begin{cases}
\frac{x}{2} & \text{if } x \equiv 0 \mod 2\\
3x+1 & \text{if } x \equiv 1 \mod 2.
\end{cases}
\]
Let $f^{(i)}(x)$ denote the result of applying $f$ $i$ times to $x$. The conjecture states that the sequence $(x, f(x), f^{(2)}(x), \dots)$ eventually reaches $1$, at which point it cycles as $(1, 4, 2, 1, \dots)$. The conjecture has been verified for very large $x$. Occasionally, the time taken to reach $1$ is very long, and has a large degree of unpredictability; for example, even starting with something as small as $x=27$, one has the sequence
\[
(27, 82, 41, \dots, 3077, 9232, 4616, \dots, 5, 16, 8, 4, 2, 1)
\]
taking $111$ steps to reach $1$. We refer to the survey by Lagarias \cite{Lagarias1985} for an excellent overview.

The connection between the Collatz conjecture (and its generalisations) and decision problems has been studied in the past. J. H. Conway \cite{Conway1972} proved that certain generalisations of the Collatz conjecture are undecidable (this has subsequently been strengthened by Kurtz \& Simon \cite{Kurtz2007}). We say that a function $g \colon \mathbb{N} \to \mathbb{N}$ is a \textit{Collatz function} if there is some integer $m$ together with some non-negative rational numbers $a_i, b_i$ $(i < m)$ such that 
\[
g(x)=
\begin{cases}
a_0 x + b_0 & \text{if } x \equiv 0 \mod m,\\
a_1 x + b_1 & \text{if } x \equiv 1 \mod m,\\
\quad\vdots \\
a_{m-1} x + b_{m-1} & \text{if } x \equiv -1 \mod m.\\
\end{cases}
\]
The usual Collatz function is of the above form for $m=2, a_0 = \frac{1}{2}, b_0 = 0, a_1 = 3, b_1 = 1$. 

\begin{theorem}[Conway, \tc{1972}]
There exists a fixed Collatz function $g_e \colon \mathbb{N} \to \mathbb{N}$ such that it is undecidable (with input $x \in \mathbb{N}$) whether $g_e^{(i)}(x) = 1$ for some $i \geq 1$.
\end{theorem}

Thus the iterative behaviour of Collatz functions is enough to encode undecidability statements. See Margenstern's survey \cite{Margenstern2000} for further details on connections between computability and the Collatz conjecture.

Guba realised that there is a connection between Collatz-like functions and the word problem for monadic one-relation monoids. We shall consider the right cancellative case $\pres{Mon}{a,b}{aub=a}$, rather than the usual left cancellative case, as this makes the formulation of the problem easier. We shall not detail how (the non-general form of) $\mathfrak{A}$ works in this case; it is entirely analogous to the left cancellative case, where \textit{prefix} is replaced by \textit{suffix}, etc. 

Let $M = \pres{Mon}{a,b}{aub = a}$. Let us say we have a pair of words $(X, Y)$, and we wish to decide whether $X = Y$ in $M$. We shall describe an iterative procedure on such pairs using the way Adian's algorithm $\mathfrak{A}$ operates. We shall indicate this action by $\longrightarrow^{\mathfrak{A}}$. We first describe the base cases. If one of $X$ and $Y$ is empty, then we terminate, and conclude that $X = Y$ in $M$ if and only if $X \equiv Y \equiv \varepsilon$. If both $X$ and $Y$ are non-empty, then we terminate if $X \equiv Y$, and conclude $X = Y$ in $M$. Otherwise, if $X \not\equiv Y$, and if $X$ (resp. $Y$) is a single letter which is a suffix of $Y$ (resp. $X$), then we terminate and conclude that $X \neq Y$ in $M$. 

Now, if $X \equiv X'a$ and $Y \equiv Y'a$, or $X \equiv X'b$ and $Y \equiv Y'b$, then we cancel these letters. This defines transformations
\begin{align*}
(X'a, Y'a) &\longrightarrow^{\mathfrak{A}} (X', Y') \\
(X'b, Y'b) &\longrightarrow^{\mathfrak{A}} (X', Y')
\end{align*}
If instead $X$ and $Y$ end in different letters, then we flip the pair $(X, Y)$, if necessary, such that it is a pair of the form $(X'b, Y'a)$. We then (as $\mathfrak{A}$ tells us to) replace the rightmost $a$ by $aub$, and cancel the right-most $b$, resulting in a transformation
\begin{align*}
(X'b, Y'a) \longrightarrow^{\mathfrak{A}} (X', Y'au).
\end{align*}
We now iterate this process, which completes the description. We conclude by Adian's theorem regarding $\mathfrak{A}$ (cf. Theorem~\ref{Thm: Algo A three cases}) that $X = Y$ in $M$ if and only if the process terminates successfully; hence, to solve the word problem is equivalent to be able to decide if the above procedure terminates on a given input $(X, Y)$. Of course, this can also be used to study the right divisibility problem, but these two problems are equivalent for $M$ by Guba's earlier result.

The insight by Guba is that one can consider the binary representation of words via $a \mapsto 1$ and $b \mapsto 0$, and that the above procedure thus produces a Collatz-like function $\mathbb{N}\times \mathbb{N} \to \mathbb{N}$.

\begin{example}
Let $M = \pres{Mon}{a,b}{abaab = a}$. Let $X \equiv aabaab$ and $Y \equiv a$. This gives the sequence
\[
(aabaab, a) \to^\mathfrak{A} (aabaa, abaa) \to^\mathfrak{A} (aaba, aba) \to^\mathfrak{A} (aab, ab) \to^\mathfrak{A} (aa, a).
\]
We terminate unsuccessfully, and conclude that $X \neq Y$ in $M$. Using the dyadic representation $a \mapsto 1$ and $b \mapsto 0$, we have that the above sequence is a sequence of transformations
\[
(110110_2, 1_2) \to^\mathfrak{A} (11011_2, 1011_2) \to^\mathfrak{A} (1101_2, 101_2) \to^\mathfrak{A} (110_2, 10_2) \to^\mathfrak{A} (11_2, 1_2)
\] 
Considered as a sequence of natural numbers, the above sequence is 
\[
(54, 1) \to^\mathfrak{A} (27, 11) \to^\mathfrak{A} (13, 5) \to^\mathfrak{A} (6, 2) \to^\mathfrak{A} (3, 1).
\]
In fact, it can be shown that no input word will give an infinite loop. 
\end{example}

The process induced by $\to^\mathfrak{A}$ is not hard to see to be Collatz-like. Note that the cancelling of final letters corresponds to removing the final bit of a binary digit, which is the same as division by $2$ and rounding down, which we hence carry out if the two words map to binary digits that are congruent $\operatorname{mod} 2$. Similarly, the transformation 
\[
(X'b, Y'a) \longrightarrow^{\mathfrak{A}} (X', Y'au)
\]
is carried out when the last bits differ, and when considered in its dyadic representation the transformation corresponds to removing the final bit of the dyadic representation of $X'b$ resp. multiplying the dyadic representation of $Y'a$ by $2^{|u|}$ and adding the binary number corresponding to $u$. Finally, we map $(X, Y) \to^\mathfrak{A} (Y, X)$ when the last bits differ and the last bit of the word corresponding to $X$ is $0$. 

Hence, given a one-relation monoid $\pres{Mon}{a,b}{aub=a}$, let $K$ be the natural number corresponding to the dyadic representation of $u$. Then the word problem for $M$ is decidable if and only if the termination problem is decidable for $G \colon \mathbb{N} \times \mathbb{N} \to \mathbb{N}$ defined by
\[
G(x, y)=
\begin{cases}
\left(\lfloor \frac{x}{2}\rfloor, \lfloor \frac{y}{2}\rfloor\right) & \text{if } x \equiv y \mod 2,\\
\left(\frac{x}{2}, 2^{|u|}y + K\right) & \text{if } x \not\equiv y \text{ and } x \equiv 0 \mod 2, \\
(y, x) & \text{if } x \not\equiv y \text{ and } x \equiv 1 \mod 2.\\
\end{cases}
\]
That is, we can solve the word problem for $M$ if and only if we can decide, for arbitrary input $(x,y)$, whether or not the sequence $((x, y), G(x, y), G^{(2)}(x, y), \dots)$ eventually terminates.

\begin{example}
Continuing the example $M_1 = \pres{Mon}{a,b}{abaab = a}$ from earlier, we find that as $baa \mapsto 011_2$, we have $K = 3$, so
\[
G_1(x, y)=
\begin{cases}
\left(\lfloor \frac{x}{2}\rfloor, \lfloor \frac{y}{2}\rfloor\right) & \text{if } x \equiv y \mod 2,\\
\left(\frac{x}{2}, 8y + 3\right) & \text{if } x \not\equiv y \text{ and } x \equiv 0 \mod 2, \\
(y, x) & \text{if } x \not\equiv y \text{ and } x \equiv 1 \mod 2.\\
\end{cases}
\]
Similarly, in the (right cancellative analogue of the) monoid from Example~\ref{Ex: Looping A}, which is given by $M_2 = \pres{Mon}{a,b}{aabbaab = a}$, we find that $u \equiv abbaa \mapsto 10011_2$, so $K = 19$. Thus we can solve the word problem in this monoid if and only if we can decide when the function 
\[
G_2(x, y)=
\begin{cases}
\left(\lfloor \frac{x}{2}\rfloor, \lfloor \frac{y}{2}\rfloor\right) & \text{if } x \equiv y \mod 2,\\
\left(\frac{x}{2}, 32y + 19 \right) & \text{if } x \not\equiv y \text{ and } x \equiv 0 \mod 2, \\
(y, x) & \text{if } x \not\equiv y \text{ and } x \equiv 1 \mod 2\\
\end{cases}
\]
terminates. Note that this function loops on input $(aaabb, a)$, as indicated by the (reversal) of the transitions given in this example; this corresponds to the infinite loop starting in $(28, 1)$, given by 
\[
(28, 1) \to^\mathfrak{A} (14, 51) \to^\mathfrak{A} (7, 1651) \to^\mathfrak{A} (3, 825) \to^\mathfrak{A} (1, 412) \to^\mathfrak{A} (412, 1).
\]
Now note that $412 \equiv 28 \mod 2^5$, i.e. $1100\textbf{11100}_2 = \textbf{11100}_2 \mod 2^5$, so this process will now loop indefinitely to produce the binary numbers $(1100)^n11100_2$ for $n \geq 0$, which gives the non-terminating sequence 
\[
(28,1) \to^\mathfrak{A} \cdots \to^\mathfrak{A} (412,1) \to^\mathfrak{A} \cdots \to^\mathfrak{A} (6556, 1) \to^\mathfrak{A} \cdots \to^\mathfrak{A} (104860, 1) \to^\mathfrak{A} \cdots 
\]
Hence we conclude that $aaabb \neq a$ in $M$ (we can also conclude that $aaabb$ is not right divisible by $a$). To be clear, if we could detect this looping behaviour for $G_2(x, y)$, we could conclude $aaabb \neq a$ in $M$ simply by checking if $G_2(x, y)$ loops on input $(28,1)$. 
\end{example}

Guba suspects that Conway's undecidability result indicates that it is probable that the word problem or one of its generalisations for monadic one-relation monoids is undecidable. The problem comes down to understanding poor behaviour of the algorithm $\mathfrak{A}$. Of particular interest is the following question.\footnote{The author thanks Victor Guba for explaining the link between this question and Collatz-like functions.}

\begin{question}[Guba, \tc{1997}]
Is there some $M = \pres{Mon}{a,b}{bua = a}$ and a word $w$ such that $w$ is left\footnote{In the English translation, \textit{left} has here been incorrectly (!) translated as \textit{right}; right is wrong, left is right.} divisible by $a^k$ in $M$ for every $k \geq 0$?
\end{question}

If this question has an affirmative answer, then this indicates that the algorithm $\mathfrak{A}$ has quite complicated behaviour. Of course, the right cancellative analogue of the question above, i.e. does there exist some $M = \pres{Mon}{a,b}{aub = a}$ and a word $w$ such that $w$ is right divisible by $a^k$ in $M$ for every $k \geq 0$, can be directly translated to a statement about the above functions by asking for the existence of $m \in \mathbb{N}$ such that the monadic Collatz-like function of $M$ when applied to the initial word $(m, 2^k-1)$ always terminates successfully for $k \geq 1$. Demonstrating the existence of such an $M$ and $m$ would indicate that $\mathfrak{A}$ (and by extension, the Collatz-like function) can behave rather poorly, and, although there is no direct implication of any kind, this would be an important first step if one intends to construct an undecidable Collatz-like function as above.

\clearpage 

\section*{Concluding remarks}

This survey has hopefully given the reader a good feel for the many intricacies, connections to other areas, and beautiful results that combine to make the word problem for one-relation monoids the fascinating problem that it continues to be. Although the many positive decidability results proved over the years seems to indicate that the problem will one day be proved to be decidable, the newly discovered links to undecidable problems by e.g. Gray and Guba might one day come to prove just the opposite.  In any case, it seems fair to say that the problem has turned out to be significantly more difficult than the early positive results in the \tc{1960}s might have suggested. The question of decidability of the word problem for one-relation monoids might be compared to the question of whether or not Thompson's group $F$ is amenable. At the door to the office of my M.Sc. advisor C. Bleak there was printed an introduction to $F$ by Cannon \& Floyd (see \cite{Cannon2011}), which I occasionally would glance at while waiting to be let in. One line stood out: ``at a recent conference devoted to the group a poll was taken. Is $F$ amenable? Twelve participants voted \textit{yes} and twelve voted \textit{no}.'' Perhaps the situation would be similar at a conference devoted to the word problem for one-relation monoids?

Given the overwhelming extent of the contributions by S. I. Adian to this area, it would not be appropriate to end this survey with such an anecdote of my own. Instead, I will mention that at a conference long ago, following a talk on the word problem for one-relation monoids, it is reported that Adian was asked: if a Western mathematician were to solve the problem, would Soviet journals publish the solution in English? Adian's response was: ``before the problem is solved, everybody will publish in English in the Soviet Union''. As Almeida \& Perrin remark, he seems to have been right \cite{Almeida2009}. There is likely no more fitting final sentence for this survey than that given by Adian in \tc{2018} at a conference at the Euler International Mathematical Institute in Saint Petersburg, regarding the word problem for one-relation monoids: \textit{its solution is a task for the future generations of mathematicians.}

\clearpage

{
\bibliography{OneRelationSurveyFinal} 

\begin{thebibliography}{100}

\bibitem{OeisBifix}
The {O}n-{L}ine {E}ncyclopedia of {I}nteger {S}equences.
\newblock \href{https://oeis.org/A094536}{A094536}.
\newblock Number of binary words of length $n$ that are not "bifix-free".

\bibitem{Aalbersberg1986}
Ijsbrand~Jan Aalbersberg and Emo Welzl.
\newblock Trace languages defined by regular string languages.
\newblock {\em RAIRO Inform. Th\'{e}or. Appl.}, 20(2):103--119, 1986.

\bibitem{Adian1955}
S.~I. Adian.
\newblock On the problem of divisibility in semigroups.
\newblock {\em Dokl. Akad. Nauk SSSR (N.S.)}, 103:747--750, 1955.

\bibitem{Adian1957}
S.~I. Adian.
\newblock The role of the cancellation law in presenting cancellation
  semi-groups by means of defining relations.
\newblock {\em Dokl. Akad. Nauk SSSR (N.S.)}, 113:1191--1194, 1957.

\bibitem{Adian1960b}
S.~I. Adian.
\newblock On the embeddability of semigroups in groups.
\newblock {\em Soviet Math. Dokl.}, 1:819--821, 1960.
\newblock [Dokl. Akad. Nauk SSSR {\bf 133}, 2 (1960)].

\bibitem{Adian1960}
S.~I. Adian.
\newblock The problem of identity in associative systems of a special form.
\newblock {\em Soviet Math. Dokl.}, 1:1360--1363, 1960.
\newblock [Dokl. Akad. Nauk SSSR {\bf 135}, 6 (1960)].

\bibitem{Adian1966}
S.~I. Adian.
\newblock {\em Defining relations and algorithmic problems for groups and
  semigroups}.
\newblock Proceedings of the Steklov Institute of Mathematics, No. 85 (1966).
  American Mathematical Society, Providence, R.I., 1966.
\newblock Translated from the Russian by M. Greendlinger.

\bibitem{Adian1976}
S.~I. Adian.
\newblock Word transformations in a semigroup that is given by a system of
  defining relations.
\newblock {\em Algebra and Logic}, 15:379--386, 1976.
\newblock [Algebra i Logika {\bf 15}, 6 (1976)].

\bibitem{Adian1993}
S.~I. Adian.
\newblock On some algorithmic problems for groups and monoids.
\newblock In {\em Rewriting techniques and applications ({M}ontreal, {PQ},
  1993)}, volume 690 of {\em Lecture Notes in Comput. Sci.}, pages 289--300.
  Springer, Berlin, 1993.

\bibitem{Adian1994}
S.~I. Adian.
\newblock On the divisibility problem for monoids defined by one relation.
\newblock {\em Math. Notes}, 55(1):3--7, 1994.
\newblock [Mat. Zametki {\bf 55}, 1 (1994)].

\bibitem{Adian2005}
S.~I. Adian.
\newblock Divisibility problem for one relator monoids.
\newblock {\em Theoret. Comput. Sci.}, 339(1):3--6, 2005.

\bibitem{Adian2018}
S.~I. Adian.
\newblock Gennadi\u{\i} {S}em\"{e}novich {M}akanin's investigations on
  algorithmic questions of group and semigroup theory.
\newblock {\em Russian Math. Surveys}, 73(3):553--568, 2018.
\newblock [Uspekhi Mat. Nauk {\bf 73}, 1 (2018)].

\bibitem{Adian2000}
S.~I. Adian and V.~G. Durnev.
\newblock Algorithmic problems for groups and semigroups.
\newblock {\em Russian Math. Surveys}, 55(2):207--296, 2000.
\newblock [Uspekhi Mat. Nauk {\bf 55}, 2 (2000)].

\bibitem{Adian1984}
S.~I. Adian and G.~S. Makanin.
\newblock Studies in algorithmic questions of algebra.
\newblock {\em Trudy Mat. Inst. Steklov.}, 168:197--217, 1984.
\newblock Algebra, mathematical logic, number theory, topology.

\bibitem{Adian1958}
S.~I. Adian and P.~S. Novikov.
\newblock The word problem for semigroups with one-sided cancellation.
\newblock {\em Z. Math. Logik Grundlagen Math.}, 4:66--88, 1958.
\newblock [Later translated in Amer. Math. Soc. Transl. (2), 46 (1965),
  193--212].

\bibitem{Adian1978}
S.~I. Adian and G.~U. Oganesian.
\newblock On problems of equality and divisibility in semigroups with a single
  defining relation.
\newblock {\em Math. USSR, Izv.}, 12:207--212, 1978.
\newblock [Izv. Akad. Nauk SSSR Ser. Mat. {\bf 42}, 2 (1978)].

\bibitem{Adian1987}
S.~I. Adian and G.~U. Oganesian.
\newblock On the word and divisibility problems for semigroups with one
  relation.
\newblock {\em Math. Notes}, 41:235--240, 1987.
\newblock [Mat. Zametki {\bf 41}, 3 (1987)].

\bibitem{Almeida2009}
Jorge Almeida and Dominique Perrin.
\newblock G\'{e}rard {L}allement (1935--2006).
\newblock {\em Semigroup Forum}, 78(3):379--383, 2009.

\bibitem{Bauer1984}
G.~Bauer and F.~Otto.
\newblock Finite complete rewriting systems and the complexity of the word
  problem.
\newblock {\em Acta Inform.}, 21(5):521--540, 1984.

\bibitem{Baumslag1971}
Gilbert Baumslag.
\newblock Positive one-relator groups.
\newblock {\em Trans. Amer. Math. Soc.}, 156:165--183, 1971.

\bibitem{Birget1994}
Jean-Camille Birget, Stuart~W. Margolis, and John~C. Meakin.
\newblock The word problem for inverse monoids presented by one idempotent
  relator.
\newblock {\em Theoret. Comput. Sci.}, 123(2):273--289, 1994.

\bibitem{Bokut1987}
L.~A. Bokut.
\newblock Imbedding of rings.
\newblock {\em Russian Math. Surveys}, 42(4):105--138, 1987.
\newblock [Uspekhi Mat. Nauk, {\bf 42}, 4 (1987)].

\bibitem{Book1983}
Ronald~V. Book.
\newblock A note on special {T}hue systems with a single defining relation.
\newblock {\em Math. Systems Theory}, 16(1):57--60, 1983.

\bibitem{Book1987}
Ronald~V. Book.
\newblock Thue systems as rewriting systems.
\newblock volume~3, pages 39--68. 1987.
\newblock Rewriting techniques and applications (Dijon, 1985).

\bibitem{Book1982}
Ronald~V. Book, Matthias Jantzen, and Celia Wrathall.
\newblock Monadic {T}hue systems.
\newblock {\em Theoret. Comput. Sci.}, 19(3):231--251, 1982.

\bibitem{Book1984}
Ronald~V. Book and Craig~C. Squier.
\newblock Almost all one-rule {T}hue systems have decidable word problems.
\newblock {\em Discrete Math.}, 49(3):237--240, 1984.

\bibitem{Boone1958}
William~W. Boone.
\newblock An analysis of {T}uring's ``{T}he word problem in semi-groups with
  cancellation''.
\newblock {\em Ann. of Math. (2)}, 67:195--202, 1958.

\bibitem{Borisov1969}
V.~V. Borisov.
\newblock Simple examples of groups with unsolvable word problem.
\newblock {\em Mat. Zametki}, 6:521--532, 1969.
\newblock [Math. Notes {\bf 6} (1969)].

\bibitem{Bouwsma1993}
J.~B Bouwsma.
\newblock {\em Semigroups presented by a single relation}.
\newblock PhD thesis, Penn. State University, 1993.

\bibitem{Brown1992}
Kenneth~S. Brown.
\newblock The geometry of rewriting systems: a proof of the
  {A}nick-{G}roves-{S}quier theorem.
\newblock In {\em Algorithms and classification in combinatorial group theory
  ({B}erkeley, {CA}, 1989)}, volume~23 of {\em Math. Sci. Res. Inst. Publ.},
  pages 137--163. Springer, New York, 1992.

\bibitem{Bucher1988}
W.~Bucher.
\newblock A note on regular classes in special {T}hue systems.
\newblock {\em Discrete Appl. Math.}, 21(3):199--205, 1988.

\bibitem{Burris1981}
Stanley Burris and H.~P. Sankappanavar.
\newblock {\em A course in universal algebra}, volume~78 of {\em Graduate Texts
  in Mathematics}.
\newblock Springer-Verlag, New York-Berlin, 1981.

\bibitem{Cain2013a}
Alan Cain and Victor Maltcev.
\newblock Monoids mon$\langle a, b \: : \: a^\alpha b^\beta a^\gamma b^\delta =
  b \rangle$ admit finite complete rewriting systems.
\newblock {\em Pre-print}, 2013.
\newblock Available at arXiv:1302.0982.

\bibitem{Cain2013b}
Alan Cain and Victor Maltcev.
\newblock Monoids mon$\langle a, b \: : \: a^\alpha b^\beta a^\gamma b^\delta
  a^\varepsilon b^\varphi = b \rangle$ admit finite complete rewriting systems.
\newblock {\em Pre-print}, 2013.
\newblock Available at arXiv:1302.2819.

\bibitem{Cain2009}
Alan~J. Cain and Victor Maltcev.
\newblock Decision problems for finitely presented and one-relation semigroups
  and monoids.
\newblock {\em Internat. J. Algebra Comput.}, 19(6):747--770, 2009.

\bibitem{Cain2009b}
Alan~J. Cain, Graham Oliver, Nik Ru\v{s}kuc, and Richard~M. Thomas.
\newblock Automatic presentations for semigroups.
\newblock {\em Inform. and Comput.}, 207(11):1156--1168, 2009.

\bibitem{Cannon2011}
J.~W. Cannon and W.~J. Floyd.
\newblock What is {$\ldots$} {T}hompson's group?
\newblock {\em Notices Amer. Math. Soc.}, 58(8):1112--1113, 2011.

\bibitem{CasalsRuiz2011}
Montserrat Casals-Ruiz and Ilya Kazachkov.
\newblock On systems of equations over free partially commutative groups.
\newblock {\em Mem. Amer. Math. Soc.}, 212(999):viii+153, 2011.

\bibitem{Chandler1982}
Bruce Chandler and Wilhelm Magnus.
\newblock {\em The history of combinatorial group theory}, volume~9 of {\em
  Studies in the History of Mathematics and Physical Sciences}.
\newblock Springer-Verlag, New York, 1982.

\bibitem{Choffrut2010}
Christian Choffrut and Robert Merca\c{s}.
\newblock Contextual partial commutations.
\newblock {\em Discrete Math. Theor. Comput. Sci.}, 12(4):59--72, 2010.

\bibitem{Chouraqui2009}
Fabienne Chouraqui.
\newblock Rewriting systems and embedding of monoids in groups.
\newblock {\em Groups Complex. Cryptol.}, 1(1):131--140, 2009.

\bibitem{Conway1972}
J.~H. Conway.
\newblock Unpredictable iterations.
\newblock In {\em Proceedings of the {N}umber {T}heory {C}onference ({U}niv.
  {C}olorado, {B}oulder, {C}olo., 1972)}, pages 49--52, 1972.

\bibitem{Cremanns1994b}
Robert Cremanns.
\newblock One-relator groups have finite derivation type.
\newblock {\em Math. Schriften Kassel}, 4(94), 1994.

\bibitem{Cremanns1994}
Robert Cremanns and Friedrich Otto.
\newblock Finite derivation type implies the homological finiteness condition
  {${\rm FP}_3$}.
\newblock {\em J. Symbolic Comput.}, 18(2):91--112, 1994.

\bibitem{Crisp2009}
John Crisp, Eddy Godelle, and Bert Wiest.
\newblock The conjugacy problem in subgroups of right-angled {A}rtin groups.
\newblock {\em J. Topol.}, 2(3):442--460, 2009.

\bibitem{Crvenkovic1995}
Sini\v{s}a Crvenkovi\'{c}.
\newblock Word problems for varieties of algebras (a survey).
\newblock {\em Filomat}, (9, part 3):427--448, 1995.
\newblock Algebra, logic \& discrete mathematics (Ni\v{s}, 1995).

\bibitem{Cummings2014}
P.~A. Cummings and D.~A. Jackson.
\newblock A solvable conjugacy problem for finitely presented {$C(3)$}
  semigroups.
\newblock {\em Semigroup Forum}, 88(1):52--66, 2014.

\bibitem{Dehn1910}
M.~Dehn.
\newblock \"{U}ber die {T}opologie des dreidimensionalen {R}aumes.
\newblock {\em Math. Ann.}, 69(1):137--168, 1910.

\bibitem{Dershowitz1989}
Nachum Dershowitz.
\newblock Completion and its applications.
\newblock In {\em Resolution of equations in algebraic structures, {V}ol. 2},
  pages 31--85. Academic Press, Boston, MA, 1989.

\bibitem{Dyck1882}
Walther Dyck.
\newblock Gruppentheoretische {S}tudien.
\newblock {\em Math. Ann.}, 20(1):1--44, 1882.

\bibitem{Fischer1972}
J.~Fischer, A.~Karrass, and D.~Solitar.
\newblock On one-relator groups having elements of finite order.
\newblock {\em Proc. Amer. Math. Soc.}, 33:297--301, 1972.

\bibitem{Garreta2021}
Albert Garreta and Robert~D. Gray.
\newblock On equations and first-order theory of one-relator monoids.
\newblock {\em Information and Computation}, 2021.

\bibitem{Gerasimov1982}
V.~N. Gerasimov.
\newblock Localization in associative rings.
\newblock {\em Sibirsk. Mat. Zh.}, 23(6):36--54, 205, 1982.

\bibitem{Gluskin1972}
L.~M. Glusk{\={\i}}n and B.~M. Schein.
\newblock {\it {T}he theory of operations as the general theory of groups} by
  {A}. {K}. {S}u\v{s}kevi\v{c}. {A}n historical review.
\newblock {\em Semigroup Forum}, 4:367--371, 1972.

\bibitem{Gray2020}
Robert~D. Gray.
\newblock Undecidability of the word problem for one-relator inverse monoids
  via right-angled {A}rtin subgroups of one-relator groups.
\newblock {\em Invent. Math.}, 219(3):987--1008, 2020.

\bibitem{Gray2021}
Robert~D. Gray and Nik Ru\v{s}kuc.
\newblock On groups of units of special and one-relator inverse monoids.
\newblock {\em Pre-print}, 2021.
\newblock Available at arXiv:2103.02995.

\bibitem{Gray2019}
Robert~D. Gray and Benjamin Steinberg.
\newblock A {L}yndon's identity theorem for one-relator monoids.
\newblock {\em Pre-print}, 2019.
\newblock Available at arXiv:1910.09914.

\bibitem{Gray2020b}
Robert~D. Gray and Benjamin Steinberg.
\newblock Free inverse monoids are not $\operatorname{FP}_2$.
\newblock {\em Pre-print}, 2020.
\newblock Available at arXiv:2002.07690.

\bibitem{Greendlinger1960}
Martin Greendlinger.
\newblock Dehn's algorithm for the word problem.
\newblock {\em Comm. Pure Appl. Math.}, 13:67--83, 1960.

\bibitem{Greendlinger1960b}
Martin Greendlinger.
\newblock {\em Dehn's {A}lgorithm for the {W}ord {P}roblem}.
\newblock PhD thesis, Thesis (Ph.D.)--New York University, 1960.

\bibitem{Guba1994}
V.~S. Guba.
\newblock Conditions for the embeddability of semigroups into groups.
\newblock {\em Math. Notes}, 56(1--2):763--769, 1994.
\newblock Mat. Zametki {\bf 56}, 2 (1994).

\bibitem{Guba1997}
V.~S. Guba.
\newblock On a relation between the word problem and the word divisibility
  problem for semigroups with one defining relation.
\newblock {\em Izv. Math.}, 61(6):1137--1169, 1997.
\newblock [Izv. Ross. Akad. Nauk Ser. Mat. {\bf 61}, 6 (1997)].

\bibitem{Guba2020}
V.~S. Guba.
\newblock On the word problem for $1$-relator monoids.
\newblock S. I. Adian Memorial Conference, May 26--27, 2020, Moscow (online),
  2020.

\bibitem{Guba1997b}
Victor Guba and Mark Sapir.
\newblock Diagram groups.
\newblock {\em Mem. Amer. Math. Soc.}, 130(620):viii+117, 1997.

\bibitem{Hermiller2010}
Susan Hermiller, Steven Lindblad, and John Meakin.
\newblock Decision problems for inverse monoids presented by a single sparse
  relator.
\newblock {\em Semigroup Forum}, 81(1):128--144, 2010.

\bibitem{Higgins1992}
Peter~M. Higgins.
\newblock {\em Techniques of semigroup theory}.
\newblock Oxford Science Publications. The Clarendon Press, Oxford University
  Press, New York, 1992.
\newblock With a foreword by G. B. Preston.

\bibitem{Hollings2009}
Christopher Hollings.
\newblock Anton {K}azimirovich {S}uschkewitsch (1889--1961).
\newblock {\em BSHM Bull.}, 24(3):172--179, 2009.

\bibitem{Hollings2014}
Christopher Hollings.
\newblock Embedding semigroups in groups: not as simple as it might seem.
\newblock {\em Arch. Hist. Exact Sci.}, 68(5):641--692, 2014.

\bibitem{Howie1986}
James Howie and Stephen~J. Pride.
\newblock The word problem for one-relator semigroups.
\newblock {\em Math. Proc. Cambridge Philos. Soc.}, 99(1):33--44, 1986.

\bibitem{Inam2017}
Muhammad Inam, John Meakin, and Robert Ruyle.
\newblock A structural property of {A}dian inverse semigroups.
\newblock {\em Semigroup Forum}, 94(1):93--103, 2017.

\bibitem{Ivanov2001}
S.~V. Ivanov, S.~W. Margolis, and J.~C. Meakin.
\newblock On one-relator inverse monoids and one-relator groups.
\newblock {\em J. Pure Appl. Algebra}, 159(1):83--111, 2001.

\bibitem{Jackson1986}
D.~A. Jackson.
\newblock Some one-relator semigroup presentations with solvable word problems.
\newblock {\em Math. Proc. Cambridge Philos. Soc.}, 99(3):433--434, 1986.

\bibitem{Jackson2001}
D.~A. Jackson.
\newblock The membership problem for {$\langle a,b\colon\ bab^2=ab\rangle$}.
\newblock {\em Semigroup Forum}, 63(1):63--70, 2001.

\bibitem{Jackson2002}
D.~A. Jackson.
\newblock Decision and separability problems for {B}aumslag-{S}olitar
  semigroups.
\newblock {\em Internat. J. Algebra Comput.}, 12(1-2):33--49, 2002.
\newblock International Conference on Geometric and Combinatorial Methods in
  Group Theory and Semigroup Theory (Lincoln, NE, 2000).

\bibitem{Jantzen1981}
M.~Jantzen.
\newblock On a special monoid with a single defining relation.
\newblock {\em Theoret. Comput. Sci.}, 16(1):61--73, 1981.

\bibitem{Jantzen1985}
M.~Jantzen.
\newblock A note on a special one-rule semi-{T}hue system.
\newblock {\em Inform. Process. Lett.}, 21(3):135--140, 1985.

\bibitem{Jantzen1988}
M.~Jantzen.
\newblock {\em Confluent string rewriting}, volume~14 of {\em EATCS Monographs
  on Theoretical Computer Science}.
\newblock Springer-Verlag, Berlin, 1988.

\bibitem{Kambites2009a}
Mark Kambites.
\newblock Small overlap monoids. {I}. {T}he word problem.
\newblock {\em J. Algebra}, 321(8):2187--2205, 2009.

\bibitem{Kambites2009b}
Mark Kambites.
\newblock Small overlap monoids. {II}. {A}utomatic structures and normal forms.
\newblock {\em J. Algebra}, 321(8):2302--2316, 2009.

\bibitem{Kambites2011}
Mark Kambites.
\newblock Generic complexity of finitely presented monoids and semigroups.
\newblock {\em Comput. Complexity}, 20(1):21--50, 2011.

\bibitem{Kapovich2005}
Ilya Kapovich, Richard Weidmann, and Alexei Miasnikov.
\newblock Foldings, graphs of groups and the membership problem.
\newblock {\em Internat. J. Algebra Comput.}, 15(1):95--128, 2005.

\bibitem{Kashintsev1978b}
E.~V. Kashintsev.
\newblock An algorithm for the solution of the conjugacy problem for certain
  semigroups.
\newblock In {\em Recursive functions ({R}ussian)}, pages 18--26. Ivanov. Gos.
  Univ., Ivanovo, 1978.

\bibitem{Kashintsev1978}
E.~V. Kashintsev.
\newblock On the word problem for special semigroups.
\newblock {\em Izv. Akad. Nauk SSSR Ser. Mat.}, 42(6):1401--1416, 1440, 1978.

\bibitem{Kashintsev1992}
E.~V. Kashintsev.
\newblock Small cancellation conditions and embeddability of semigroups in
  groups.
\newblock {\em Internat. J. Algebra Comput.}, 2(4):433--441, 1992.

\bibitem{Kashintsev1993}
E.~V. Kashintsev.
\newblock On the satisfiability of the conditions {$C'(\frac 13)$} and {$C(4)$}
  for special homogeneous semigroups with defining words-degrees.
\newblock {\em Mat. Zametki}, 54(3):40--47, 158, 1993.

\bibitem{Kashkarev2013}
Ilya Kashkarev.
\newblock A generalization of {A}djan's theorem on embeddings of semigroups.
\newblock {\em Asian-Eur. J. Math.}, 6(2), 2013.

\bibitem{Kilibarda1997}
Vesna Kilibarda.
\newblock On the algebra of semigroup diagrams.
\newblock {\em Internat. J. Algebra Comput.}, 7(3):313--338, 1997.

\bibitem{Kobayashi1990}
Yuji Kobayashi.
\newblock Complete rewriting systems and homology of monoid algebras.
\newblock {\em J. Pure Appl. Algebra}, 65(3):263--275, 1990.

\bibitem{Kobayashi1998}
Yuji Kobayashi.
\newblock Homotopy reduction systems for monoid presentations: asphericity and
  low-dimensional homology.
\newblock {\em J. Pure Appl. Algebra}, 130(2):159--195, 1998.

\bibitem{Kobayashi2000}
Yuji Kobayashi.
\newblock Finite homotopy bases of one-relator monoids.
\newblock {\em J. Algebra}, 229(2):547--569, 2000.

\bibitem{Kurth1990}
W.~Kurth.
\newblock {\em Termination und {C}onfluenz von {S}emi-{T}hue-{S}ystems mit nur
  einer {R}egel}.
\newblock PhD thesis, Technische Universit\"at Clausthal, 1990.

\bibitem{Kurtz2007}
Stuart~A. Kurtz and Janos Simon.
\newblock The undecidability of the generalized {C}ollatz problem.
\newblock In {\em Theory and applications of models of computation}, volume
  4484 of {\em Lecture Notes in Comput. Sci.}, pages 542--553. Springer,
  Berlin, 2007.

\bibitem{Lagarias1985}
Jeffrey~C. Lagarias.
\newblock The {$3x+1$} problem and its generalizations.
\newblock {\em Amer. Math. Monthly}, 92(1):3--23, 1985.

\bibitem{Lallement1974}
G\'{e}rard Lallement.
\newblock On monoids presented by a single relation.
\newblock {\em J. Algebra}, 32:370--388, 1974.

\bibitem{Lallement1988}
G\'{e}rard Lallement.
\newblock Some algorithms for semigroups and monoids presented by a single
  relation.
\newblock In {\em Semigroups, theory and applications ({O}berwolfach, 1986)},
  volume 1320 of {\em Lecture Notes in Math.}, pages 176--182. Springer,
  Berlin, 1988.

\bibitem{Lallement1993}
G\'{e}rard Lallement.
\newblock The word problem for semigroups presented by one relation.
\newblock In {\em Semigroups ({L}uino, 1992)}, pages 167--173. World Sci.
  Publ., River Edge, NJ, 1993.

\bibitem{Lallement1995}
G\'{e}rard Lallement.
\newblock The word problem for {T}hue rewriting systems.
\newblock In {\em Term rewriting ({F}ont {R}omeux, 1993)}, volume 909 of {\em
  Lecture Notes in Comput. Sci.}, pages 27--38. Springer, Berlin, 1995.

\bibitem{Lallement1994}
G\'{e}rard Lallement and Laurent Rosaz.
\newblock Residual finiteness of a class of semigroups presented by a single
  relation.
\newblock {\em Semigroup Forum}, 48(2):169--179, 1994.

\bibitem{Magnus1932}
W.~Magnus.
\newblock Das {I}dentit\"{a}tsproblem f\"{u}r {G}ruppen mit einer definierenden
  {R}elation.
\newblock {\em Math. Ann.}, 106(1):295--307, 1932.

\bibitem{Magnus1930}
Wilhelm Magnus.
\newblock \"{U}ber diskontinuierliche {G}ruppen mit einer definierenden
  {R}elation. ({D}er {F}reiheitssatz).
\newblock {\em J. Reine Angew. Math.}, 163:141--165, 1930.

\bibitem{Makanin1966b}
G.~S. Makanin.
\newblock {\em On the {I}dentity {P}roblem for {F}initely {P}resented Groups
  and {S}emigroups}.
\newblock PhD thesis, Steklov Mathematical Institute, Moscow, 1966.

\bibitem{Makanin1966}
G.~S. Makanin.
\newblock On the identity problem in finitely defined semigroups.
\newblock {\em Soviet Math. Dokl.}, 7:1478--1480, 1966.
\newblock [Dokl. Akad. Nauk SSSR {\bf 171} (1966)].

\bibitem{Maltsev1937}
A.~I. Maltsev.
\newblock On the immersion of an algebraic ring into a field.
\newblock {\em Math. Ann.}, 113(1):686--691, 1937.

\bibitem{Maltsev1939}
A.~I. Maltsev.
\newblock {O}n the embedding of associative systems into groups ({R}ussian).
\newblock {\em Mat. Sbornik (Receuil Math\'em.)}, 6(48)(2):331--336, 1939.

\bibitem{Maltsev1940}
A.~I. Maltsev.
\newblock {O}n the embedding of associative systems into groups {II}
  ({R}ussian).
\newblock {\em Mat. Sbornik (Receuil Math\'em.)}, 8(50)(2):251--264, 1940.

\bibitem{Maltsev1965}
A.~I. Maltsev.
\newblock {\em Algoritmy i rekursivnye funktsii [{A}lgorithms and Recursive
  Functions]}.
\newblock Izdat. ``Nauka'', Moscow, 1965.
\newblock [Translated into English by Leo F. Boron, with the collaboration of
  Luis E. Sanchis, John Stillwell and Kiyoshi Is\'eki. Wolters-Noordhoff
  Publishing, Groningen (1970) 372 pp.].

\bibitem{Margenstern2000}
Maurice Margenstern.
\newblock Frontier between decidability and undecidability: a survey.
\newblock {\em Theoret. Comput. Sci.}, 231(2):217--251, 2000.
\newblock Universal machines and computations (Metz, 1998).

\bibitem{Margolis2005}
Stuart~W. Margolis, John Meakin, and Zoran {\v{S}}uni\'{k}.
\newblock Distortion functions and the membership problem for submonoids of
  groups and monoids.
\newblock In {\em Geometric methods in group theory}, volume 372 of {\em
  Contemp. Math.}, pages 109--129. Amer. Math. Soc., Providence, RI, 2005.

\bibitem{Margolis1993}
Stuart~W. Margolis and John~C. Meakin.
\newblock Inverse monoids, trees and context-free languages.
\newblock {\em Trans. Amer. Math. Soc.}, 335(1):259--276, 1993.

\bibitem{Margolis1987}
Stuart~W. Margolis, John~C. Meakin, and Joseph~B. Stephen.
\newblock Some decision problems for inverse monoid presentations.
\newblock In {\em Semigroups and their applications ({C}hico, {C}alif., 1986)},
  pages 99--110. Reidel, Dordrecht, 1987.

\bibitem{Markov1947}
A.~A. Markov.
\newblock On the impossibility of certain algorithms in the theory of
  associative systems.
\newblock {\em Doklady Akad. Nauk SSSR (N.S.)}, 55:583--586, 1947.

\bibitem{Markov1947a}
A.~A. Markov.
\newblock On the impossibility of certain algorithms in the theory of
  associative systems. {II}.
\newblock {\em Doklady Akad. Nauk SSSR (N.S.)}, 58:353--356, 1947.

\bibitem{Matijasevic1967}
Ju.~V. Matijasevi\v{c}.
\newblock Simple examples of unsolvable associative calculi.
\newblock {\em Soviet Math. Dokl.}, 8:555--557, 1967.
\newblock [Dokl. Akad. Nauk SSSR {\bf 173} (1967)].

\bibitem{Matijasevic1967b}
Ju.~V. Matijasevi\v{c}.
\newblock Simple examples of unsolvable canonical calculi.
\newblock {\em Trudy Mat. Inst. Steklov}, 93:50--88, 1967.

\bibitem{McCool1973}
James McCool and Paul~E. Schupp.
\newblock On one relator groups and {${\rm HNN}$} extensions.
\newblock {\em J. Austral. Math. Soc.}, 16:249--256, 1973.
\newblock Collection of articles dedicated to the memory of Hanna Neumann, II.

\bibitem{McNaughton1987}
R.~McNaughton and P.~Narendran.
\newblock Special monoids and special {T}hue systems.
\newblock {\em J. Algebra}, 108(1):248--255, 1987.

\bibitem{Metivier1985}
Yves M{\'{e}}tivier.
\newblock Calcul de longueurs de cha\^{\i}nes de r\'{e}\'{e}criture dans le
  mono\"{\i}de libre.
\newblock {\em Theoret. Comput. Sci.}, 35(1):71--87, 1985.

\bibitem{Mitchell2021}
J.~D. Mitchell and M.~Tsalakou.
\newblock An explicit algorithm for normal forms in small overlap monoids.
\newblock {\em Pre-print}, 2021.
\newblock Available at arXiv:2105.12125.

\bibitem{Munn1974}
W.~D. Munn.
\newblock Free inverse semigroups.
\newblock {\em Proc. London Math. Soc. (3)}, 29:385--404, 1974.

\bibitem{Narendran1991}
Paliath Narendran, Colm {\'{O}}'D{\'{u}}nlaing, and Friedrich Otto.
\newblock It is undecidable whether a finite special string-rewriting system
  presents a group.
\newblock {\em Discrete Math.}, 98(2):153--159, 1991.

\bibitem{Newman1968}
B.~B. Newman.
\newblock {\em Some aspects of one-relator groups}.
\newblock PhD thesis, University College of Townsville (later James Cook
  University), 1968.

\bibitem{Newman1942}
M.~H.~A. Newman.
\newblock On theories with a combinatorial definition of ``equivalence.''.
\newblock {\em Ann. of Math. (2)}, 43:223--243, 1942.

\bibitem{Nielsen1973}
P.~Tolstrup Nielsen.
\newblock A note on bifix-free sequences.
\newblock {\em IEEE Trans. Inform. Theory}, IT-19:704--706, 1973.

\bibitem{NybergBrodda2020a}
C.-F. Nyberg-Brodda.
\newblock The geometry of special monoids.
\newblock {\em Pre-print (submitted)}, 2020.
\newblock Available at arXiv:2011.04536.

\bibitem{NybergBrodda2020c}
C.-F. Nyberg-Brodda.
\newblock On the word problem for compressible monoids.
\newblock {\em Pre-print (submitted)}, 2020.
\newblock Available at arXiv:2012.01402.

\bibitem{NybergBrodda2020b}
C.-F. Nyberg-Brodda.
\newblock On the word problem for special monoids.
\newblock {\em Pre-print (submitted)}, 2020.
\newblock Available at arXiv:2011.09466.

\bibitem{NybergBrodda2021a}
C.-F. Nyberg-Brodda.
\newblock The {B.} {B.} {N}ewman spelling theorem.
\newblock {\em The British Journal for the History of Mathematics}, 36:2, 2021.

\bibitem{NybergBrodda2021c}
C.-F. Nyberg-Brodda.
\newblock On computing groups of units.
\newblock {\em Pre-print}, 2021.
\newblock In preparation.

\bibitem{NybergBrodda2021b}
C.-F. Nyberg-Brodda.
\newblock On the {I}dentity {P}roblem for {F}initely {P}resented groups and
  {S}emigroups, 2021.
\newblock English translation of ``K Probleme Tozhdestva v
  Konechno-opredelennyh Gruppah i Polugruppah'', Ph.D. Thesis by G. S. Makanin
  (1966) (Available online at arXiv:2102.00745).

\bibitem{NybergBrodda2021d}
C.-F. Nyberg-Brodda.
\newblock {\em Theory of generalised groups}.
\newblock 2021.
\newblock English translation of ``The theory of generalised groups'', by A. K.
  Sushkevich (1937). In preparation.

\bibitem{Oganesian1978}
G.~U. Oganesian.
\newblock A class of semigroups with a solvable word problem.
\newblock {\em Math. Notes}, 23(5):640--643, 1978.
\newblock [Mat. Zametki {\bf 24}, 2 (1978)].

\bibitem{Oganesian1978b}
G.~U. Oganesian.
\newblock Problems of equality and divisibility in a semigroup with a defining
  relation of the form {$a=bA$}.
\newblock {\em Math. USSR-Izv.}, 12(3):557--566, 1978.
\newblock [Izv. Akad. Nauk SSSR Ser. Mat. {\bf 42}, 3 (1978)].

\bibitem{Oganesian1979}
G.~U. Oganesian.
\newblock The solvability of the word problem for semigroups with a defining
  relation of the form {$A=BtC$}.
\newblock {\em Izv. Akad. Nauk Armyan. SSR Ser. Mat.}, 14(4):288--291, 315,
  1979.

\bibitem{Oganesian1982}
G.~U. Oganesian.
\newblock Semigroups with one relation and semigroups without cycles.
\newblock {\em Math. USSR, Izv.}, 20:89--95, 1983.
\newblock [Izv. Akad. Nauk SSSR Ser. Mat. {\bf 46} (1982)].

\bibitem{Oganesian1984}
G.~U. Oganesian.
\newblock The isomorphism problem for semigroups with one defining relation.
\newblock {\em Math. Notes}, 35:360--363, 1984.
\newblock [Mat. Zametki {\bf 35}, 5 (1984).

\bibitem{Osipova1968}
V.~A. Osipova.
\newblock The word problem for finitely defined semigroups.
\newblock {\em Sov. Math. Dokl.}, 9:237--240, 1968.
\newblock [Dokl. Akad. Nauk SSSR {\bf 9} (1968)].

\bibitem{Osipova1972}
V.~A. Osipova.
\newblock Equations with one unknown in semigroups with a restricted measure of
  overlap of defining words.
\newblock {\em Sov. Math. Dokl.}, 13:542--545, 1972.
\newblock [Dokl. Akad. Nauk SSSR {\bf 203} (1972)].

\bibitem{Osipova1973}
V.~A. Osipova.
\newblock An algorithm for recognizing the solvability of equations with one
  unknown in semigroups with a measure of overlap of the defining words that is
  less that {$1/3$}.
\newblock {\em Mat. Sb. (N.S.)}, 92(134):3--33, 165, 1973.

\bibitem{Otto1984a}
Friedrich Otto.
\newblock Conjugacy in monoids with a special {C}hurch-{R}osser presentation is
  decidable.
\newblock {\em Semigroup Forum}, 29(1-2):223--240, 1984.

\bibitem{Otto1984}
Friedrich Otto.
\newblock Finite complete rewriting systems for the {J}antzen monoid and the
  {G}reendlinger group.
\newblock {\em Theoret. Comput. Sci.}, 32(3):249--260, 1984.

\bibitem{Otto1988}
Friedrich Otto.
\newblock An example of a one-relator group that is not a one-relation monoid.
\newblock {\em Discrete Math.}, 69(1):101--103, 1988.

\bibitem{Otto1992}
Friedrich Otto.
\newblock Completing a finite special string-rewriting system on the congruence
  class of the empty word.
\newblock {\em Appl. Algebra Engrg. Comm. Comput.}, 2(4):257--274, 1992.

\bibitem{Otto1995}
Friedrich Otto.
\newblock Solvability of word equations modulo finite special and confluent
  string-rewriting systems is undecidable in general.
\newblock {\em Inform. Process. Lett.}, 53(5):237--242, 1995.

\bibitem{Otto1991}
Friedrich Otto and Louxin Zhang.
\newblock Decision problems for finite special string-rewriting systems that
  are confluent on some congruence class.
\newblock {\em Acta Inform.}, 28(5):477--510, 1991.

\bibitem{Paris2002}
Luis Paris.
\newblock Artin monoids inject in their groups.
\newblock {\em Comment. Math. Helv.}, 77(3):609--637, 2002.

\bibitem{Pedersen1984}
John Pedersen.
\newblock {\em Confluence {M}ethods and the {W}ord {P}roblem in {U}niversal
  {A}lgebra}.
\newblock PhD thesis, Emory University, Australia, 1984.

\bibitem{Pedersen1989}
John Pedersen.
\newblock Morphocompletion for one-relation monoids.
\newblock In Nachum Dershowitz, editor, {\em Rewriting Techniques and
  Applications, 3rd International Conference, North Carolina, USA, April 3-5,
  1989, Proceedings}, volume 355 of {\em Lecture Notes in Computer Science},
  pages 574--578. Springer, 1989.

\bibitem{Perrin1984}
Dominique Perrin and Paul Schupp.
\newblock Sur les mono\"{\i}des \`a un relateur qui sont des groupes.
\newblock {\em Theoret. Comput. Sci.}, 33(2-3):331--334, 1984.

\bibitem{Post1947}
Emil~L. Post.
\newblock Recursive unsolvability of a problem of {T}hue.
\newblock {\em J. Symbolic Logic}, 12:1--11, 1947.

\bibitem{Potts1984}
D.~H. Potts.
\newblock Remarks on an example of {J}antzen.
\newblock {\em Theoret. Comput. Sci.}, 29(3):277--284, 1984.

\bibitem{Pride1995}
Stephen~J. Pride.
\newblock Low-dimensional homotopy theory for monoids.
\newblock {\em Internat. J. Algebra Comput.}, 5(6):631--649, 1995.

\bibitem{Remmers1980}
John~H. Remmers.
\newblock On the geometry of semigroup presentations.
\newblock {\em Adv. in Math.}, 36(3):283--296, 1980.

\bibitem{Sarkisian1976}
O.~A. Sarkisian.
\newblock The relation between algorithmic problems in groups and semigroups.
\newblock {\em Sov. Math. Dokl.}, 17:615--617, 1976.
\newblock [Dokl. Akad. Nauk SSSR {\bf 227} (1976)].

\bibitem{Sarkisian1979}
O.~A. Sarkisian.
\newblock Some relations between word problems and divisibility problems in
  groups and semigroups.
\newblock {\em Math. USSR}, 15:161--171, 1980.
\newblock [Izv. Akad. Nauk SSSR Ser. Mat. {\bf 43} (1979)].

\bibitem{Sarkisian1981}
O.~A. Sarkisian.
\newblock On the word and divisibility problems in semigroups and groups
  without cycles.
\newblock {\em Math. USSR Izv.}, 19:643--656, 1982.
\newblock [Izv. Akad. Nauk SSSR Ser. Mat. {\bf 45} (1981)].

\bibitem{Schein1975}
B.~M. Schein.
\newblock Free inverse semigroups are not finitely presentable.
\newblock {\em Acta Math. Acad. Sci. Hungar.}, 26:41--52, 1975.

\bibitem{Schein2015}
Boris~M. Schein.
\newblock Transitive representations of inverse semigroups.
\newblock {\em J. Algebra}, 441:108--124, 2015.

\bibitem{Scott1956}
D.~S. Scott.
\newblock A short recursively unsolvable problem.
\newblock {\em J. Symbolic Logic}, 21:111--112, 1956.

\bibitem{Shestakov2005}
S.~L. Shestakov.
\newblock The equation {$[x,y]=g$} in partially commutative groups.
\newblock {\em Siberian Math. J.}, 46(2):364--372, 2005.
\newblock [Sibirsk. Math. Zh. {\bf 46} (2005)].

\bibitem{Shestakov2006}
S.~L. Shestakov.
\newblock The equation {$x^2y^2=g$} in partially commutative groups.
\newblock {\em Siberian Math. J.}, 47(2):383--390, 2006.
\newblock [Sibirsk. Math. Zh. {\bf 47} (2006)].

\bibitem{Shevrin1985}
L.~N. Shevrin and M.~V. Volkov.
\newblock Identities of semigroups.
\newblock {\em Soviet Math. (Iz. VUZ)}, (29):1--64, 1985.
\newblock [Izv. Vyssh. Uchebn. Zaved. Mat. {\bf 11} (1985)].

\bibitem{Shneerson1972}
L.~M. Shneerson.
\newblock Identities in semigroups with one defining relation.
\newblock {\em Logic, Algebra, and Computational Mathematics, Ivanovo
  Pedagogical Institute}, 1(1--2):139--156, 1972.

\bibitem{Shneerson1972a}
L.~M. Shneerson.
\newblock Identities in semigroups with one defining relation. {II}.
\newblock {\em Logic, Algebra, and Computational Mathematics, Ivanovo
  Pedagogical Institute}, 1(3--4):112--124, 1972.

\bibitem{Squier1983}
C.~Squier and C.~Wrathall.
\newblock The {F}reiheitssatz for one-relation monoids.
\newblock {\em Proc. Amer. Math. Soc.}, 89(3):423--424, 1983.

\bibitem{Squier1987}
Craig~C. Squier.
\newblock Word problems and a homological finiteness condition for monoids.
\newblock {\em J. Pure Appl. Algebra}, 49(1-2):201--217, 1987.

\bibitem{Stephen1987}
J.~B. Stephen.
\newblock {\em Applications of automata theory to presentations of monoids and
  inverse monoids}.
\newblock ProQuest LLC, Ann Arbor, MI, 1987.
\newblock Thesis (Ph.D.)--The University of Nebraska - Lincoln.

\bibitem{Stephen1990}
J.~B. Stephen.
\newblock Presentations of inverse monoids.
\newblock {\em J. Pure Appl. Algebra}, 63(1):81--112, 1990.

\bibitem{Sushkevich1922}
A.~K. Sushkevich.
\newblock {\em The Theory of Operations as the General Theory of Groups}.
\newblock PhD thesis, Voronezh State University, 1922.

\bibitem{Sushkevich1935}
A.~K. Sushkevich.
\newblock On the extension of a semigroup to the whole group.
\newblock {\em Zap. Khark. Mat. Obshch.}, 12:81--87, 1935.
\newblock in Ukrainian.

\bibitem{Sushkevich1937}
A.~K. Sushkevich.
\newblock {\em Theory of generalised groups}.
\newblock DNTVU, Kharkiv, Kiev, 1937.
\newblock in Russian.

\bibitem{Thue1914}
Axel Thue.
\newblock Probleme \"uber {V}er\"anderungen von {Z}eichenreihennach gegebenen
  {R}egeln.
\newblock {\em Christiana Videnskabs-Selskabs Skrifter}, 10, 1914.

\bibitem{Tseitin1958}
G.~S. Tseitin.
\newblock An associative calculus with an insoluble problem of equivalence.
\newblock {\em Trudy Mat. Inst. Steklov.}, 52:172--189, 1958.

\bibitem{Turing1950}
A.~M. Turing.
\newblock The word problem in semi-groups with cancellation.
\newblock {\em Ann. of Math. (2)}, 52:491--505, 1950.

\bibitem{Vazhenin1983}
Yu.~M. Vazhenin.
\newblock Semigroups with one defining relation whose elementary theories are
  decidable.
\newblock {\em Sib. Math. J.}, 24:33--41, 1983.
\newblock [Sibirsk. Mat. Zh. {\bf 24} (1983)].

\bibitem{Wagner1952}
V.~V Wagner.
\newblock Generalised groups.
\newblock {\em Dokl. Akad. Nauk SSSR}, 84:1119--1122, 1952.

\bibitem{Watier1996}
Guillaume Watier.
\newblock Left-divisibility and word problems in single relation monoids.
\newblock {\em Semigroup Forum}, 53(2):194--207, 1996.

\bibitem{Watier1997}
Guillaume Watier.
\newblock On the word problem for single relation monoids with an unbordered
  relator.
\newblock {\em Internat. J. Algebra Comput.}, 7(6):749--770, 1997.

\bibitem{Wise2012}
Daniel~T. Wise.
\newblock {\em From riches to raags: 3-manifolds, right-angled {A}rtin groups,
  and cubical geometry}, volume 117 of {\em CBMS Regional Conference Series in
  Mathematics}.
\newblock Published for the Conference Board of the Mathematical Sciences,
  Washington, DC; by the American Mathematical Society, Providence, RI, 2012.

\bibitem{Wrathall1995}
C.~Wrathall and V.~Diekert.
\newblock On confluence of one-rule trace-rewriting systems.
\newblock {\em Math. Systems Theory}, 28(4):341--361, 1995.

\bibitem{Wrathall1992}
C.~Wrathall, V.~Diekert, and F.~Otto.
\newblock One-rule trace-rewriting systems and confluence.
\newblock In {\em Mathematical foundations of computer science 1992 ({P}rague,
  1992)}, volume 629 of {\em Lecture Notes in Comput. Sci.}, pages 511--521.
  Springer, Berlin, 1992.

\bibitem{Yasuhara1970}
Ann Yasuhara.
\newblock The solvability of the word problem for certain semigroups.
\newblock {\em Proc. Amer. Math. Soc.}, 26:645--650, 1970.

\bibitem{Zhang1991}
Louxin Zhang.
\newblock Conjugacy in special monoids.
\newblock {\em J. Algebra}, 143(2):487--497, 1991.

\bibitem{Zhang1992a}
Louxin Zhang.
\newblock Applying rewriting methods to special monoids.
\newblock {\em Math. Proc. Cambridge Philos. Soc.}, 112(3):495--505, 1992.

\bibitem{Zhang1992d}
Louxin Zhang.
\newblock Congruential languages specified by special string-rewriting systems.
\newblock In {\em Words, languages and combinatorics ({K}yoto, 1990)}, pages
  551--563. World Sci. Publ., River Edge, NJ, 1992.

\bibitem{Zhang1992c}
Louxin Zhang.
\newblock On the conjugacy problem for one-relator monoids with elements of
  finite order.
\newblock {\em Internat. J. Algebra Comput.}, 2(2):209--220, 1992.

\bibitem{Zhang1992b}
Louxin Zhang.
\newblock A short proof of a theorem of {A}djan.
\newblock {\em Proc. Amer. Math. Soc.}, 116(1):1--3, 1992.

\bibitem{Zhang1996}
Louxin Zhang, Lian Li, and Jinzhao Wu.
\newblock On the descriptive power of special {T}hue systems.
\newblock {\em Discrete Math.}, 160(1-3):291--297, 1996.

\end{thebibliography}
\bibliographystyle{plain}
}
\end{document}